\documentclass[reqno,11pt]{amsart}
\usepackage{amsmath,amssymb}
\usepackage{geometry}
\usepackage{amsmath, amsfonts,amsthm,amssymb,amsbsy,upref,color,graphicx,amscd,hyperref,enumerate}
\usepackage[active]{srcltx}
\usepackage[latin1]{inputenc}
\usepackage{bbm}
\usepackage{algorithm}
\usepackage{stmaryrd}
\usepackage{array}
\usepackage{algorithmicx}
\usepackage{algpseudocode}
\usepackage{graphicx}
\usepackage{multirow}
\usepackage{wrapfig}
\usepackage[justification=centering]{caption}
\usepackage[position=b]{subcaption}
\graphicspath{{./pics_eps/}{./pics/}}

\numberwithin{equation}{section}
\theoremstyle{plain}
\newtheorem{thm}{Theorem}[section]

\theoremstyle{definition}
\newtheorem{rem}[thm]{Remark}

\newcommand{\gconv}{\stackrel{\Gamma}{\longrightarrow}}

\newcommand{\di}{\operatorname{div}}
\newcommand{\Per}{\operatorname{Per}}
\newcommand{\ds}{\displaystyle}
\renewcommand{\Bbb}{\mathbb}
\newcommand{\bo}[1]{{\bf #1}}
\providecommand{\e}[1]{\ensuremath{\times 10^{#1}}}

\author[B. Bogosel]{Beniamin Bogosel}

\author[\'E. Oudet]{\'Edouard Oudet}

\address[Beniamin Bogosel, Edouard Oudet]{Laboratoire Jean Kuntzmann, Universit\' e Grenoble Alpes, Bâtiment IMAG, 700 avenue centrale, 38400 Saint Martin d'Hères France}

\email[Beniamin Bogosel]{beniamin.bogosel@univ-savoie.fr}
\email[\'Edouard Oudet]{edouard.oudet@imag.fr}

\title{Partitions of minimal length on manifolds}

\begin{document}
\begin{abstract}
We study partitions on three dimensional manifolds which minimize the total geodesic perimeter. We propose a relaxed framework based on a $\Gamma$-convergence result and we show some numerical results. We compare our results to those already present in the literature in the case of the sphere. For general surfaces we provide an optimization algorithm on meshes which can give a good approximation of the optimal cost, starting from the results obtained using the relaxed formulation.
\end{abstract}

\maketitle

\section{Introduction}

In this article we propose a theoretical and numerical framework for the study of the partitions $(\omega_i)_{i=1}^n$ of a surface $S \subset \Bbb{R}^3$ which minimize the total geodesic perimeter while keeping a prescribed area for each cell.
Thus, we are interested in minimizing $\mathcal{H}^1(\cup_{i=1}^n \partial_S \omega_i)$
or equivalently
\[ \Per(\omega_1)+...+\Per(\omega_n)\]
in the class of partitions $(\omega_i)$ of the surface $S$ such that $|\omega_i| = c_i$, with the compatibility constraint $c_1+...+c_n = |S|$. 
Here $\partial_S \omega$ denotes the boundary of a set $\omega$ as a subset of the surface $S$, $\Per(\omega)$ denotes the geodesic perimeter of $\omega$, i.e. the perimeter of $\omega$ regarded as a subset of the surface $S$ and $|\omega|$ is the area of the subset $\omega$. General theoretical results concerning these minimal partitioning problems are presented by Morgan in \cite{morgan-bubbles}. This theoretical result states that boundaries of a minimal-perimeter partition are arcs of constant geodesic curvature and the boundaries of the sets meet in threes with angles of measure $2\pi/3$.

The more specific case concerning the minimal perimeter partitions of sphere with cells of equal areas was intensively studied from both theoretical and numerical points of view. In the case $n=2$ the solution is the partition into two half-spheres. This was proved by Bernstein in 1905 \cite{bernstein-sphere}. In the case $n=3$ the optimal candidate is the partition of the sphere into three slices corresponding to an angle of $2\pi/3$. This was proved by Masters in \cite{masters-sphere}. The case $n = 12$ was solved by Hales in \cite{hales-sphere} using methods similar to the ones involved in the proof of the honeycomb conjecture \cite{hales}. The case $n=4$ was treated by Engelstein in \cite{engelstein-four} and the corresponding optimal partition is the one associated to the regular tetrahedron. 

The case of the sphere has been studied numerically by Cox and Flikkema \cite{cox-partitions} using the Surface Evolver software \cite{evolver}. They perform computations for $n \in \llbracket 2,32\rrbracket$ and they confirm the natural conjecture for $n=6$: the optimal partition in this case is probably the one associated to the cube. Their algorithm performs the perimeter optimization after choosing a topological structure for the partition. Thus, the optimization algorithm has to know \emph{a priori} the topological structure in order to find the corresponding local minimum. In the end we keep the configuration which gives the best optimal cost among the admissible combinatorial possibilities.

The algorithm we propose is a generalization of the ideas in \cite{oudet} to the case of surfaces. First, there is a theoretical result, similar to the theorem of Modica and Mortola, which we present in Section \ref{theoretical-result}. This theoretical result justifies the use of the functional
\[ \mathcal{J}_\varepsilon(u) = \varepsilon\int_S|\nabla_\tau u |^2 +\frac{1}{\varepsilon} \int_S u^2(1-u)^2\]
as an approximation of the perimeter as $\varepsilon \to 0$. The direct consequence of the $\Gamma$-convergence result is that a sequence of minimizers $u_\varepsilon$ for $\mathcal{J}_\varepsilon$ under the constraint $\int_S u_\varepsilon = c$ converges to a minimizer of the geodesic perimeter under area constraint. For the partitioning case we prove that functionals of the type
\[ \sum_{i=1}^n \mathcal J_\varepsilon(u_i)\]
approximate the perimeter as $\varepsilon \to 0$, where $u_i$ are functions associated to the sets $\omega_i$ which satisfy some integral and non-overlapping constraints. We implement an optimization algorithm which is able to solve the above problem on a large class of surfaces. This is an advantage over the methods used in \cite{cox-partitions} which can be used only in the case of the sphere.

Working with the relaxed formulation does not provide an exact representation of the contours. Thus, we cannot directly provide the associated cost once we have the relaxed optimal partitions. The particular case of the sphere can be solved directly by noting that boundaries between two cells have constant geodesic curvature \cite{morgan-bubbles} and are, thus, arcs of circles. We recover all the results presented in \cite{cox-partitions} in the case of the sphere. On more complex surfaces it is complicated to explicitly work with curves of constant geodesic curvature. Nevertheless, we can extract the contours from the density representation in order to compute the total perimeter. Since the extracted contours are not smooth, we perform a constrained optimization stage on the triangulated surface preserving the topology to obtain reliable approximations of the optimal costs.

\section{Theoretical result}
\label{theoretical-result}

As in \cite{oudet} we would like to have a rigorous theoretical framework which justifies our numerical method. In the euclidean case it was an adapted version of the Modica-Mortola theorem to the case of partitions which provided the needed result. In the case of surfaces we did not find an equivalent result in the literature. We did find the results in \cite{baldo-manifolds} which suggest that the relaxation we consider is the right one on general manifolds. In the above reference a the authors do not prove a $\Gamma$-convergence result, but only the convergence of minimisers. We are concerned here only with smooth manifolds of codimension one and in this particular case it is possible to adapt classical methods in order to prove a $\Gamma$-convergence result.

We start by defining the space of functions of bounded variations on a $d-1$ dimensional surface in $\Bbb{R}^d$. Let $S$ be a smooth $d-1$ dimensional manifold without boundary in $\Bbb{R}^d$. In the following we consider the tangential gradient of a function $u$ defined on $S$ to be 
\[ \nabla_\tau u = \nabla \tilde u - (\nabla \tilde u . n)n,\]
where $\tilde u$ is a regular extension of $u$ in a neighbourhood of $S$ $n$ denotes the normal vector to the surface. In the same way we define the tangential divergence of a vector field $w \in C^1(S;\Bbb{R}^d)$ by 
\[ \di_\tau w = \text{tr}(D_\tau w)\]
where the matrix $D_\tau w$ contains on line $i$ the tangential gradient of the $i$-th component of $w$, i.e. $\nabla_\tau w_i$. See \cite[Section 5.4] {henrot-pierre} for further details.

We consider the space of functions with bounded variation on $S$ 
\[ BV(S) = \{ u \in L^1(S) : TV(u) <\infty\}\]
where
\[ TV(u) = \sup\{ \int_S u \di_\tau g : |g|_\infty\leq 1\}.\]
Using the divergence theorem on manifolds (see \cite[Section 5.4]{henrot-pierre}), we obtain that if $u$ is $C^1(S)$ then 
\[ TV(u) = \int_S |\nabla_\tau u|.\] 
If $\omega$ is a subset of $S$ we define its generalized perimeter as $\Per(\omega) = TV(\chi_\omega)$, where $\chi_\omega$ represents the characteristic function of $\omega$. By mimicking the proof in the euclidean case we can prove that the total variation is lower semi-continuous for the $L^1(S)$ convergence. We refer to \cite{braides2} for more details. 

Let $(C_i)$ is a set of local charts which cover $S$ such that each $C_i$ is diffeomorphic to a connected and bounded open subset  $D_i$ of $\Bbb{R}^{d-1}$. We denote by $\theta_i: D_i \to C_i$ these diffeomorphisms. Then it is possible to transfer a function $u$ from $C_i$ to $D_i$ using the transformation $\tilde u_i = u\circ \theta_i$. These new functions $\tilde u_i$, which lie now in Euclidean spaces, are functions of bounded variation. Therefore, it is possible to transfer some of the theory of BV functions from Euclidean spaces to manifolds of co-dimension $1$ by using local charts and partitions of unity. In particular, it is possible to approximate finite perimeter sets $\omega \subset S$ with smooth sets $\omega_n \subset S$ such that $\omega_n \to \omega$ in the $L^1(S)$ topology and $\Per(\omega_n) \to \Per(\omega)$.

We are now ready to state the relaxation result in the case of a single phase, which will be generalized later to the case of a partition. To derive the theorem below we follow the approach provided by Buttazzo in \cite{buttazzogconv} and Alberti in \cite{gammaconvalberti}.

\begin{thm}
Define $F_\varepsilon,F : L^1(S) \to [0,+\infty]$ as follows:
\[ F_\varepsilon(u) = \begin{cases}
\ds \int_S \left( \varepsilon |\nabla_\tau u|^2 +\frac{1}{\varepsilon} u^2(1-u)^2 \right)d\sigma & \text{ if } u \in H^1(S),\ \int_S u = c\\
+\infty & \text{ otherwise.}
\end{cases}
\]
\[ F(u) = \begin{cases}
\frac{1}{3}\Per(\{u = 1\}) & \text{ if } u \in BV(S,\{0,1\}),\ \int_S u = c\\
+\infty & \text{ otherwise.}
\end{cases}
\]
Then $F_\varepsilon \gconv F$ in the $L^1(S)$ topology.
\label{gconv1}
\end{thm}

\emph{Proof:} We define $\phi(t) = \int_0^t |s(1-s)| ds$. We consider a sequence $(u_\varepsilon) \to u$ in $L^1(S)$ such that $\liminf_{\varepsilon \to 0} F_\varepsilon(u_\varepsilon)<+\infty$. Since $F_\varepsilon(u_\varepsilon) \geq \frac{1}{\varepsilon} \int_S u_\varepsilon^2(1-u_\varepsilon)^2$, if we take a subsequence of $u_\varepsilon$ which converges almost everywhere to $u$ we obtain that
\[ \int_S u^2(1-u)^2 = 0,\]
and thus $u \in \{0,1\}$ almost everywhere in $S$. Note that truncating $u_\varepsilon$ between $0$ and $1$ decreases the value of $F_\varepsilon(u_\varepsilon)$ while preserving the fact that $u_\varepsilon \to u$ in $L^1(S)$. Also note that $\phi$ is Lipschitz on $[0,1]$ so we can conclude that $\phi \circ u_\varepsilon \to \phi \circ u$ in $L^1(S)$. By applying the classical inequality $a^2+b^2\geq 2ab$ we get that
\[ F_\varepsilon(u_\varepsilon) \geq 2\int_S |\nabla_\tau u| \phi'(u_\varepsilon) = 2 \int_S |\nabla_\tau (\phi\circ u_\varepsilon)|.\]
Taking $\liminf$ in the above inequality and using the semi-continuity of the total variation with respect to the $L^1(S)$ convergence we obtain that
\[ \liminf_{\varepsilon \to 0} F_\varepsilon(u_\varepsilon) \geq2 TV(\phi \circ u) = 2\phi(1) TV(u).\]
Since $u$ is a characteristic function, it follows that the perimeter of $\{u=1\}$ is bounded and therefore $u \in BV(S,\{0,1\})$. Note that $\phi(1) = 1/6$ and thus we recover the desired constant in front of the perimeter. It is obvious that the integral condition is also preserved in the limit. This concludes the proof of the $\Gamma-\liminf$ part of the theorem.

For the $\Gamma-\limsup$ part we need to exhibit a recovery sequence for each $u$ such that $F(u) <+\infty$. By a classical argument it is enough to find a recovery sequence only for functions $u$ which are characteristic functions of smooth sets in $S$. See \cite{braides2} for more details concerning the reduction to regular sets and \cite[Theorem 3.42]{ambrosiofuscopallara} for the BV approximation of finite perimeter sets with smooth sets. 

Let's consider now $u = \chi_\omega$ where $\omega \subset S$ is a set with smooth boundary relative to $S$. We consider the signed distance function $d_\omega : S \to \Bbb{R}$ defined by
\[ d_\omega(x) = d_\tau(x,S\setminus \omega)-d_\tau(x,\omega),\]
where $d_\tau$ is the geodesic distance on $S$. Note that $d_\omega$ is positive outside $\omega$ and negative inside. Consider the optimal profile problem
\[ c = \min\{ \int_{\Bbb{R}} (W(v)+|v'|^2 : v(-\infty) = 0,\ v(+\infty) = 1\}.\]
Any solution of this minimizing problem satisfies $v' = \sqrt{W(v)}$ and we can impose the initial condition $v(0)=1/2$ in order to have a symmetric behaviour. We can see that the optimal value is $c = 2\int_0^1\sqrt{W(s)}ds$. In our problem we have chosen $W(s) = s^2(1-s)^2$. In order to have a function which goes from $0$ to $1$ in finite time we may choose 
\[ v^\eta = \min\{ \max\{0,(1+2\eta)v-\eta\},1\}.\]
We see that as $\eta \to 0$ we have 
\[ c^\eta = \int_{\Bbb{R}}(W(v^\eta)+|(v^\eta)'|^2) \to c \text{ as } \eta \to 0.\]
All these considerations are inspired from \cite{braides2}. We can define
\[ u_\varepsilon(x) = v^\eta(d_\omega(x)/\varepsilon).\]
We can see that 
\begin{align*}
F_\varepsilon(u_\varepsilon) & = \int_S\left( \varepsilon |\nabla_\tau u|^2 +\frac{1}{\varepsilon} W(u) \right) \\
  & = \int_{-T\varepsilon}^{T/\varepsilon} \int_{d_\omega(x) = t} \left(\varepsilon |(v^\eta)'(d_\omega(x)/\varepsilon)|^2\frac{|\nabla_\tau d_\omega(x)|^2}{\varepsilon^2}+\frac{1}{\varepsilon} W(v^\eta(d_\omega/\varepsilon))\right)d\mathcal{H}^{d-2}(x) dt\\
  & = \int_{-T/\varepsilon}^{T/\varepsilon} \int_{d_\omega(x)=t}\frac{1}{\varepsilon} (|(v^\eta)'(t/\varepsilon)|^2+W(v^\eta(t/\varepsilon)))d\mathcal{H}^{d-2}(x) dt  \\
  & = \int_{-T/\varepsilon}^{T/\varepsilon} \Per(d_\omega(x)=t) \frac{1}{\varepsilon} (|(v^\eta)'(t/\varepsilon)|^2+W(v^\eta(t/\varepsilon)))dt \\
  & = \int_{-T}^{T} \Per(d_\omega(x) = t\varepsilon) (|(v^\eta)'(t)|^2+W(v^\eta(t)))dt
\end{align*}
where we have applied the co-area formula and $T$ is chosen such that the support of $v^\eta$ is inside $[-T,T]$. Since $\lim_{s \to 0} \Per(\{d_\omega(x) = s\}) = \Per(\omega)$ we see that for $\varepsilon$ small enough there exists $\delta$ such that $\Per(d_\omega(x) = s) <\Per(\omega)+\delta$ when $|s|<T\varepsilon$. Therefore
\[ \limsup_{\varepsilon \to 0}F_\varepsilon(u_\varepsilon) \leq (\Per(\omega)+\delta)\int_{-T}^T (|(v^\eta)'(t)|^2+W(v^\eta(t)))dt=(\Per(\omega)+\delta)c_\eta.\]
Since this is true for any $\delta,\eta$ small enough, by letting $\delta,\eta \to 0$ we obtain the desired result.

In order to have a fixed integral equal to $\int_S \chi_\omega = c$ it is enough to consider a shift in the definition of $u_\varepsilon$:
\[ u_\varepsilon(x) = v^\eta((d_\omega(x)+s_\varepsilon)/\varepsilon),\]
where $s_\varepsilon \in [-T\varepsilon,T\varepsilon].$ We can see that for $s_\varepsilon = T\varepsilon$ we have $u_\varepsilon = 1$ on $\omega$ and thus $\int_S u_\varepsilon>c$ while for $s_\varepsilon = -T\varepsilon$ the support of $u_\varepsilon$ is included in $\omega$ and we have the opposite inequality. Thus, for each $\varepsilon$ small enough we can change the definition of $u_\varepsilon$ so that $\int_S u_\varepsilon =c$. The estimates presented above are carried with no difficulty in this setting.
\hfill $\square$

We can now state the result in the partitioning case. We denote by $\bo u$ an element in $(L^1(S))^n$. In order to simplify the notations we introduce the space
\[ X = \{ \bo u \in (L^1(S))^n : \int_S u_i = c_i,\ \sum_{i=1}^n u_i = 1\}\]
where $c_i$ satisfy the compatibility condition $\sum_{i=1}^n c_i = \mathcal{H}^{d-1}(S)$. It is easy to see that $X$ is closed under the convergence in $(L^1(S))^n$.

\begin{thm}
Define $F_\varepsilon,F : (L^1(S)^n \to [0,+\infty]$ as follows:
\[ F_\varepsilon(\bo u) = \begin{cases}
\ds \sum_{i=1}^n\int_S \left( \varepsilon |\nabla_\tau u_i|^2 +\frac{1}{\varepsilon} u_i^2(1-u_i)^2 \right)d\sigma & \text{ if } \bo u \in (H^1(S))^n \cap X\\
+\infty & \text{ otherwise}
\end{cases}
\]
\[ F(\bo u) = \begin{cases}
\frac{1}{3}\sum_{i=1}^n\Per(\{u_i = 1\}) & \text{ if } \bo u \in (BV(S,\{0,1\}))^n  \cap X\\
+\infty & \text{ otherwise}
\end{cases}
\]
Then $F_\varepsilon \gconv F$ in the $(L^1(S))^n$ topology.
\label{gconv2}
\end{thm}

\emph{Proof:} It is easy to see that the $\Gamma-\liminf$ part follows at once from Theorem \ref{gconv1} and from the fact that $X$ is closed under the topology of $(L^1(S))^n$. 

In order to construct the recovery sequence we reduce the problem to the case where the limit $\bo u$ is consists of piecewise smooth parts in $S$. In this case we define $u_i = v^\eta(d_{\omega_i}(x)/\varepsilon)$ as in the one phase case. Thus on each $\omega_i$ we have $u_i \geq 1/2$ which implies that $\sum_{i=1}^n u_i \geq 1/2$. There are two points which need to be addressed:
\begin{enumerate}
\item The sum equal to $1$ condition. Due to the symmetry of the optimal profile we deduce that there is only one zone where the sum condition is not satisfied and that is in the neighborhood of singular points. Since an $\varepsilon$-neighborhood of the singular set is of order $\varepsilon^{d-1}$. Replacing each $u_i$ by $u_i/(\sum_{i=1}^n u_i)$ in these problematic regions we preserve the regularity of each $u_i$ and we note that the functions have bounded gradient of order $O(1/\varepsilon)$. We immediately find that the corresponding energy
\[ \int_{N_\varepsilon} \left(\varepsilon |\nabla_\tau u_i|^2 +\frac{1}{\varepsilon} u_i^2(1-u_i)^2 \right)\] 
vanishes as $\varepsilon \to 0$.
\item We also need to modify the functions $u_i$ so that they have the same integral over $S$. In order to do this we apply a procedure found in \cite{ambrosio} where we consider a family of balls in regions where $u_i \in \{0,1\}$. On each such ball we can consider modifications of $u_i$ such that the sum is preserved and the integrals have the right value. As above, the sum of energies on these balls will be negligible in the limit.
\end{enumerate}
Once these points are addressed, the $\limsup$ estimates follows just like in the one dimensional case and the proof of the theorem is completed.
\hfill $\square$


\section{Finite Element framework}
We wish to use this relaxation by $\Gamma$-convergence to perform numerical computations so we need a framework which allows us to compute the quantity
\[ \varepsilon \int_S |\nabla_\tau u|^2 +\frac{1}{\varepsilon} \int_S u^2(1-u)^2,\]
in fast, efficient way. In order to do this we triangulate the surface $S$ and we compute the mass matrix $M$ and the stiffness matrix $K$ associated to the $P_1$ finite elements on this triangulation. Then, if for the sake of simplicity, we use the same notation $u$ for the $P_1$ finite element approximation of $u$, we have
\[ \int_S |\nabla_\tau u|^2 = u^T Ku\]
and 
\[ \int_S u^2(1-u)^2 = v^T M v,\]
where $v = u.^2.\times (1-u).^2$. We have used the Matlab convention that adding a point before an operation means that we are doing component-wise vector computations. Note that once the matrices $K,M$ are computed, we only have to perform matrix-vector multiplications, which is really fast. In this setting we use the discrete gradients of the above expressions given by:
\[ \nabla_u u^T Ku = 2Ku,\]
\[ \nabla_u v^T M v = 2Mv.\times (1-2u).\]

The partition condition and the equal areas constraint are imposed by making an orthogonal projection on the linear constraints as follows. We write the discrete vectors representing $P_1$ discretization of the density functions in the following matrix form
\[ M = (\varphi^1 \  \varphi^2 \ ... \ \varphi^n).\]
The partition constraint implies that the sum of the elements on every line of $M$ is equal to $1$ and the equal area constraint implies that for every column of the matrix $M$ we have the relation
\[ \langle v,\varphi^i \rangle = A/n, \text{ where } v = {\bf 1}_{1\times N}\cdot M.\]
Here the constant $A$ is the total area of the surface, $N$ is the total number of points in the triangulation and the notation ${\bf 1}_{p\times q}$ represents the $p\times q$ matrix whose entries are all equal to $1$. These conditions are discretizations in the finite element setting of the conditions that the integrals of the density functions $u_i$ are all equal to $A/n$. Indeed, given a triangulation $\mathcal{T}$ of $S$ and its associated mass matrix $M$, we have $\ds \int_S 1\cdot u_i = {\bf 1}_{1\times N}\cdot M \cdot \varphi^i$, where $\varphi^i$ is the vector containing the values of $u_i$ at the vertices of the triangulation.
 The projection routine can be found in Algorithm \ref{projection-perimeter}. 
 
 \begin{algorithm}
 \caption{Orthogonal projection on the partition and area constraints}
 \label{projection-perimeter}
 \begin{algorithmic}[1]
 \Require $A = (a_{ij}) \in \Bbb{R}_{N\times n} $, $c \in \Bbb{R}_{1\times n}$, $d \in \Bbb{R}_{N \times 1}$, $v$
 \State $(e_i) = \sum_j a_{ij}-c_i$    (line sum error; $N \times 1$ column vector)
 \State $(f_i) = \sum_i v_i a_{ij}-d_j$ (column scalar product error; $n\times 1$ column vector)
 \State Define the matrix $C$ of size $n \times n$ by
  \[ \begin{cases}
  c_{kl} = \|v\|_2^2/n & \text{ if } k \neq l \\
  c_{kk} = \|v\|_2^2 - \|v\|_2^2/n
  \end{cases}\]
 
 \State $(q_j) = (f_j) -  \langle v, e\rangle/n$ ($n \times 1$ column vector)
 \State Compute $(\lambda_j) \in \Bbb{R}_{n \times 1}$ with $\lambda_n = 0$ such that $C|_{(n-1)\times (n-1)} (\lambda_j)|_{n-1}=(q_j)|_{n-1}$. The indices indicate a sub-matrix with the first $n-1$ lines and columns, or the sub-vector formed by the first $n-1$ components.
 \State $S = \sum_j \lambda_j$
 \State $\eta_i = (e_i-S\cdot v_i)/n$ ($N\times 1$ column vector)
 \State $A_{\text{orth}} = (\eta_i)\cdot {\bf 1}_{1\times n}+v\cdot (\lambda_j)^T$, where $ {\bf 1}_{p\times q}$ is the $p\times q$ matrix with all entries equal to $1$
 \State $A = A-A_{\text{orth}}$
 
 \Return $A$
 \end{algorithmic}
 \end{algorithm}

Once we have this discrete formulation we use an optimized LBFGS gradient descent procedure \cite{lbfgs} to compute the numerical minimizers. In order to avoid local minima where one of the phases $\varphi^l$ is constant, which arise often when the number of phases is greater than $5$, we add a Lagrange multiplier which penalizes the constant functions. In this way, we optimize
\[ \sum_{i=1}^n \varepsilon \int_S |\nabla_\tau \varphi^i|^2+\frac{1}{\varepsilon} \int_S (\varphi^i)^2(1-\varphi^i)^2 +\lambda (\text{std}(\varphi^i)-\text{starget})^2,\]
where $\text{std}(\varphi^l)$ is the standard deviation of $\varphi^l$ and $\text{starget}$ is the standard deviation of a characteristic function of area $\text{Area}(S)/n$. 

\begin{figure}
\centering
\includegraphics[width= 0.19\textwidth]{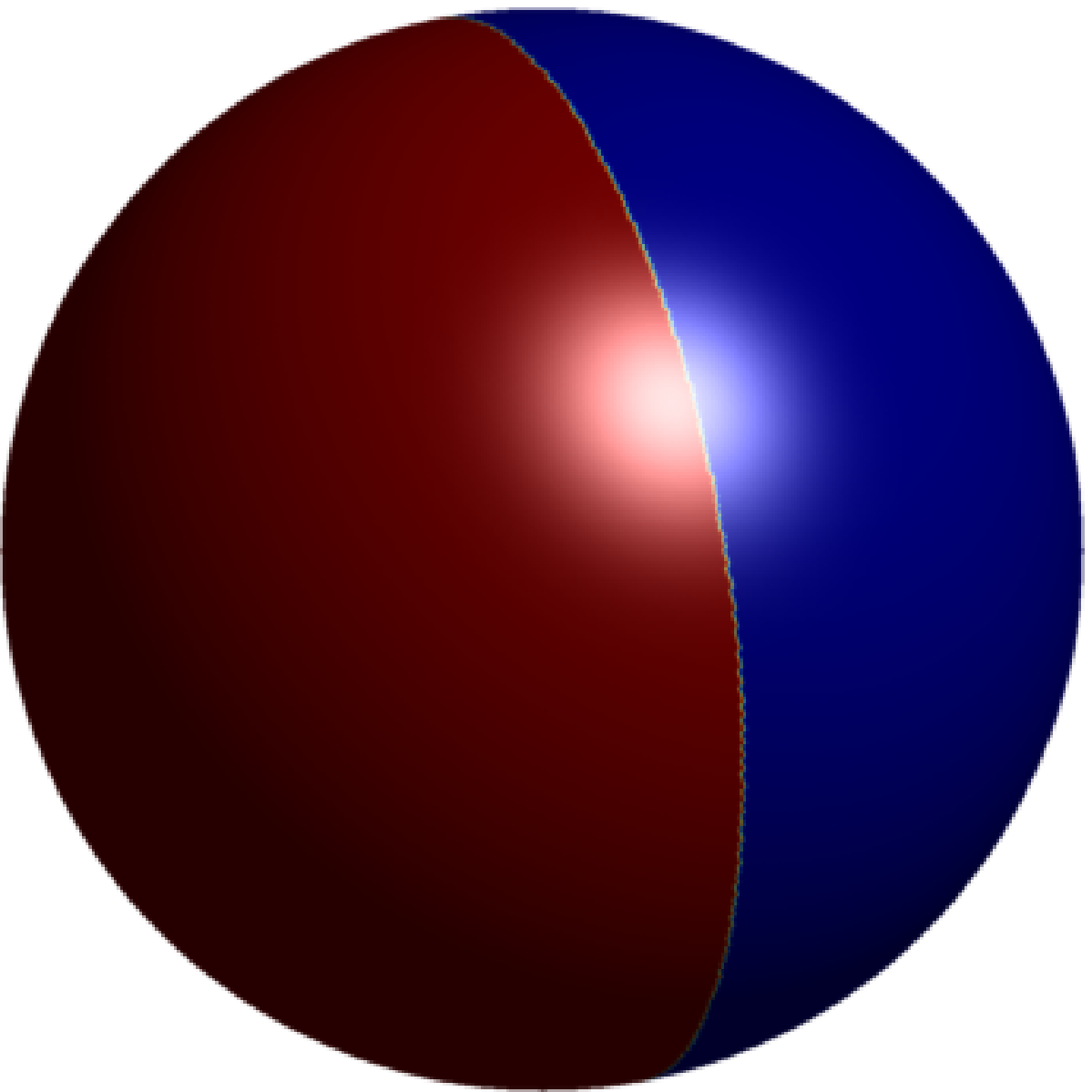}~
\includegraphics[width= 0.19\textwidth]{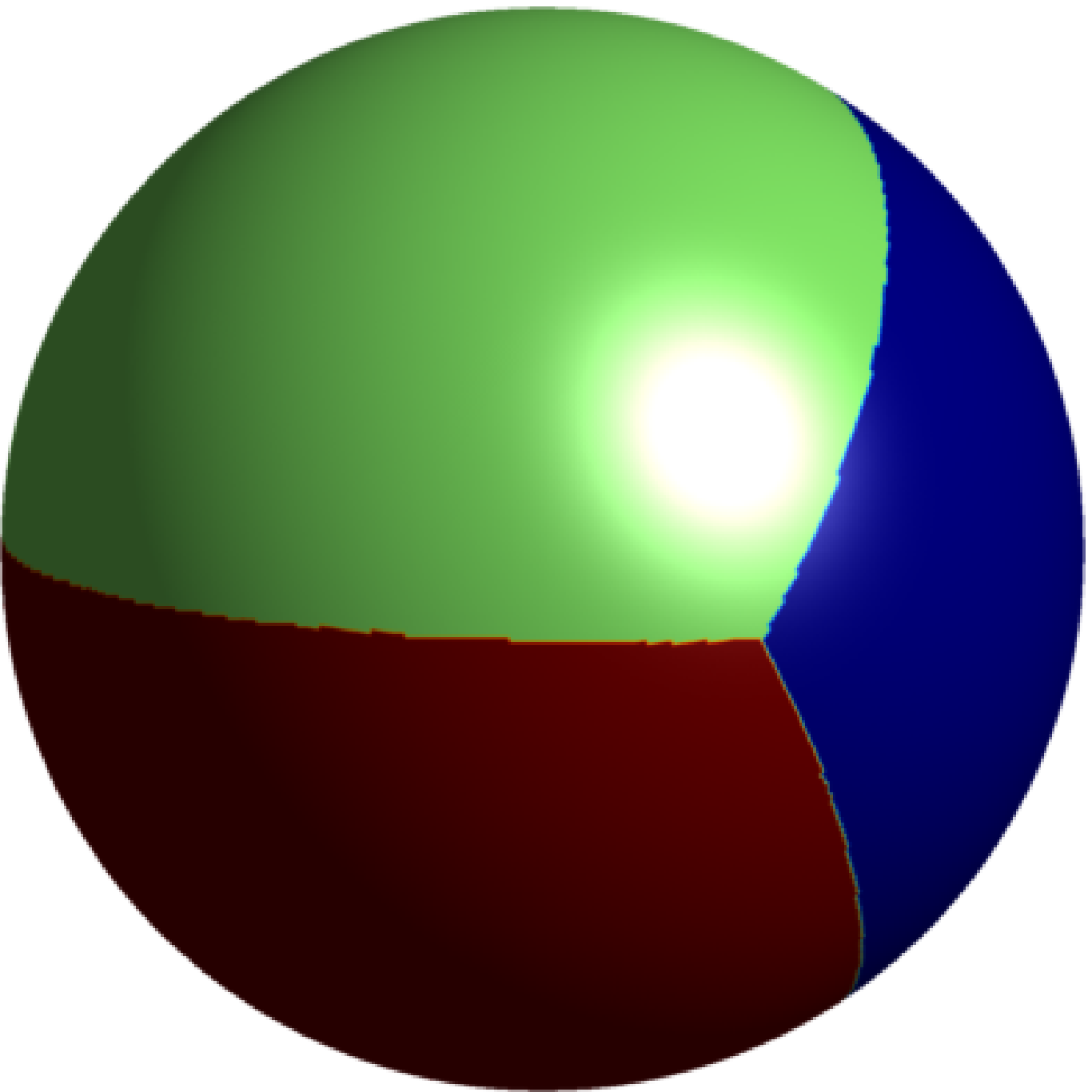}~
\includegraphics[width= 0.19\textwidth]{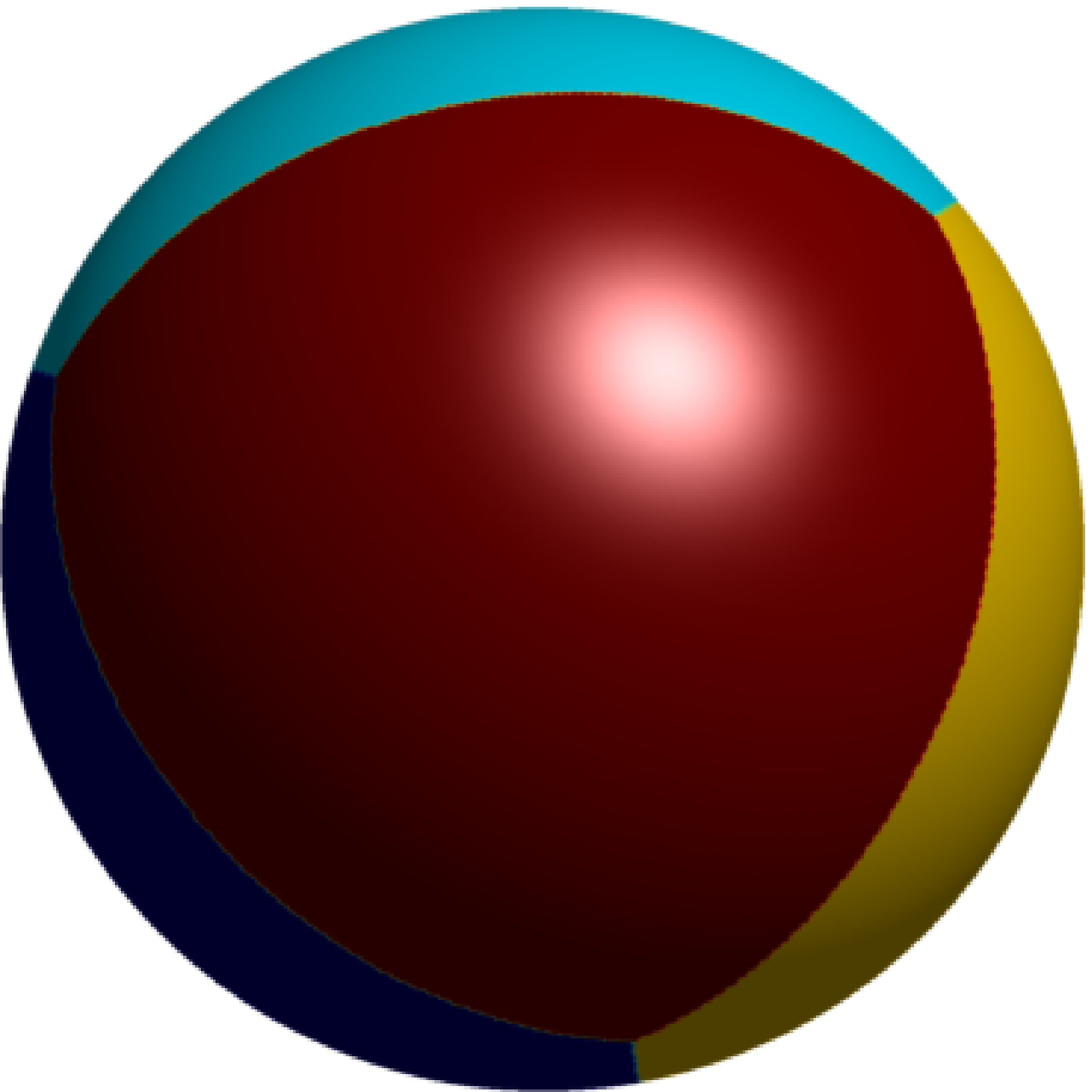}~
\includegraphics[width= 0.19\textwidth]{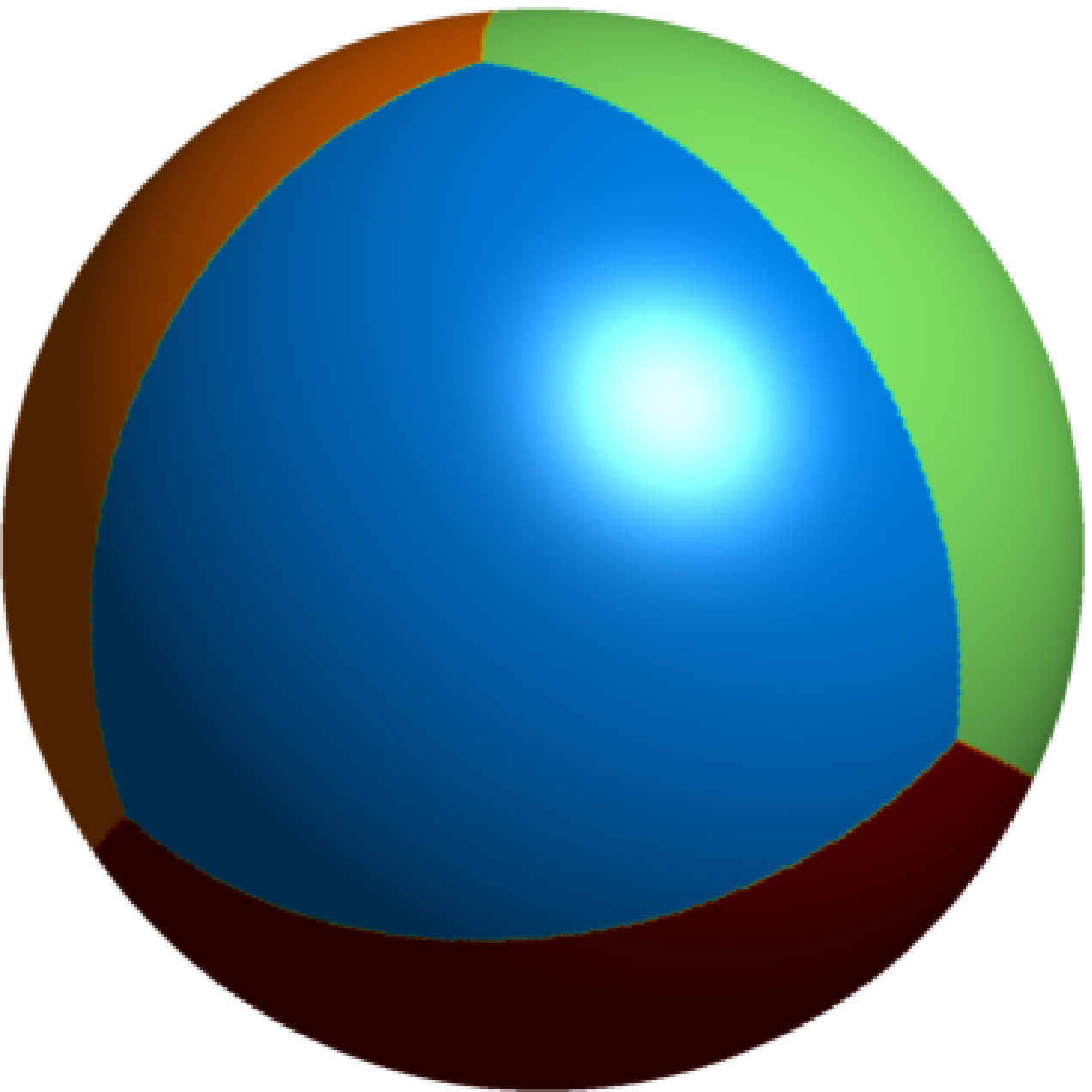}~
\vspace{0.1cm}

\includegraphics[width= 0.19\textwidth]{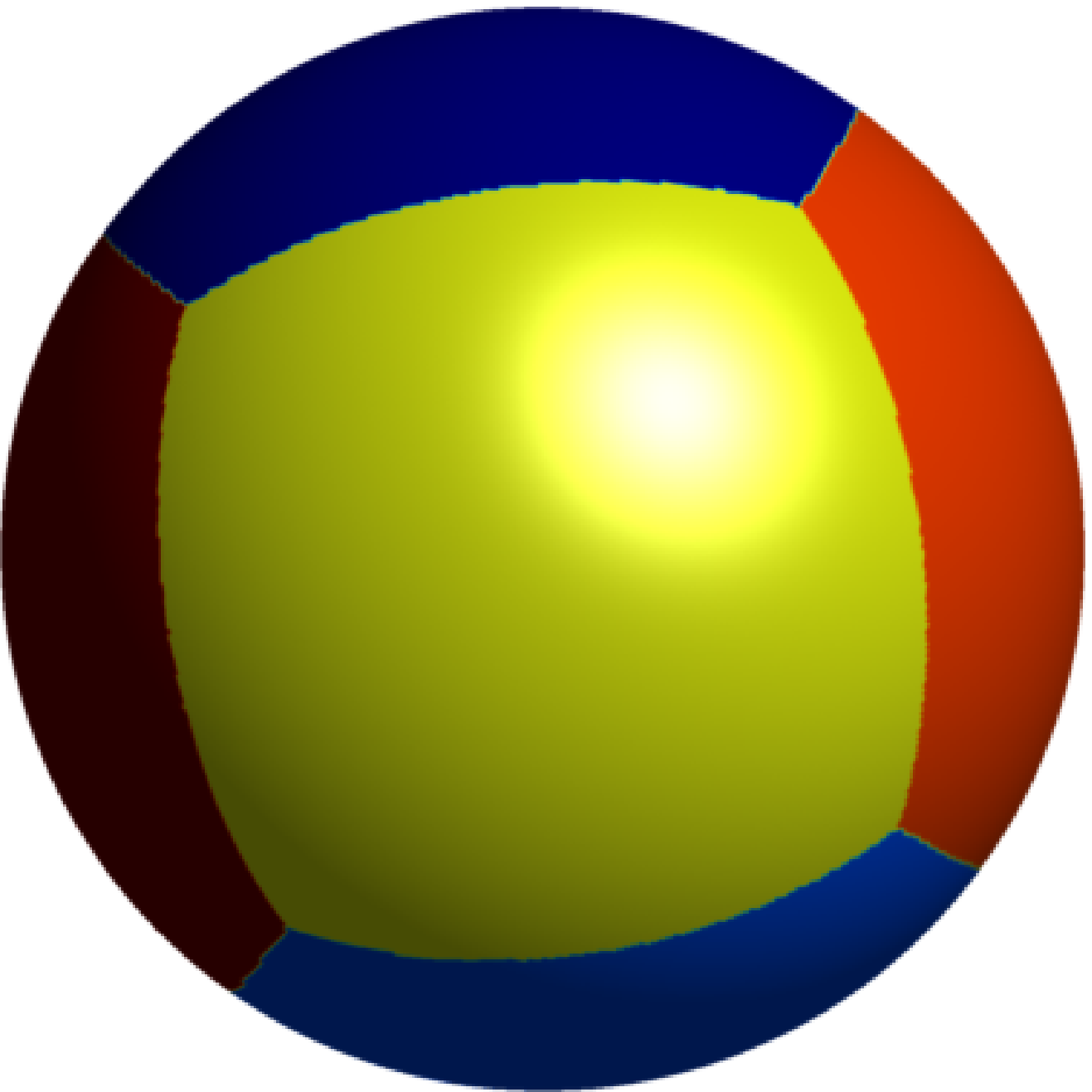}~
\includegraphics[width= 0.19\textwidth]{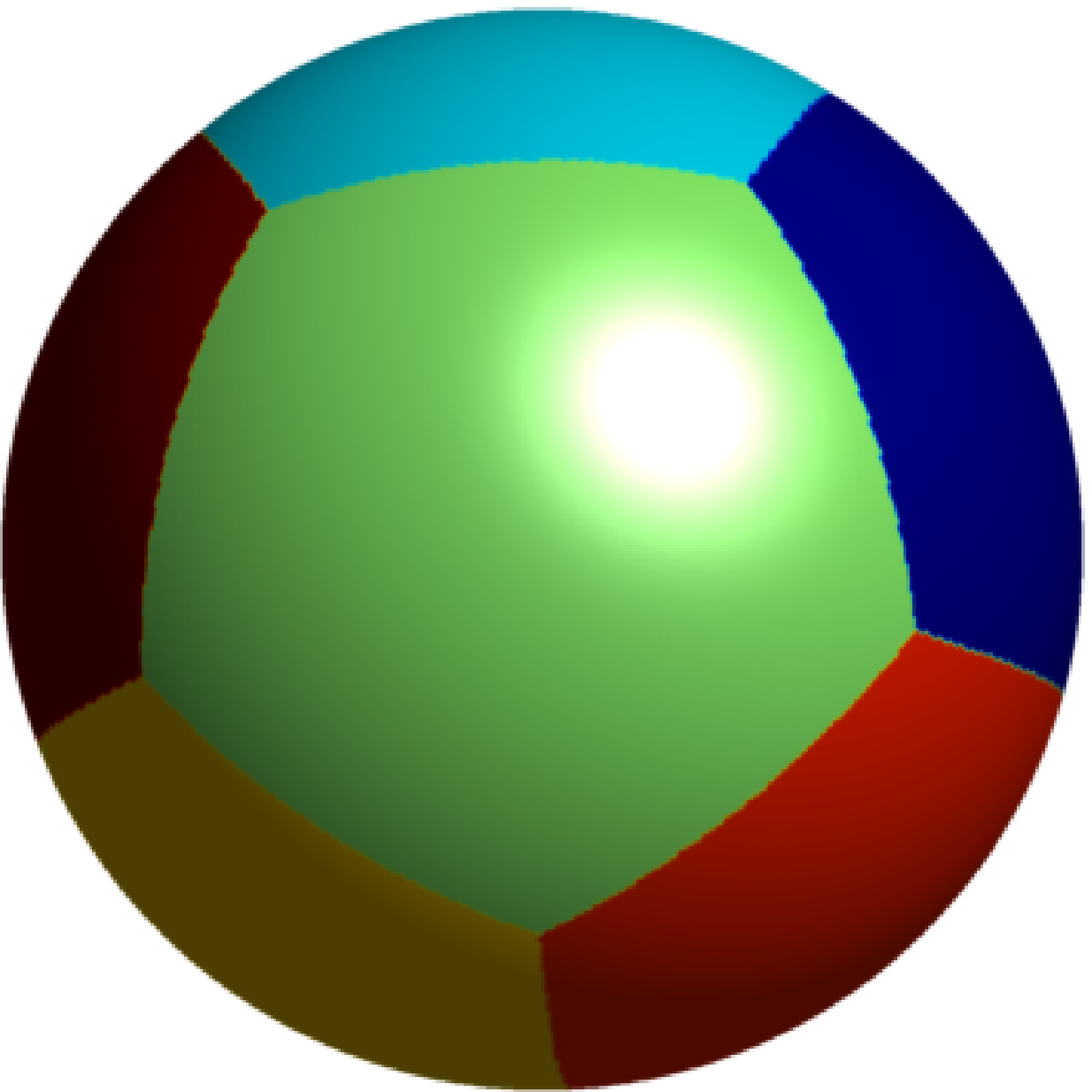}~
\includegraphics[width= 0.19\textwidth]{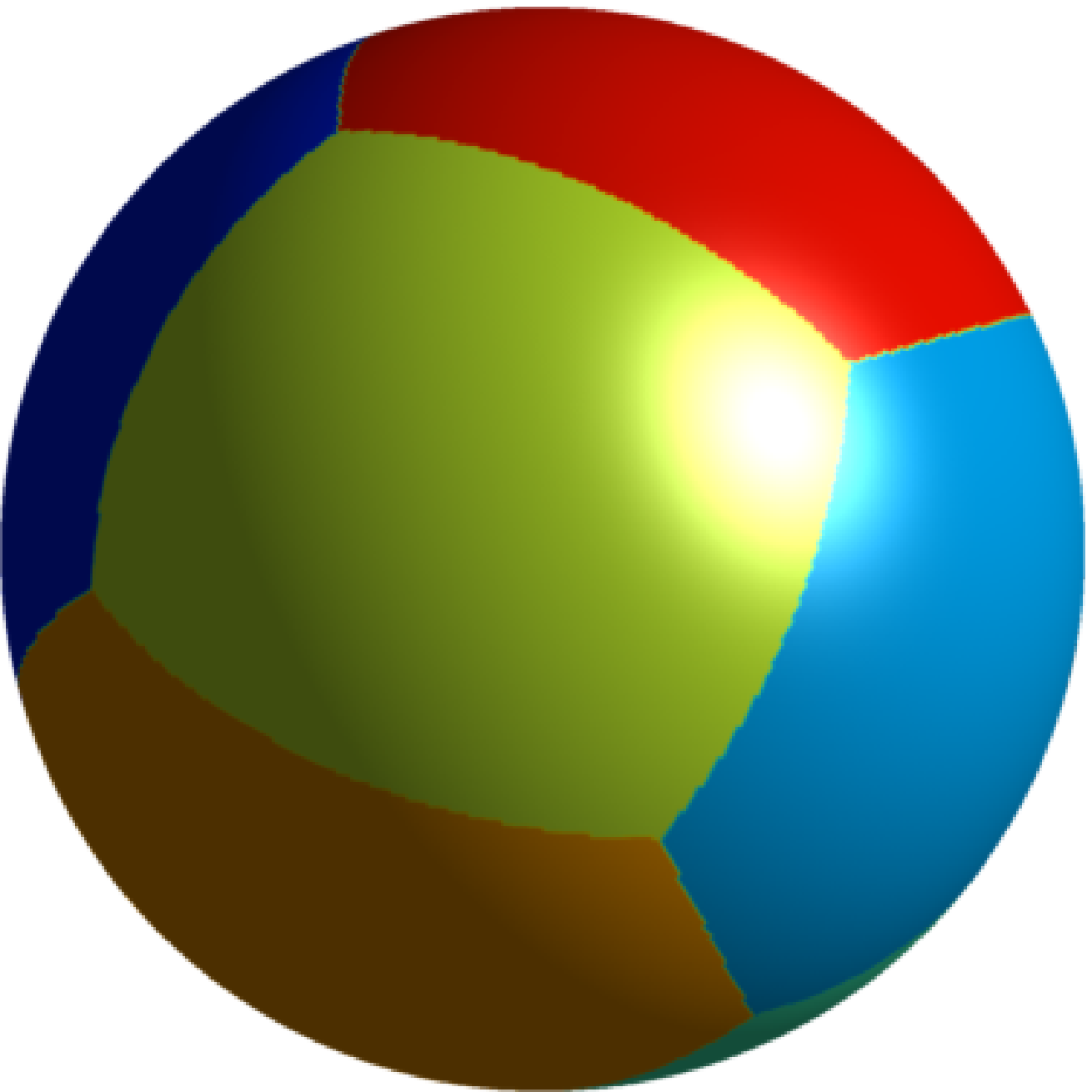}~
\includegraphics[width= 0.19\textwidth]{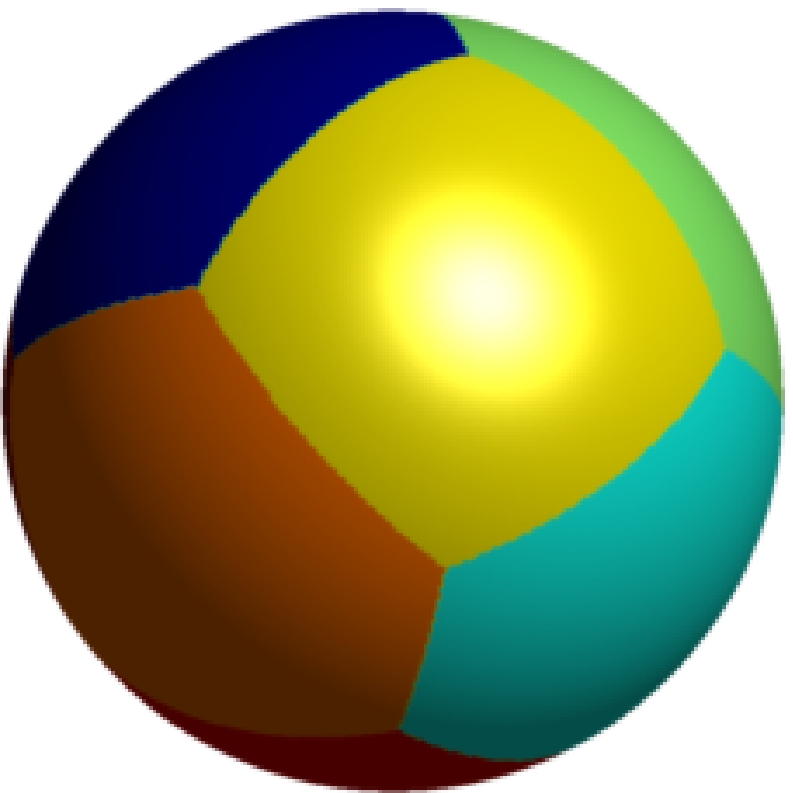}~
\vspace{0.1cm}

\includegraphics[width= 0.19\textwidth]{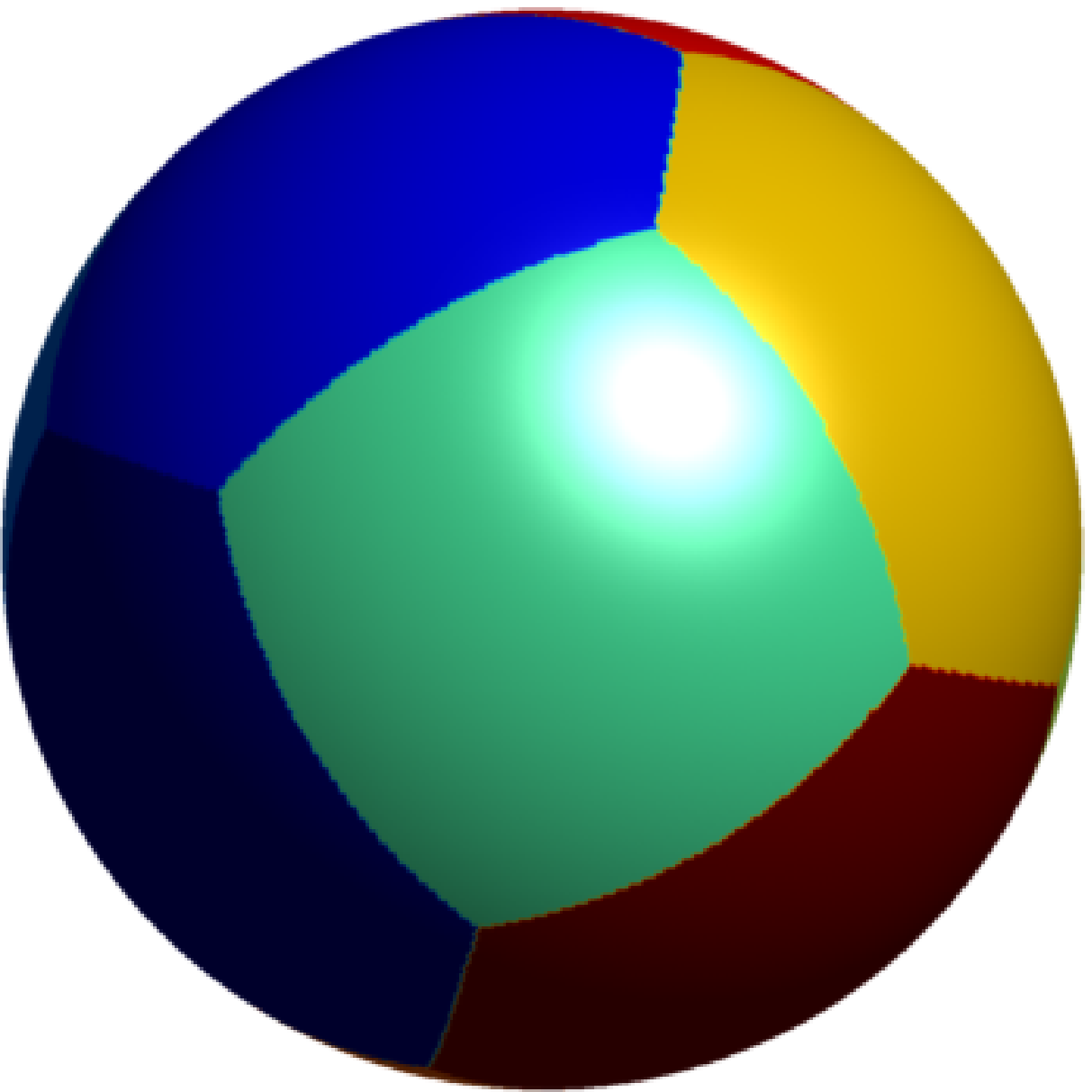}~
\includegraphics[width= 0.19\textwidth]{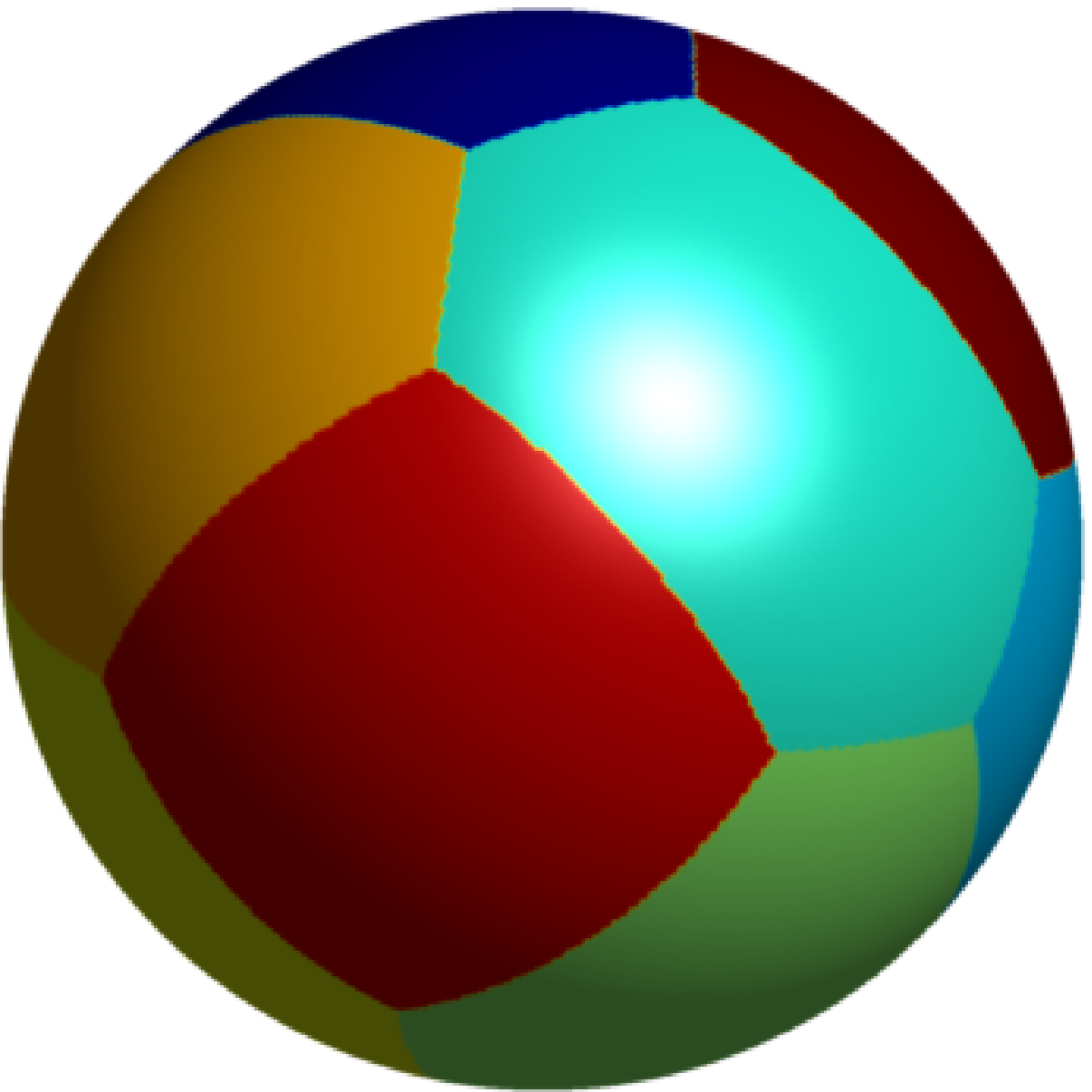}~
\includegraphics[width= 0.19\textwidth]{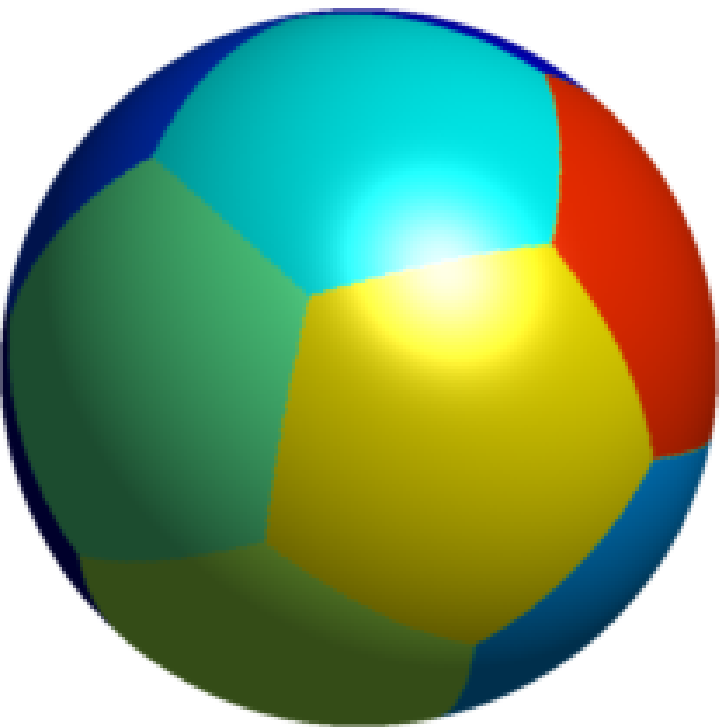}~
\includegraphics[width= 0.19\textwidth]{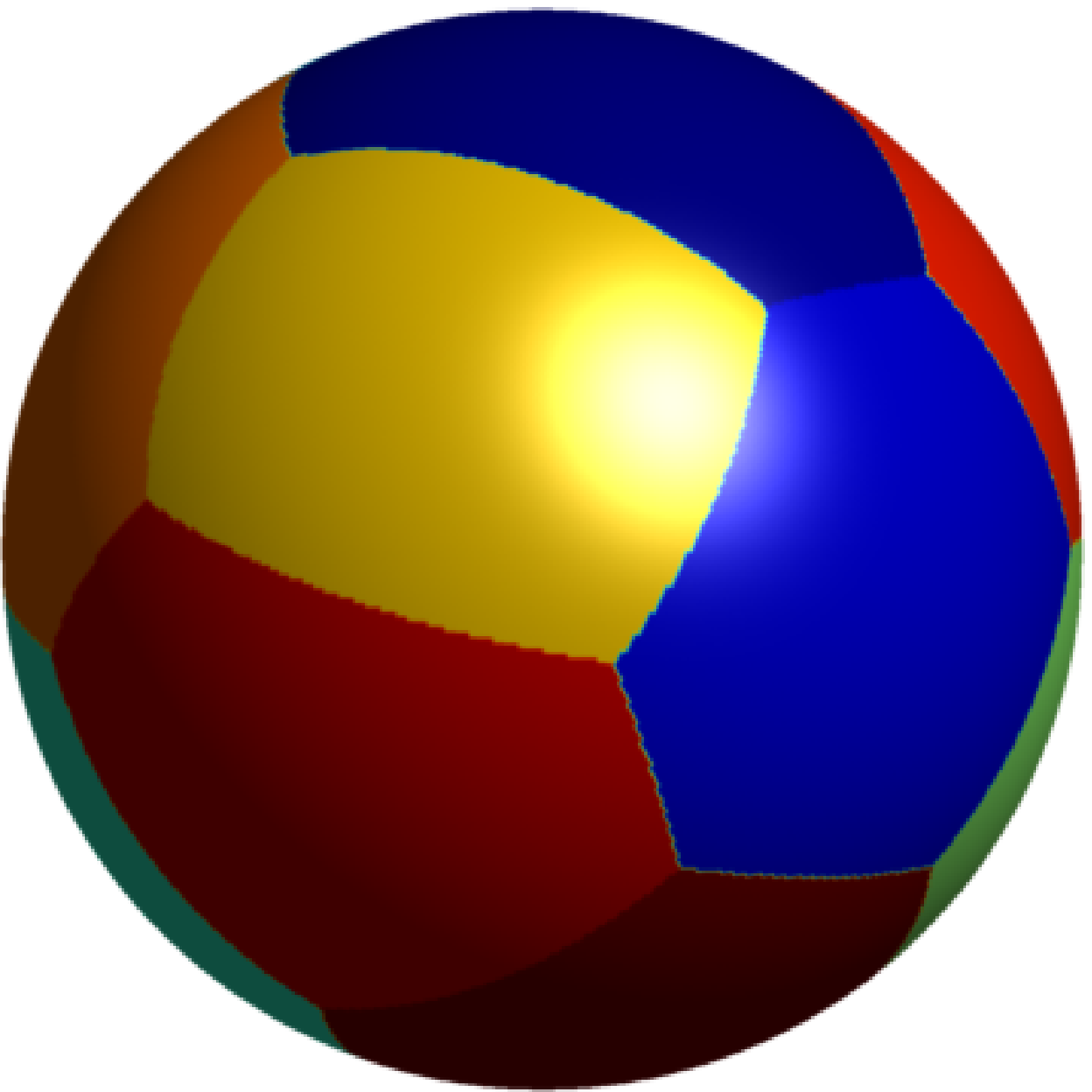}~
\vspace{0.1cm}

\includegraphics[width= 0.19\textwidth]{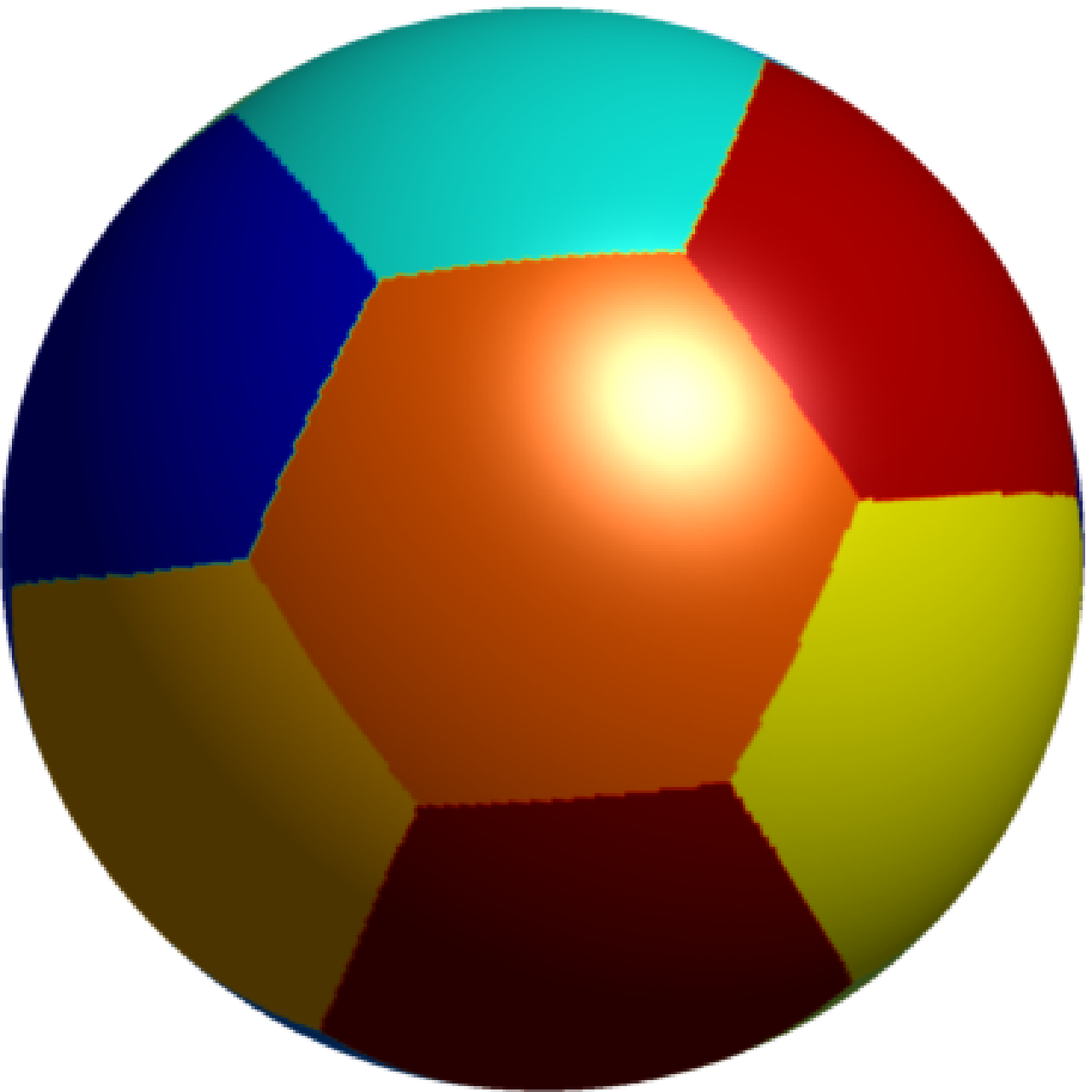}~
\includegraphics[width= 0.19\textwidth]{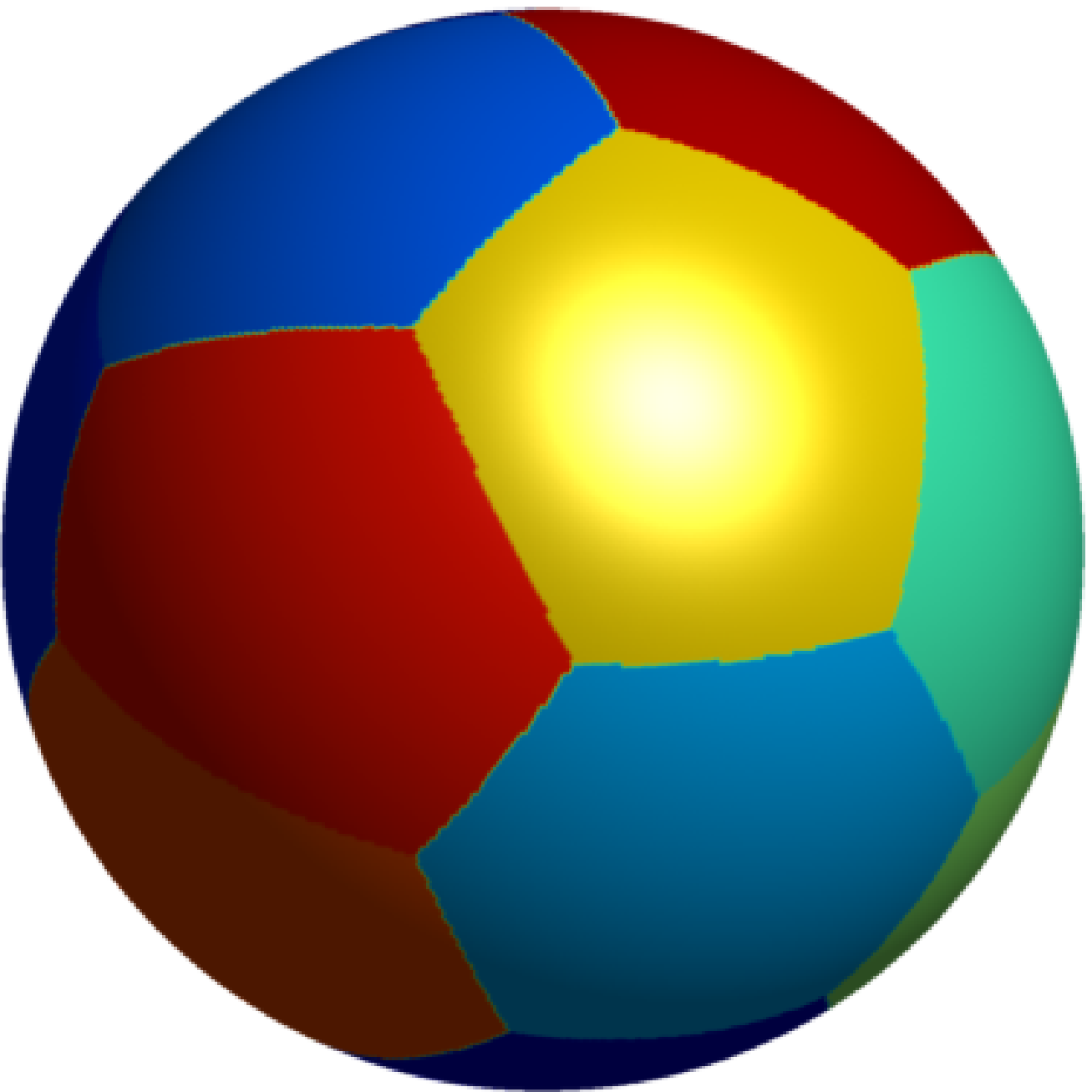}~
\includegraphics[width= 0.19\textwidth]{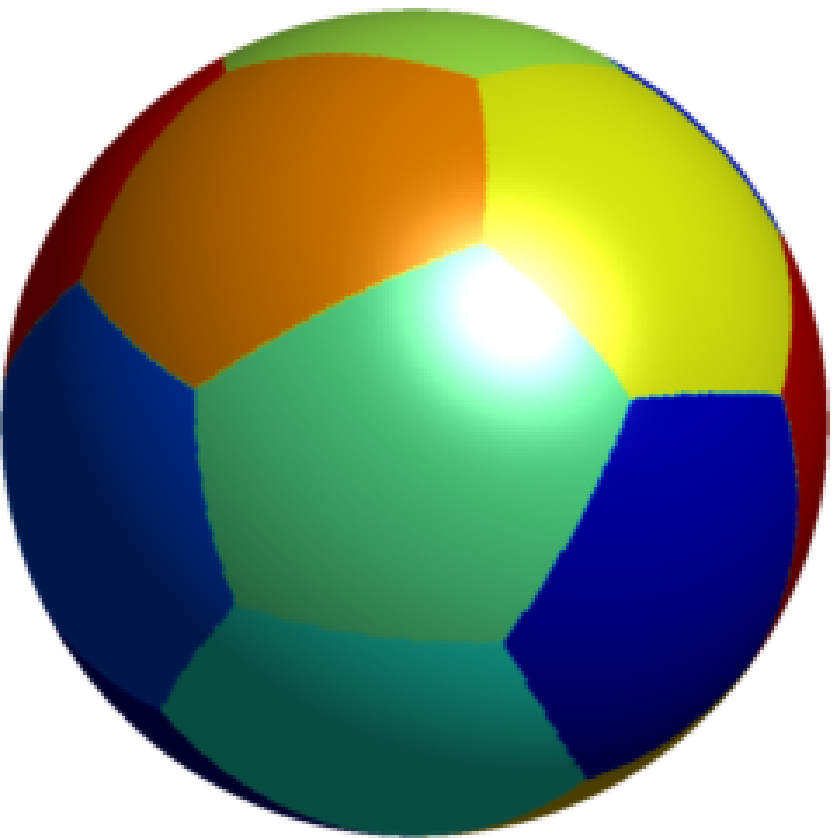}~
\includegraphics[width= 0.19\textwidth]{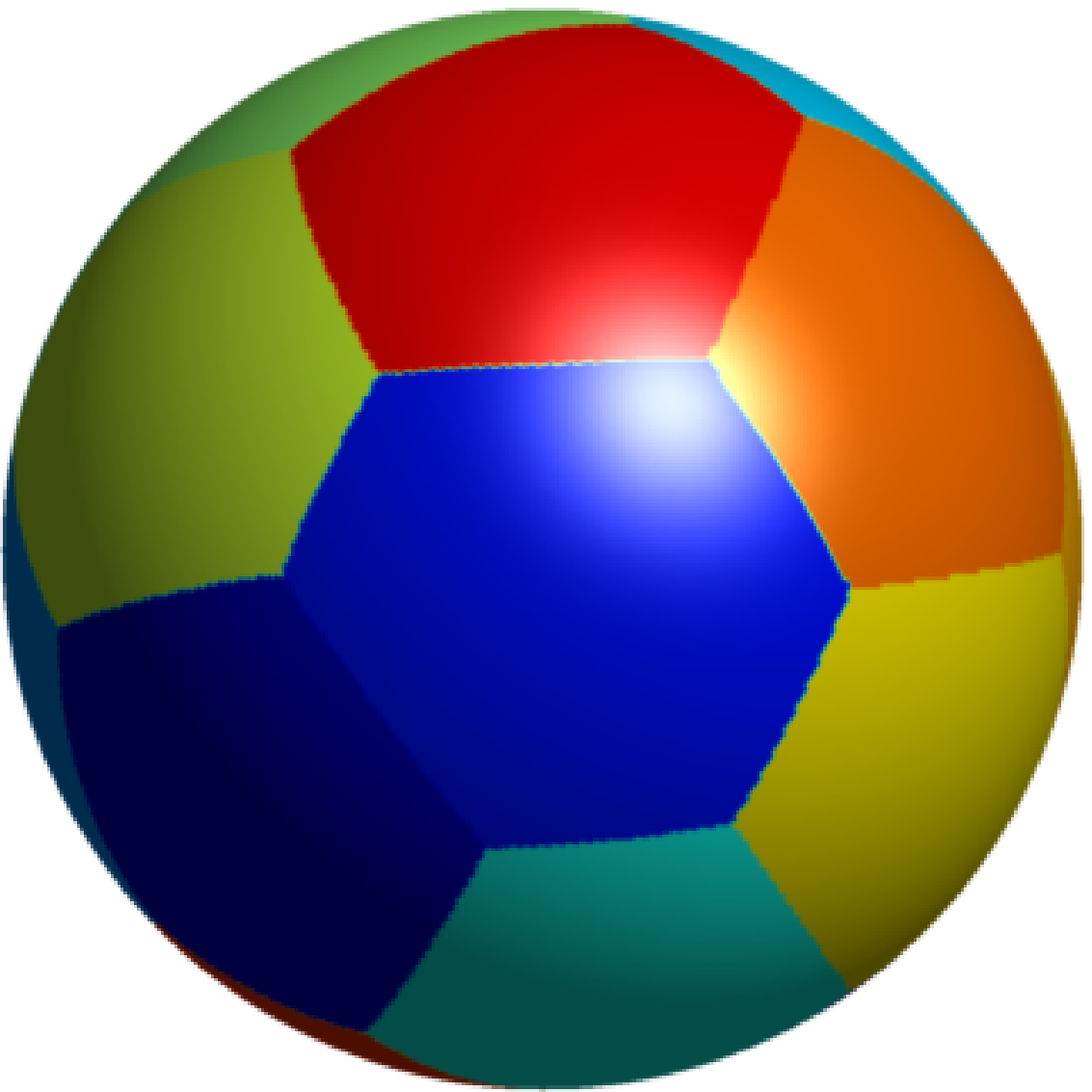}~
\vspace{0.1cm}

\includegraphics[width= 0.19\textwidth]{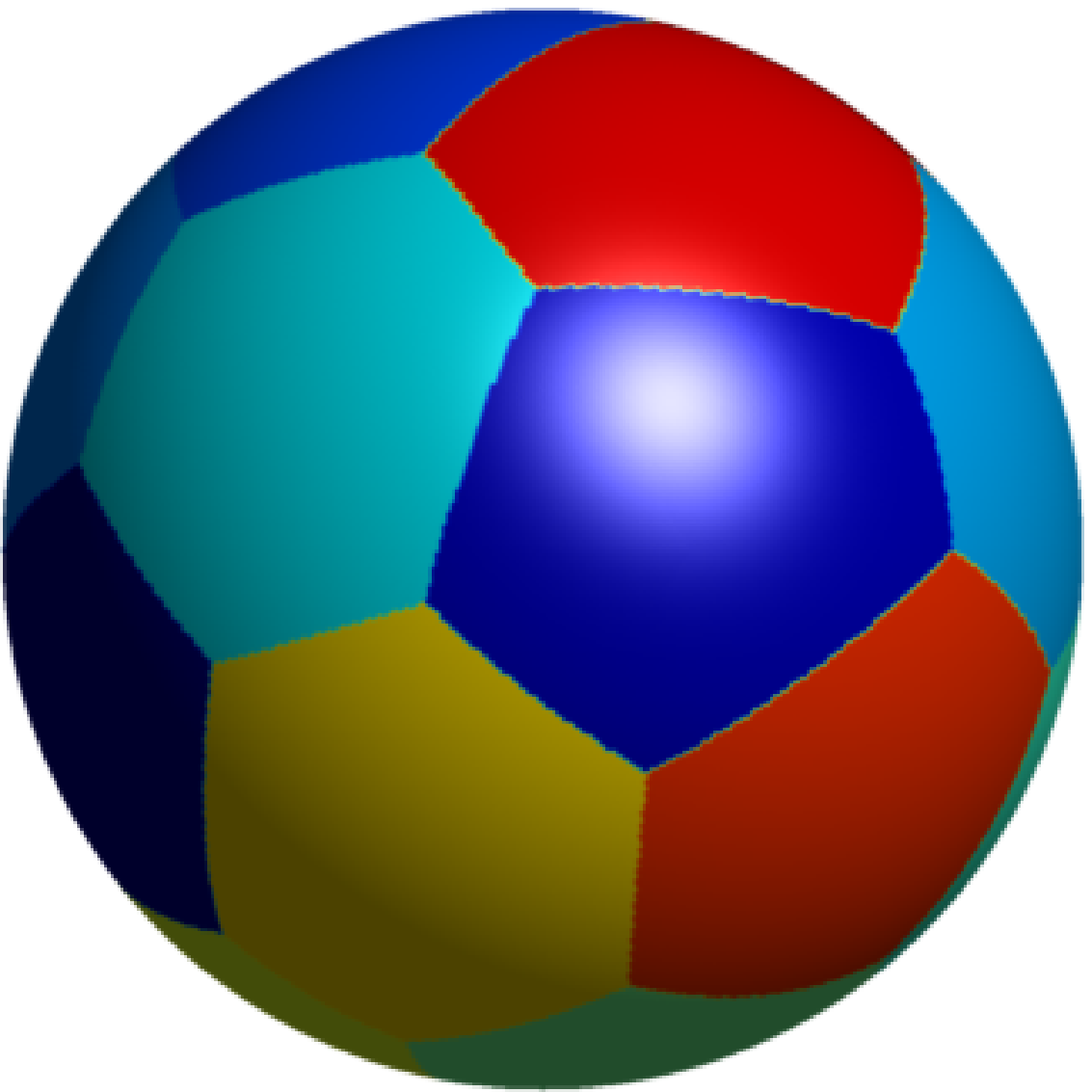}~
\includegraphics[width= 0.19\textwidth]{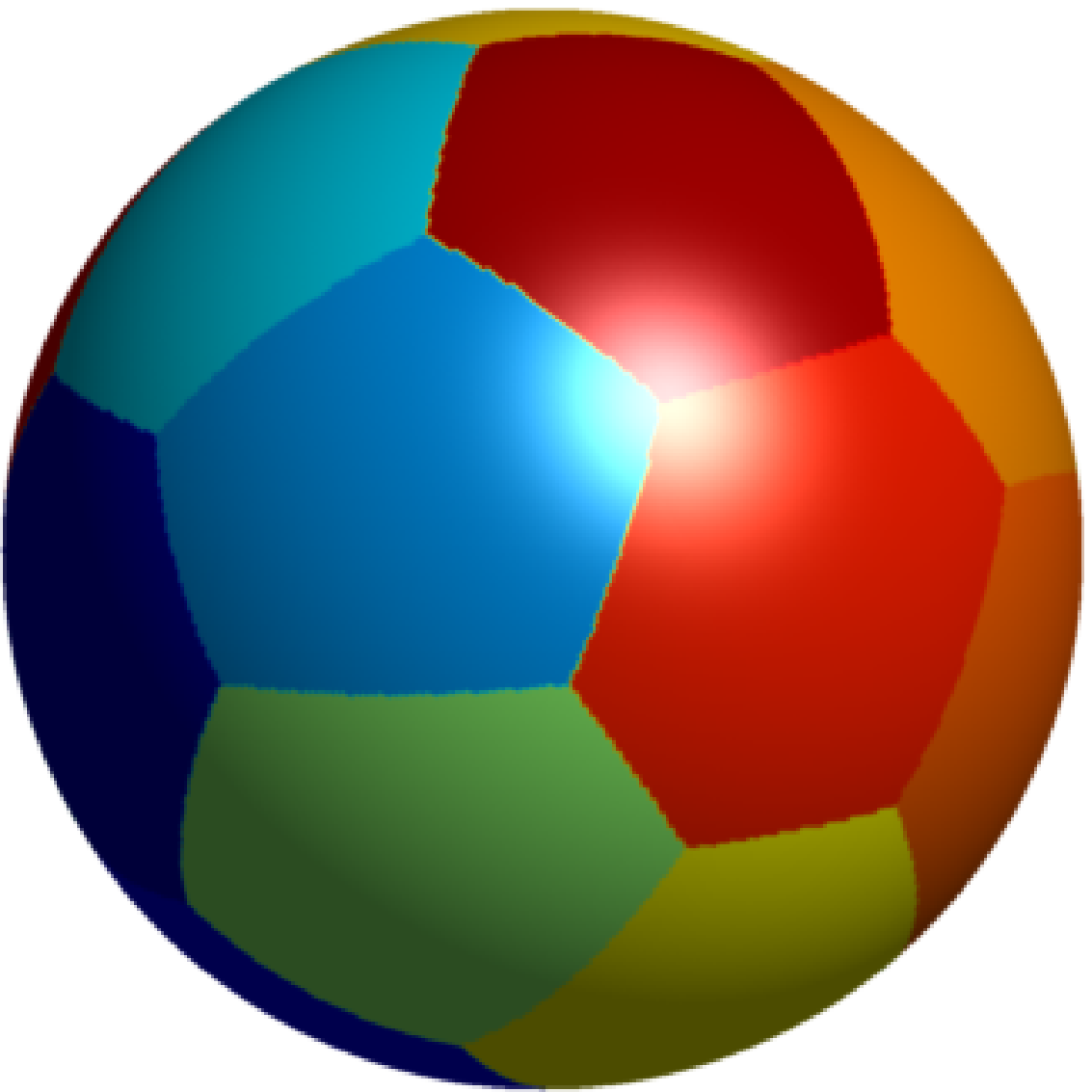}~
\includegraphics[width= 0.19\textwidth]{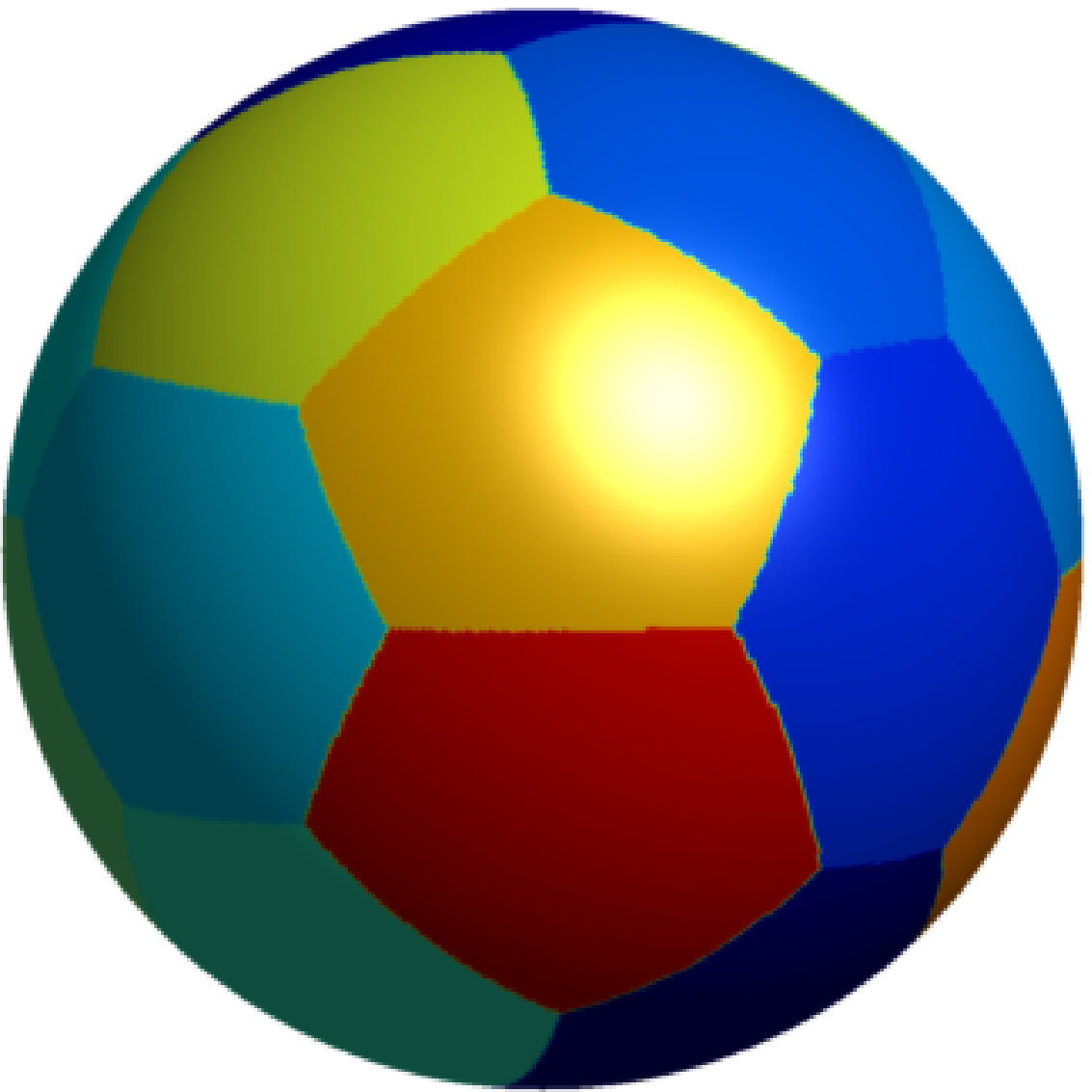}~
\includegraphics[width= 0.19\textwidth]{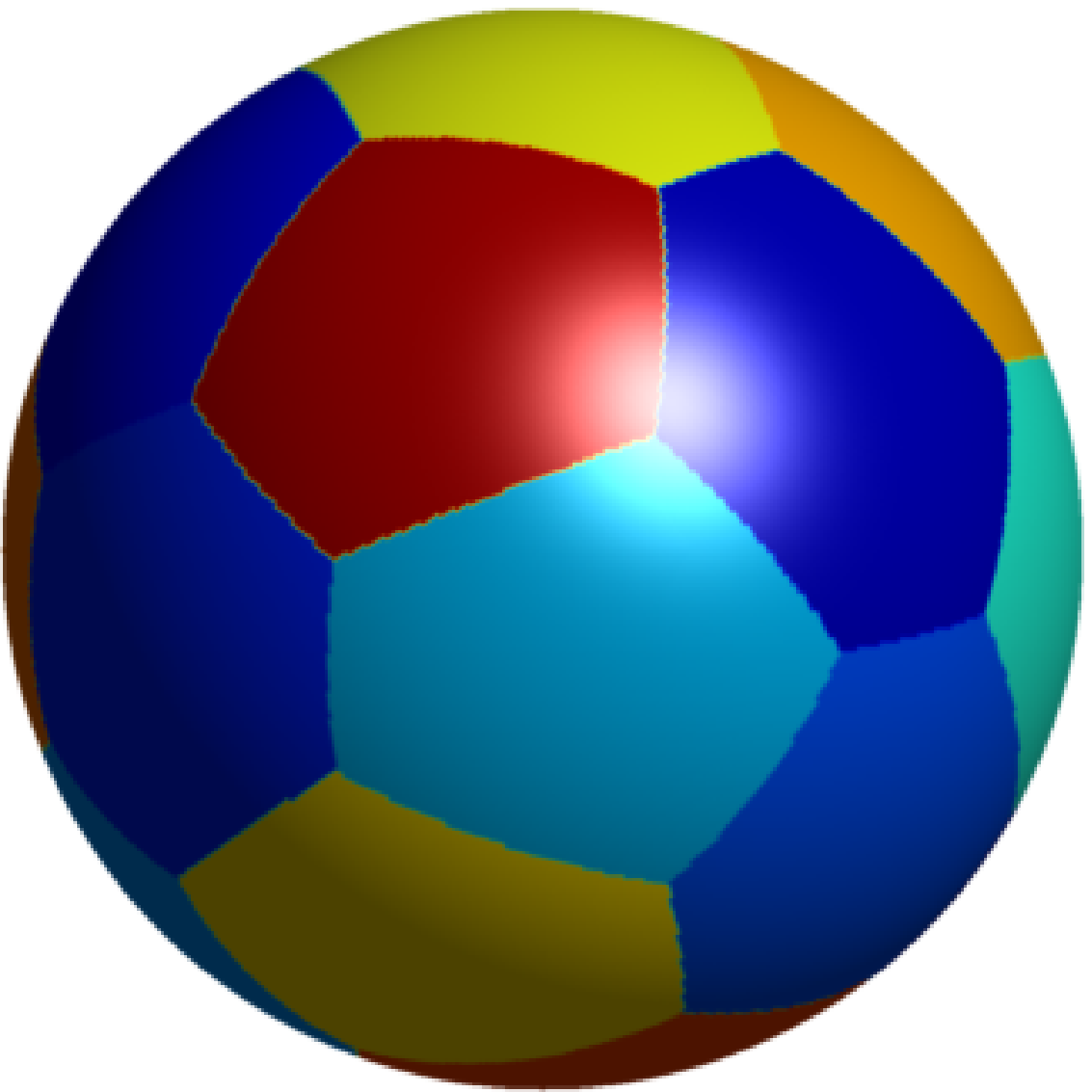}~
\vspace{0.1cm}

\includegraphics[width= 0.19\textwidth]{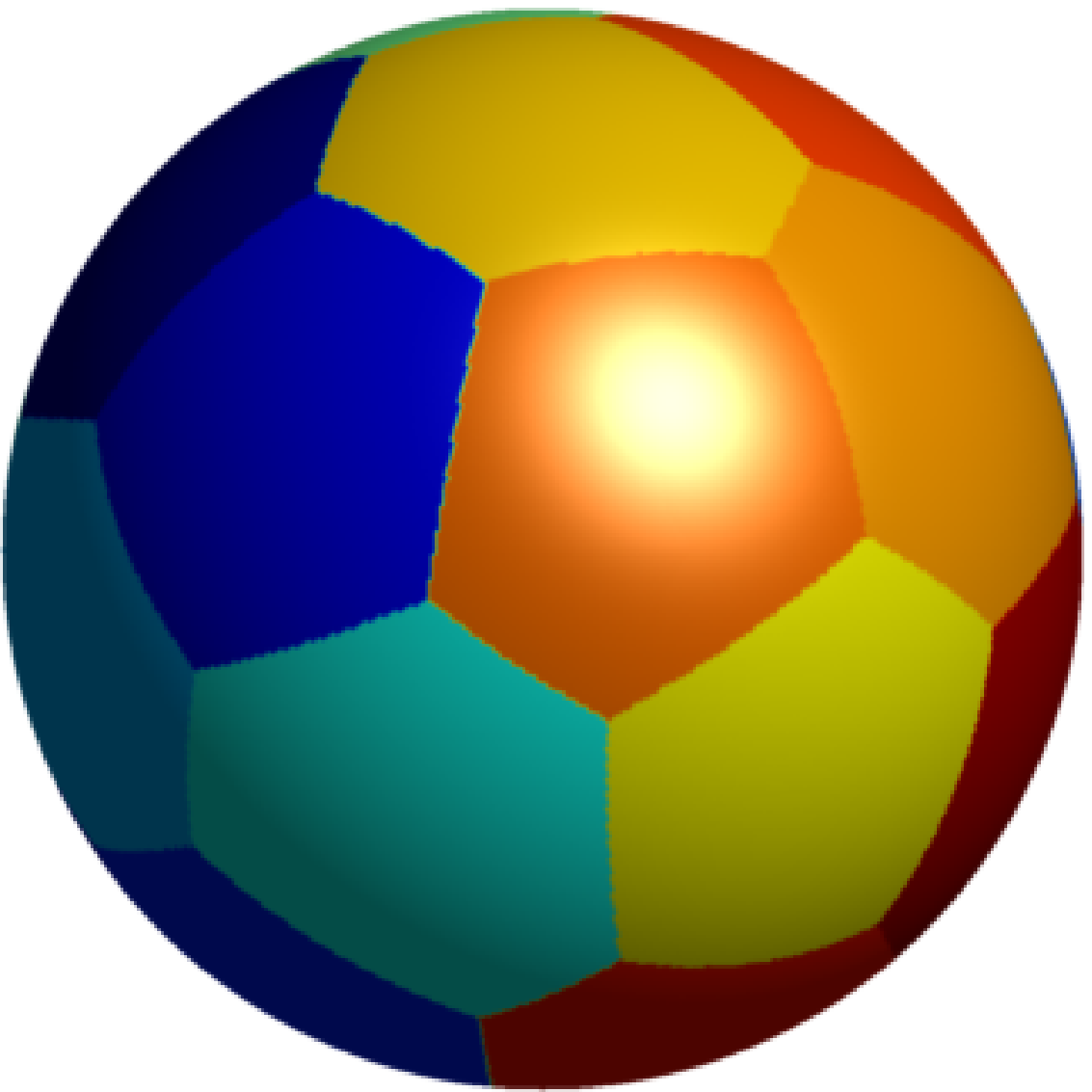}~
\includegraphics[width= 0.19\textwidth]{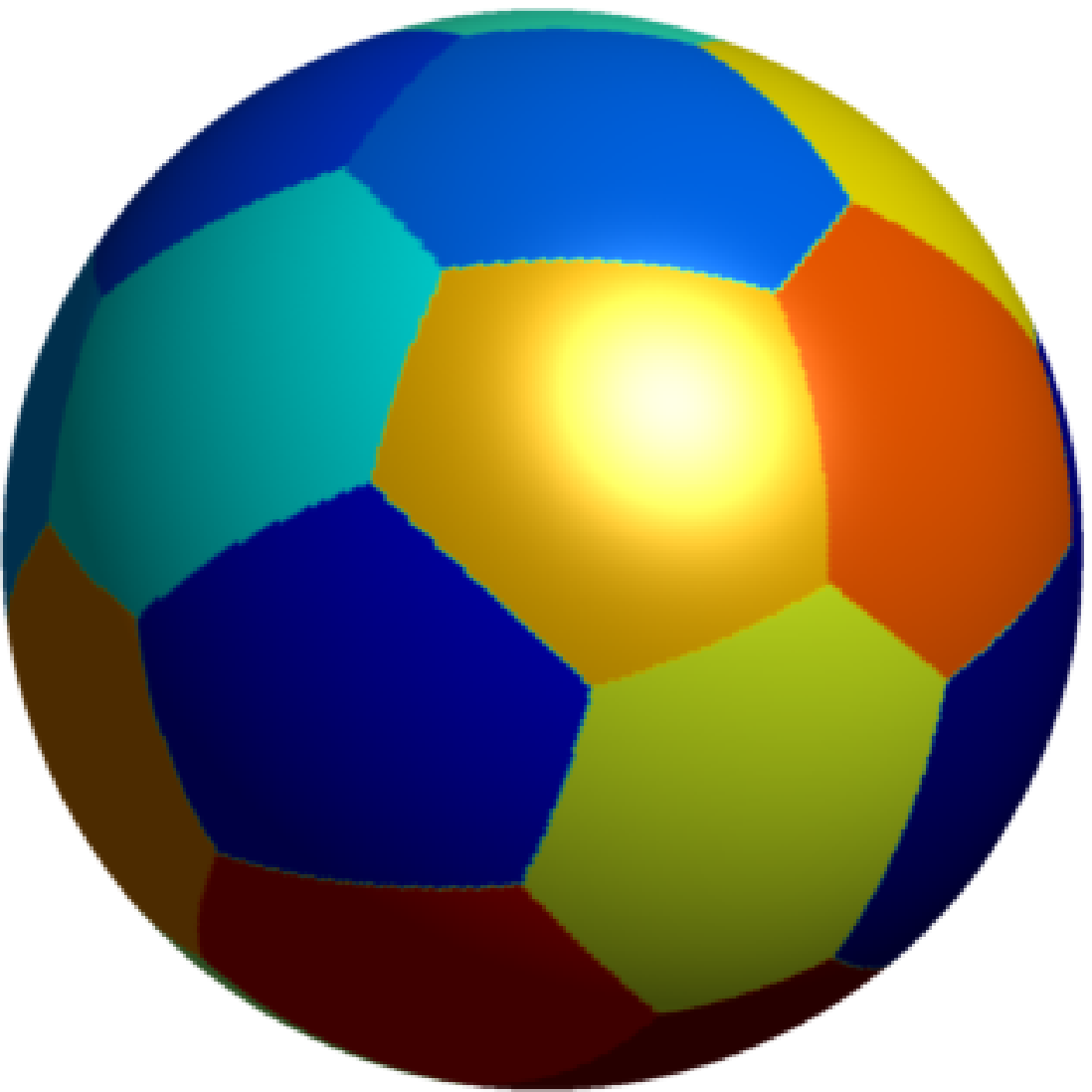}~
\includegraphics[width= 0.19\textwidth]{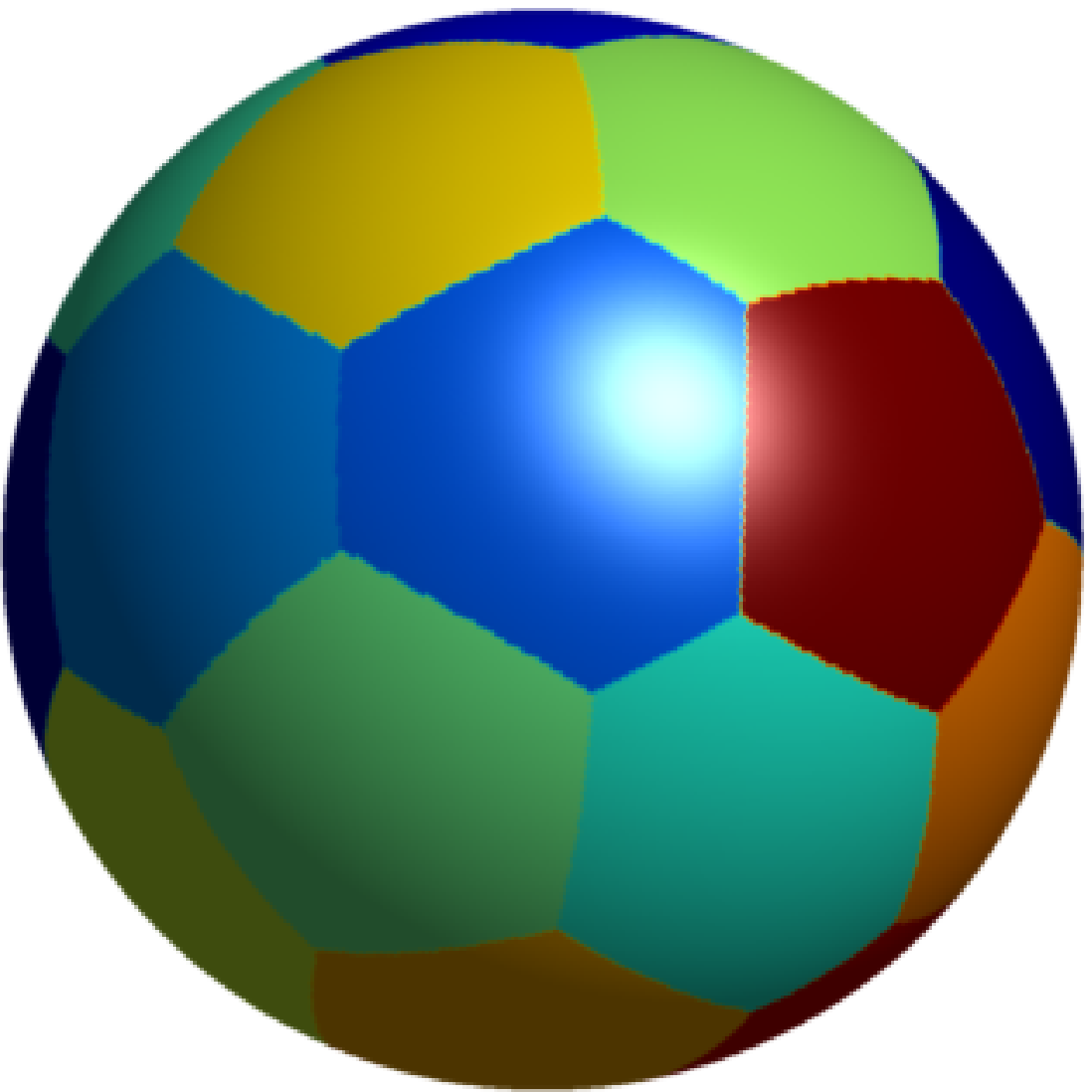}~
\includegraphics[width= 0.19\textwidth]{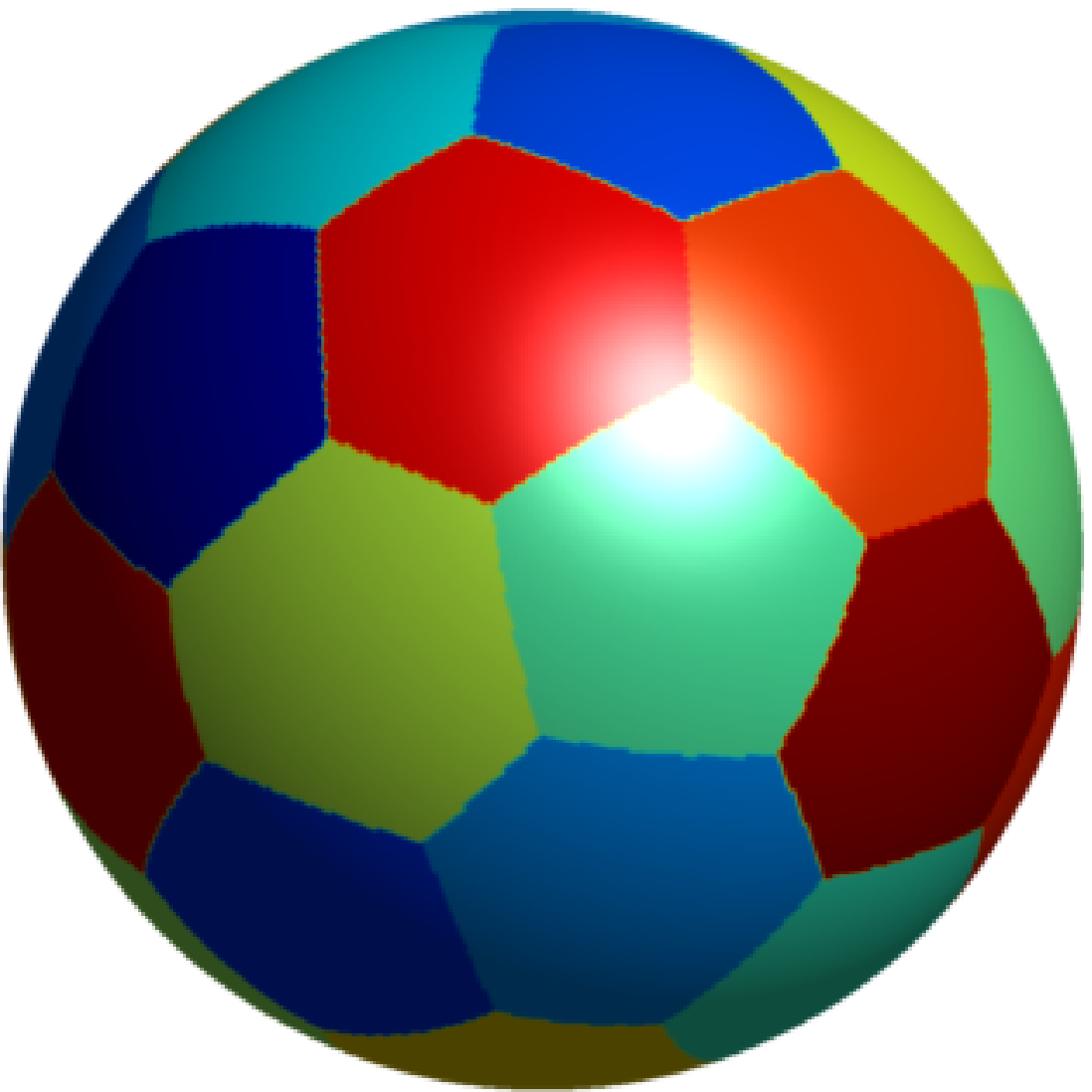}~

\caption{Minimal perimeter partitions on the sphere into $n$ equal area cells for $n \in \{2,3,...,24,32\}$.}
\label{sphere-perim}
\end{figure}

In order to have a good approximation of the optimal partition, we want do decrease $\varepsilon$ so that the width of the interface is small. We notice that if we chose $\varepsilon$ of the same order as the sides of the mesh triangles the algorithm converges. Furthermore, we cannot make $\varepsilon$ smaller, since then the gradient term will not contain any real information, as the width of the interface is of size $\varepsilon$. In order to avoid this problem, we consider refined meshes associated to each $\varepsilon$. At each step where we decrease $\varepsilon$ we interpolate the values of the previous optimizer on a refined mesh and we consider these interpolated densities as starting point for the descent algorithm on the new mesh. In the case of the sphere we make four refinements ranging from $10000$ to $160000$ points. Some optimal configurations, in the case of the sphere, are presented in Figure \ref{sphere-perim}. A detailed study of the case of the sphere along with a comparison with the known results of Cox and Flikkema \cite{cox-partitions} are presented in the next section.

\begin{figure}
\centering
\includegraphics[width= 0.19\textwidth]{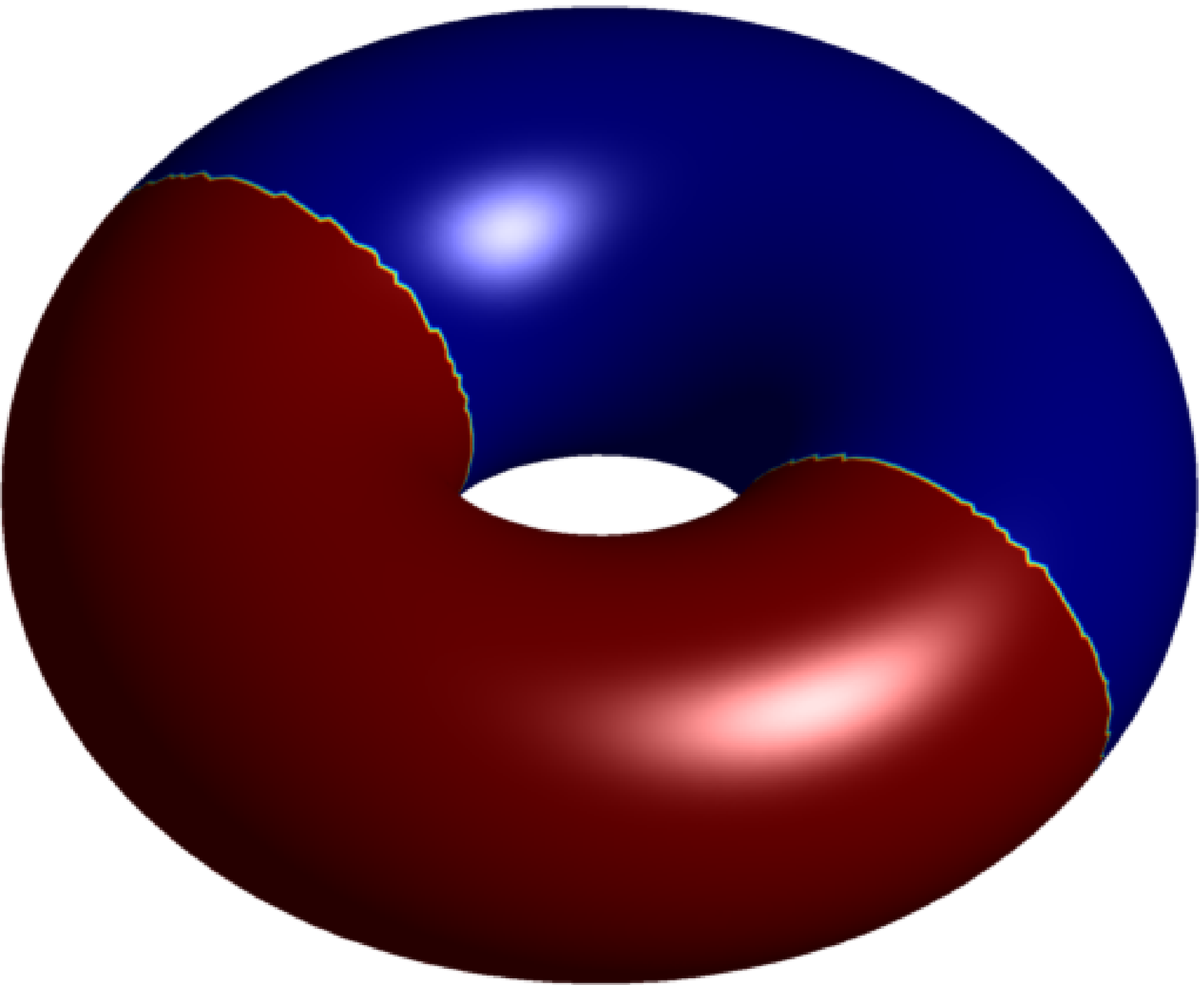}~
\includegraphics[width= 0.19\textwidth]{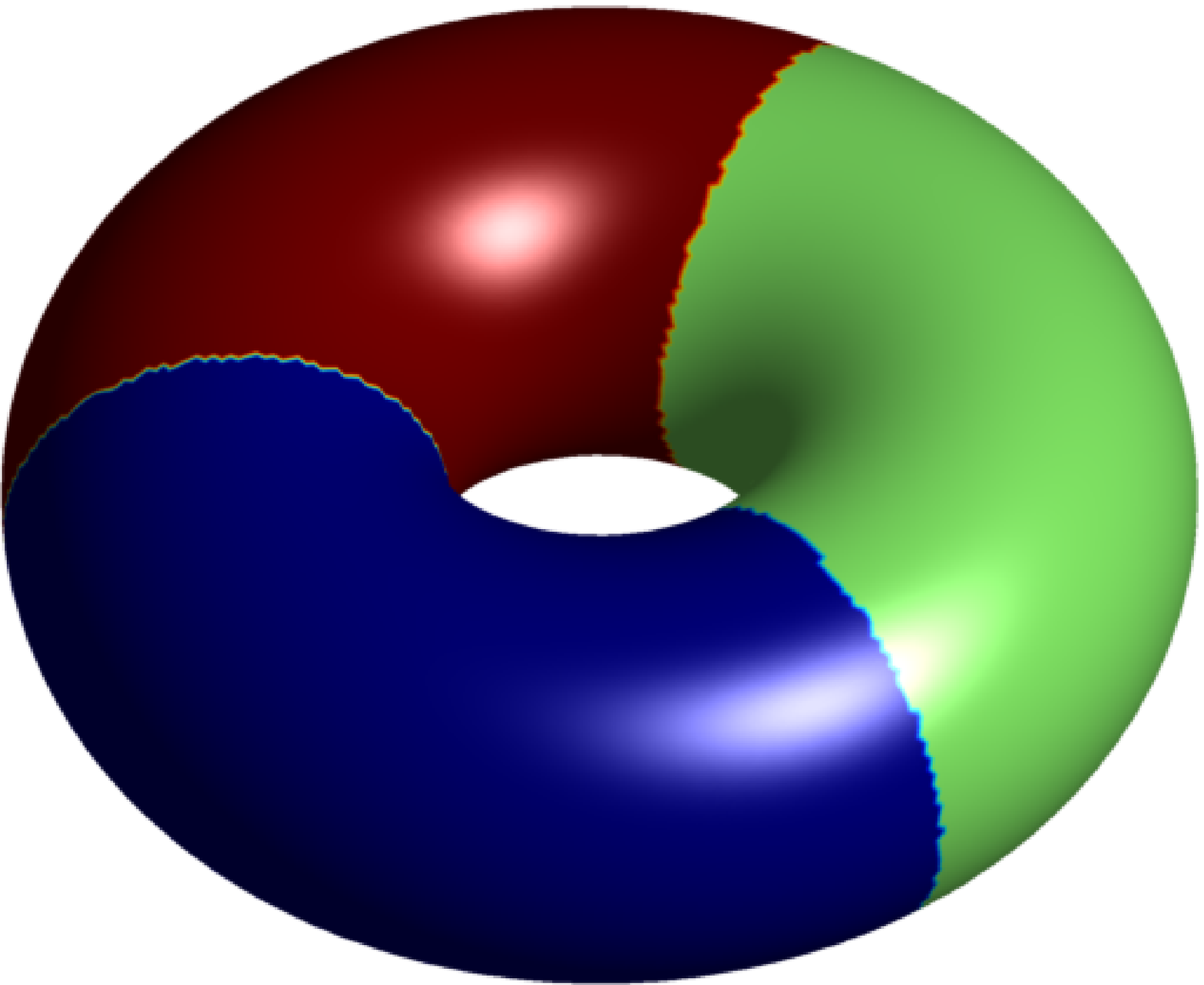}~
\includegraphics[width= 0.19\textwidth]{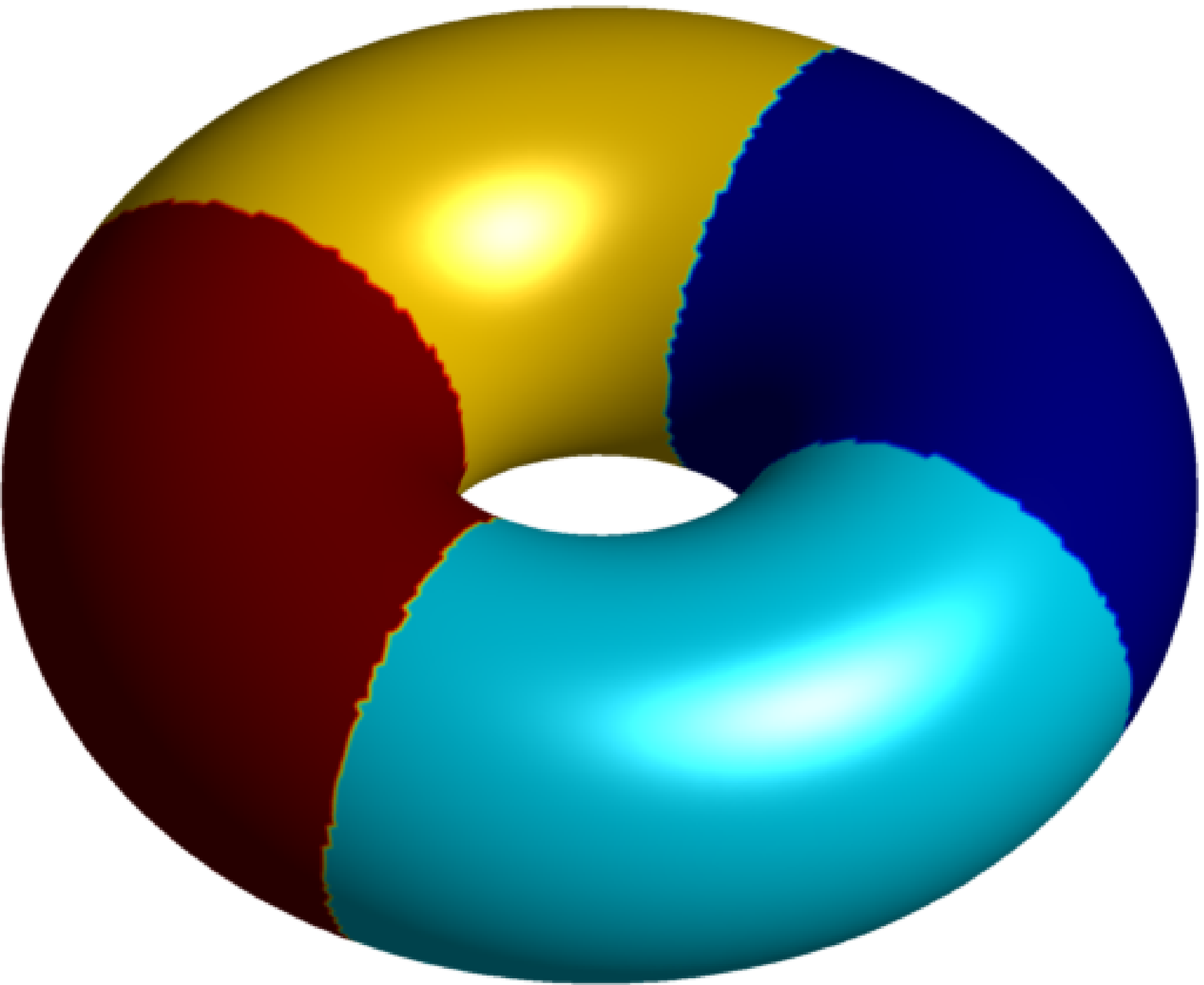}~
\includegraphics[width= 0.19\textwidth]{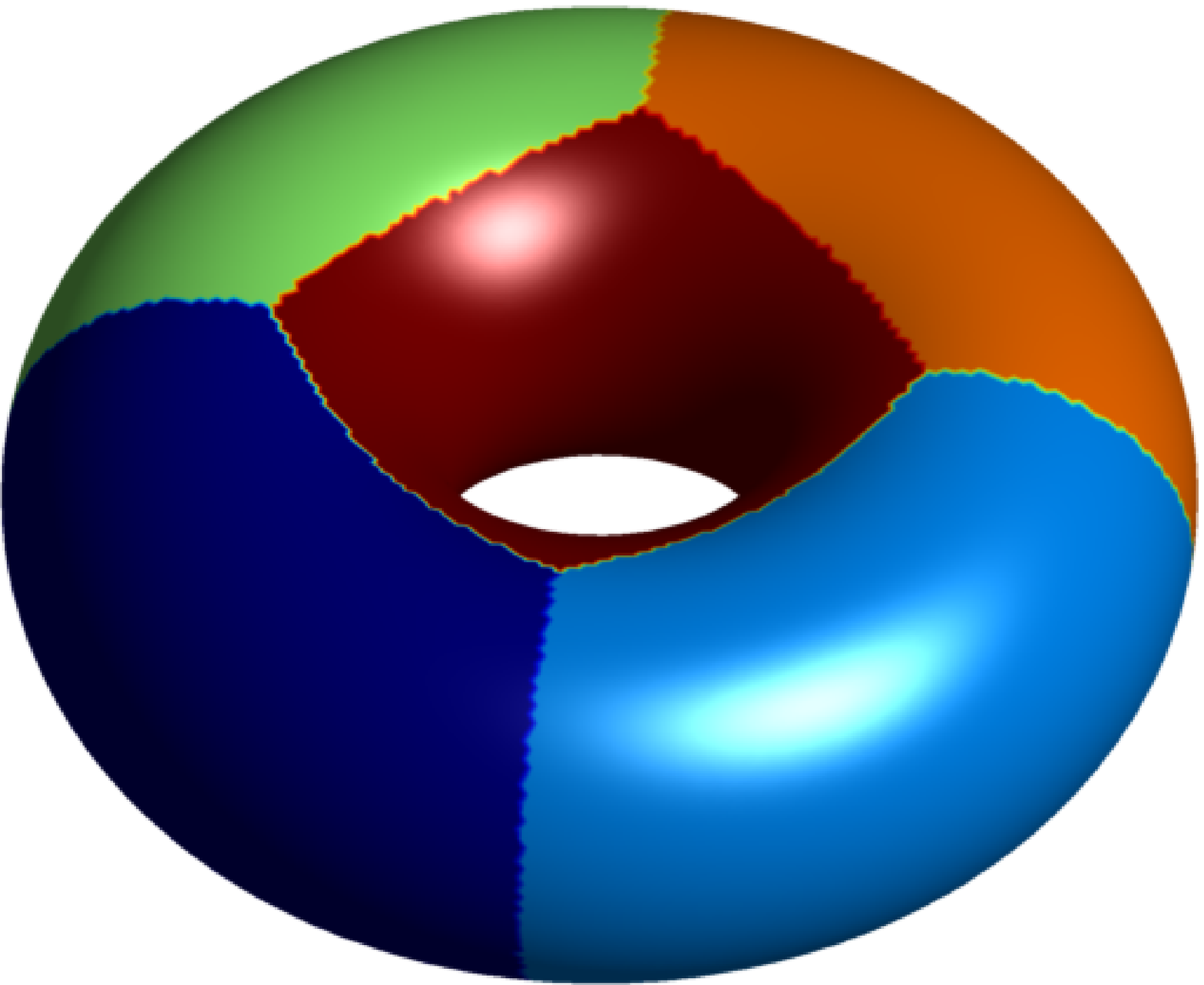}~
\includegraphics[width= 0.19\textwidth]{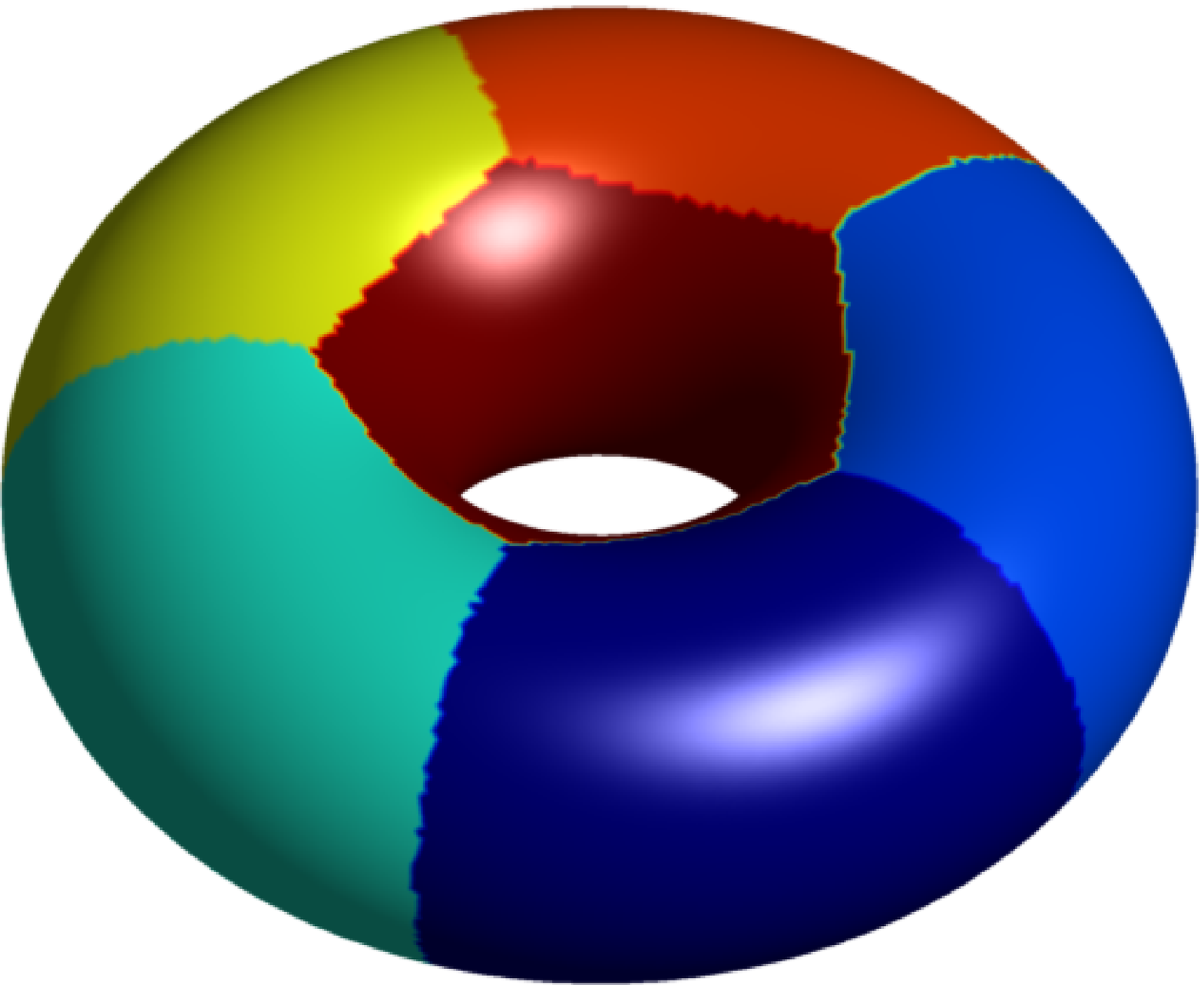}

\includegraphics[width= 0.19\textwidth]{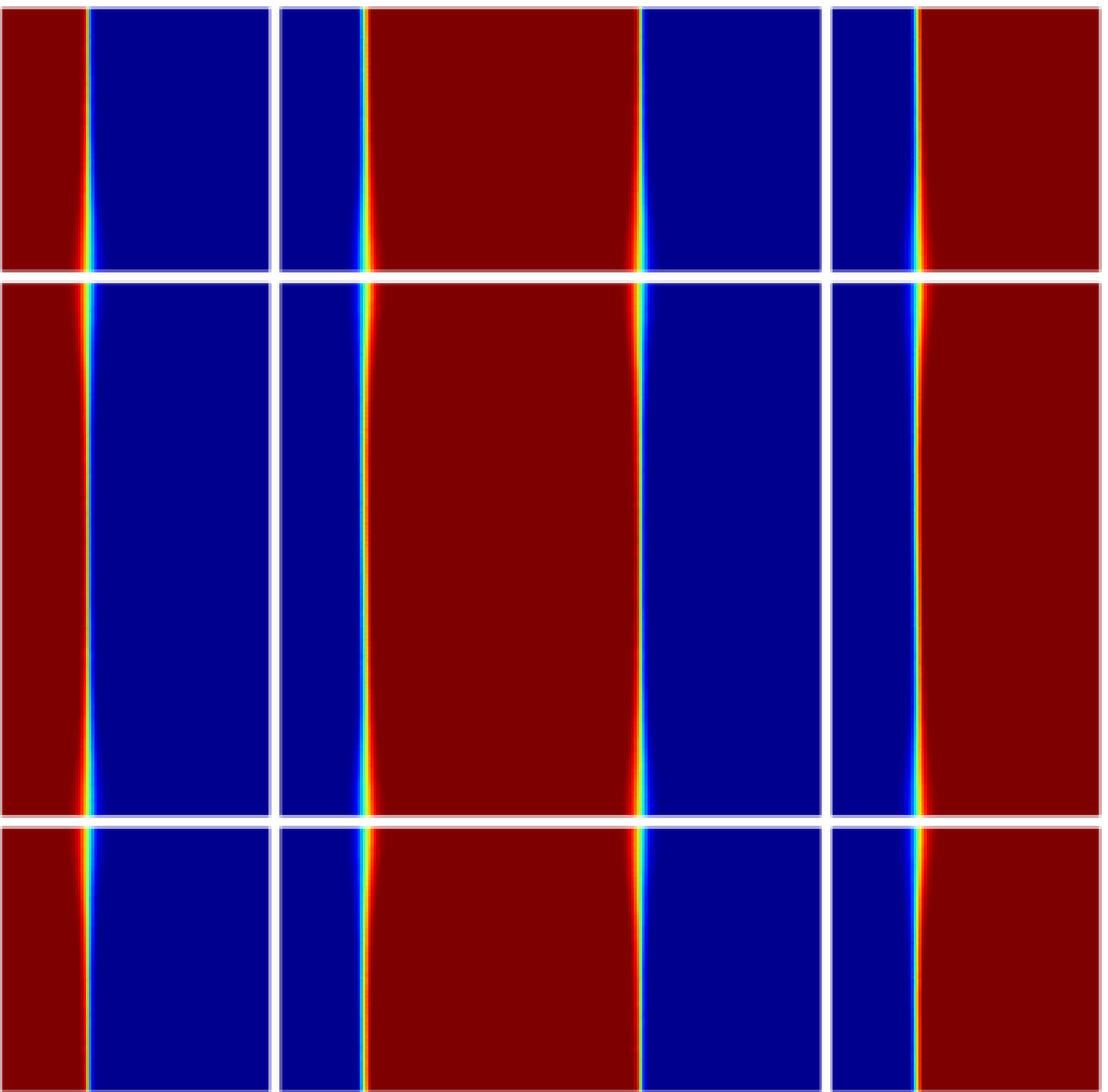}~
\includegraphics[width= 0.19\textwidth]{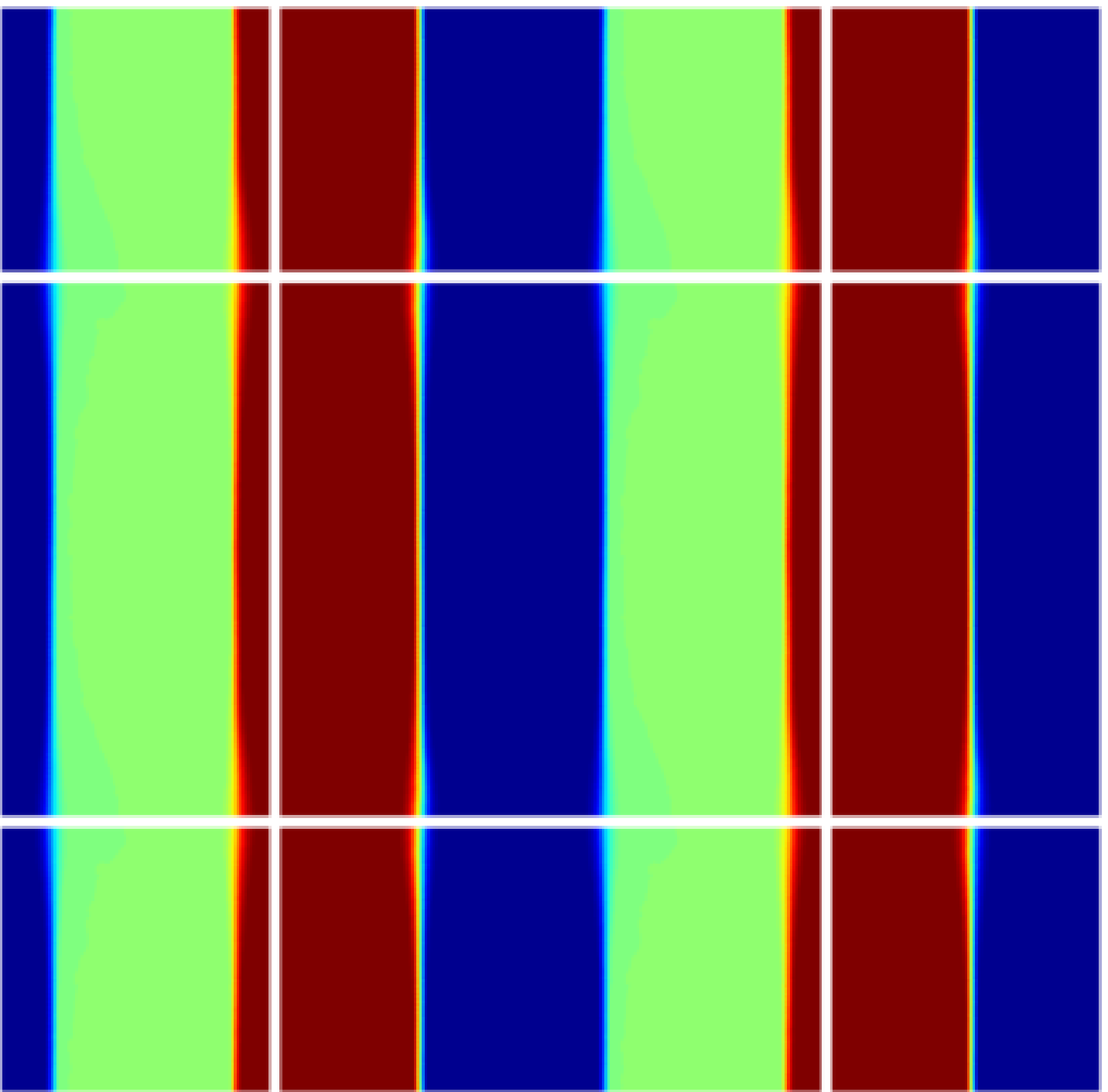}~
\includegraphics[width= 0.19\textwidth]{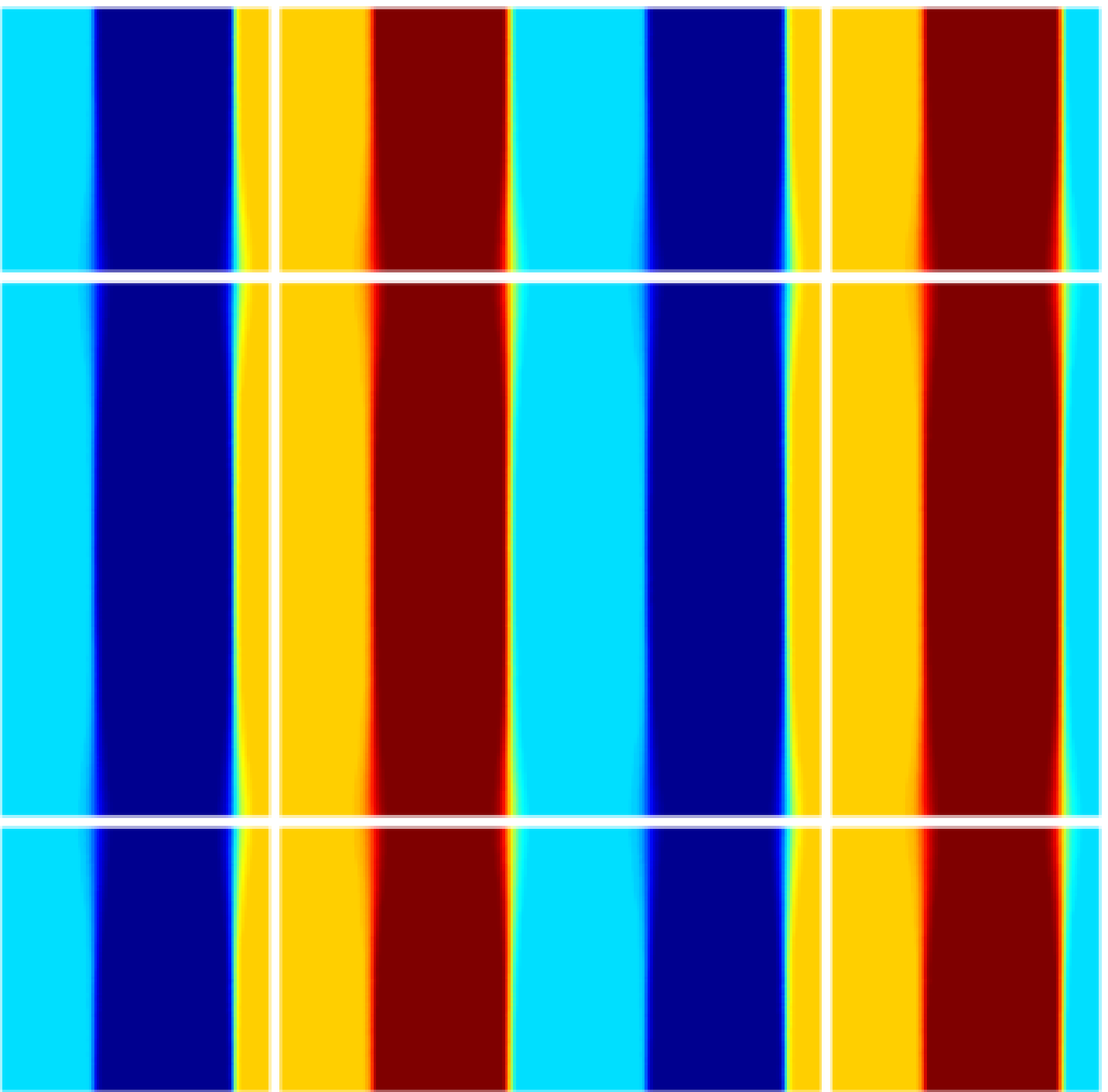}~
\includegraphics[width= 0.19\textwidth]{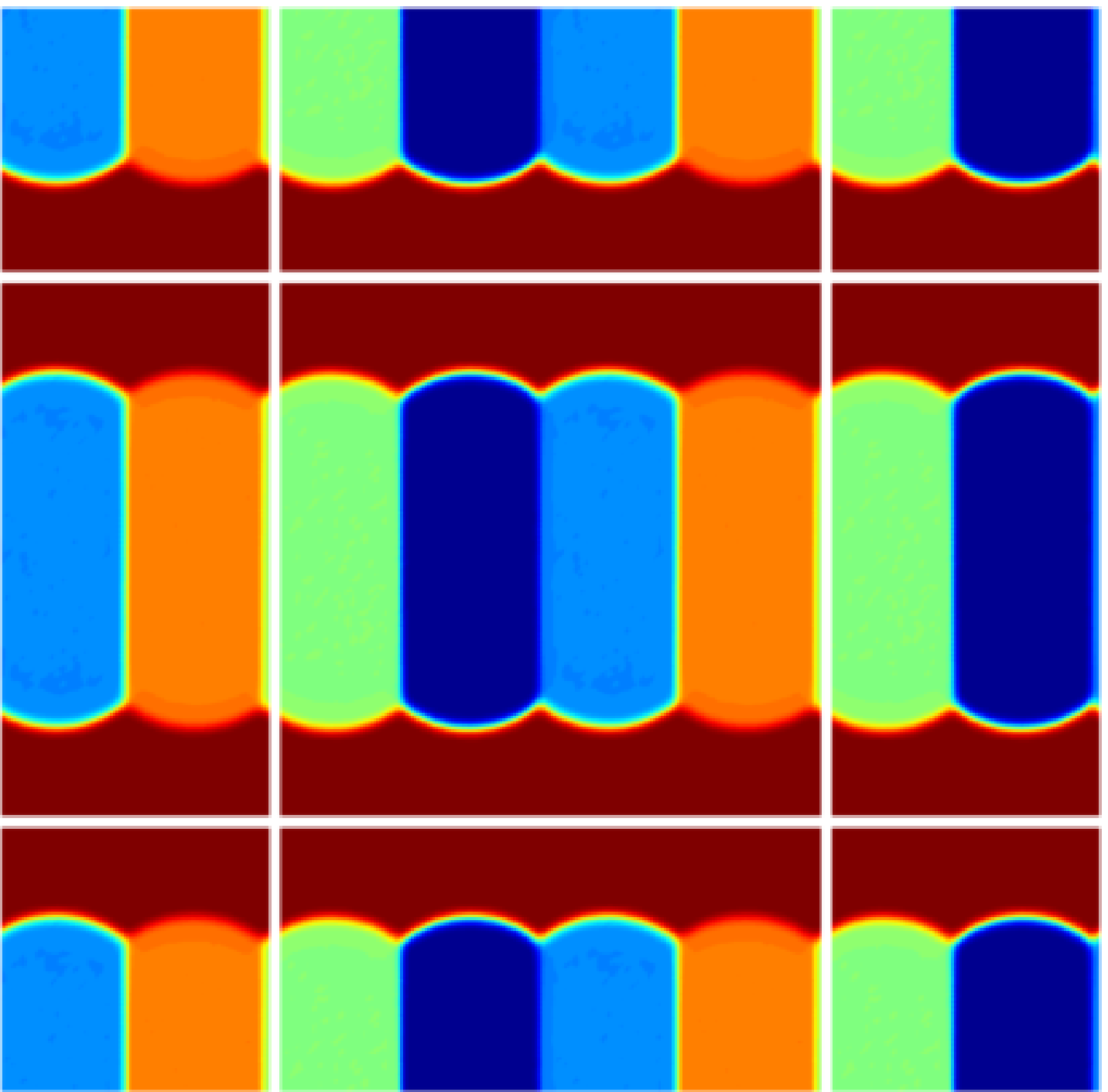}~
\includegraphics[width= 0.19\textwidth]{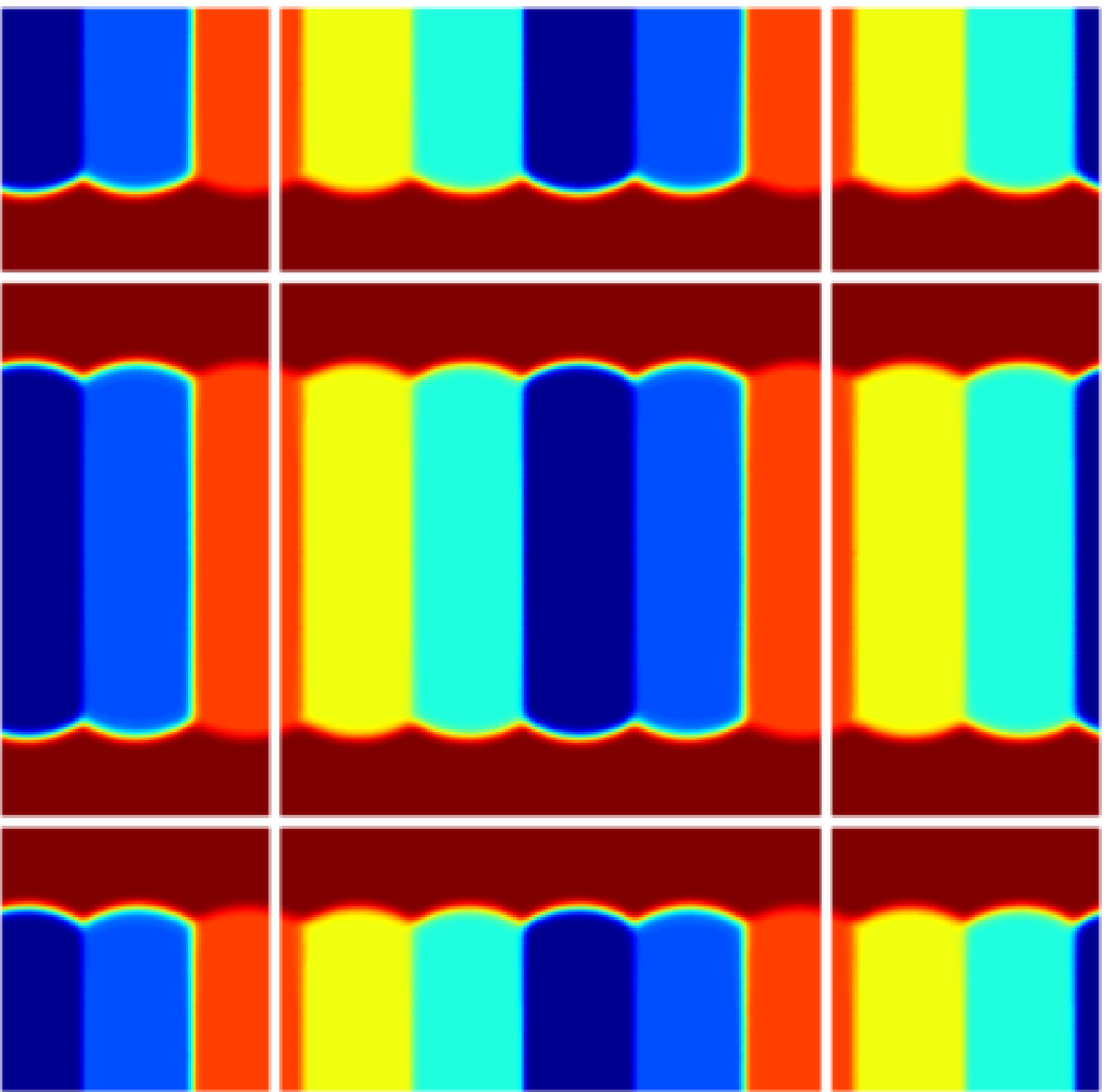}

\includegraphics[width= 0.19\textwidth]{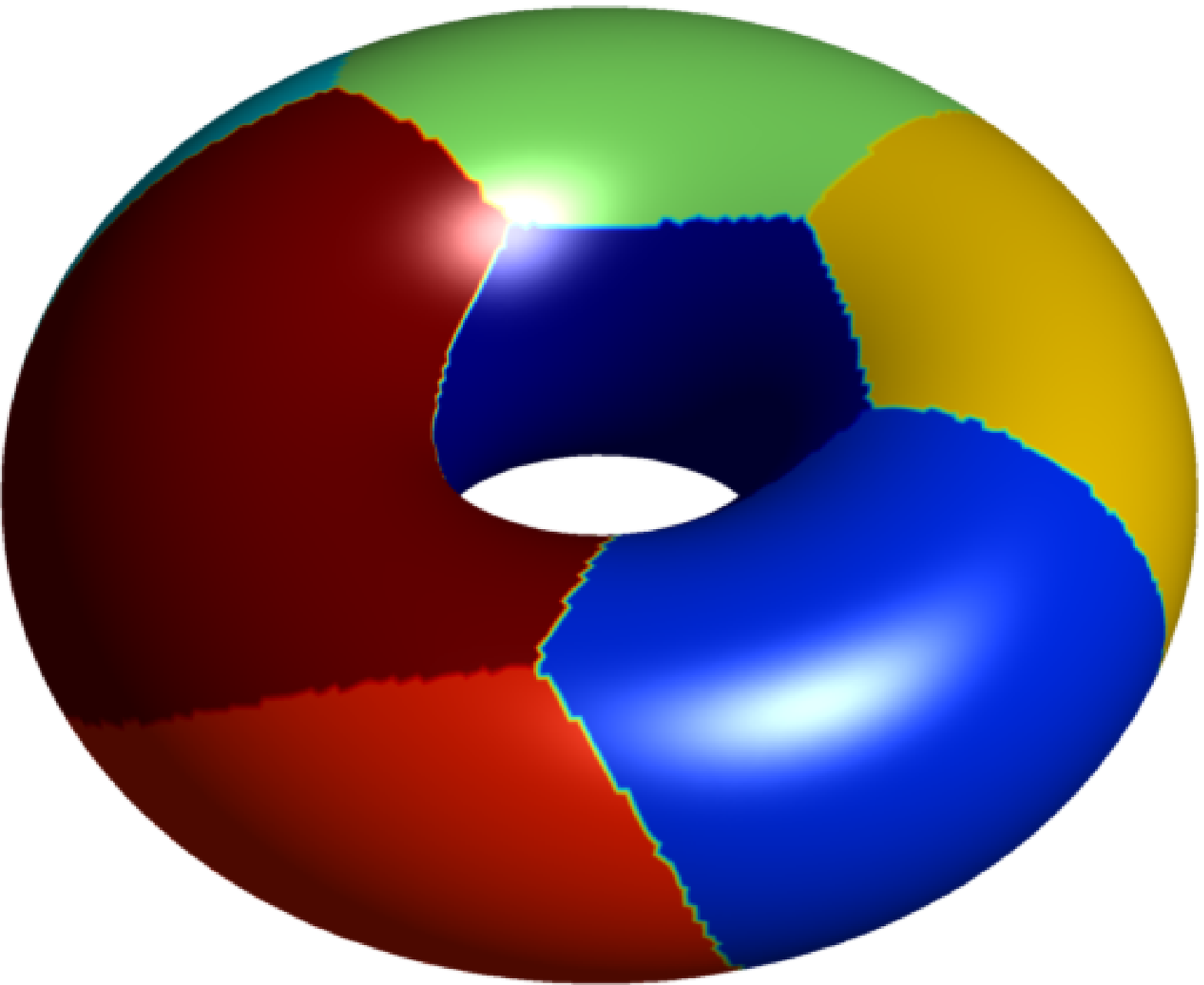}~
\includegraphics[width= 0.19\textwidth]{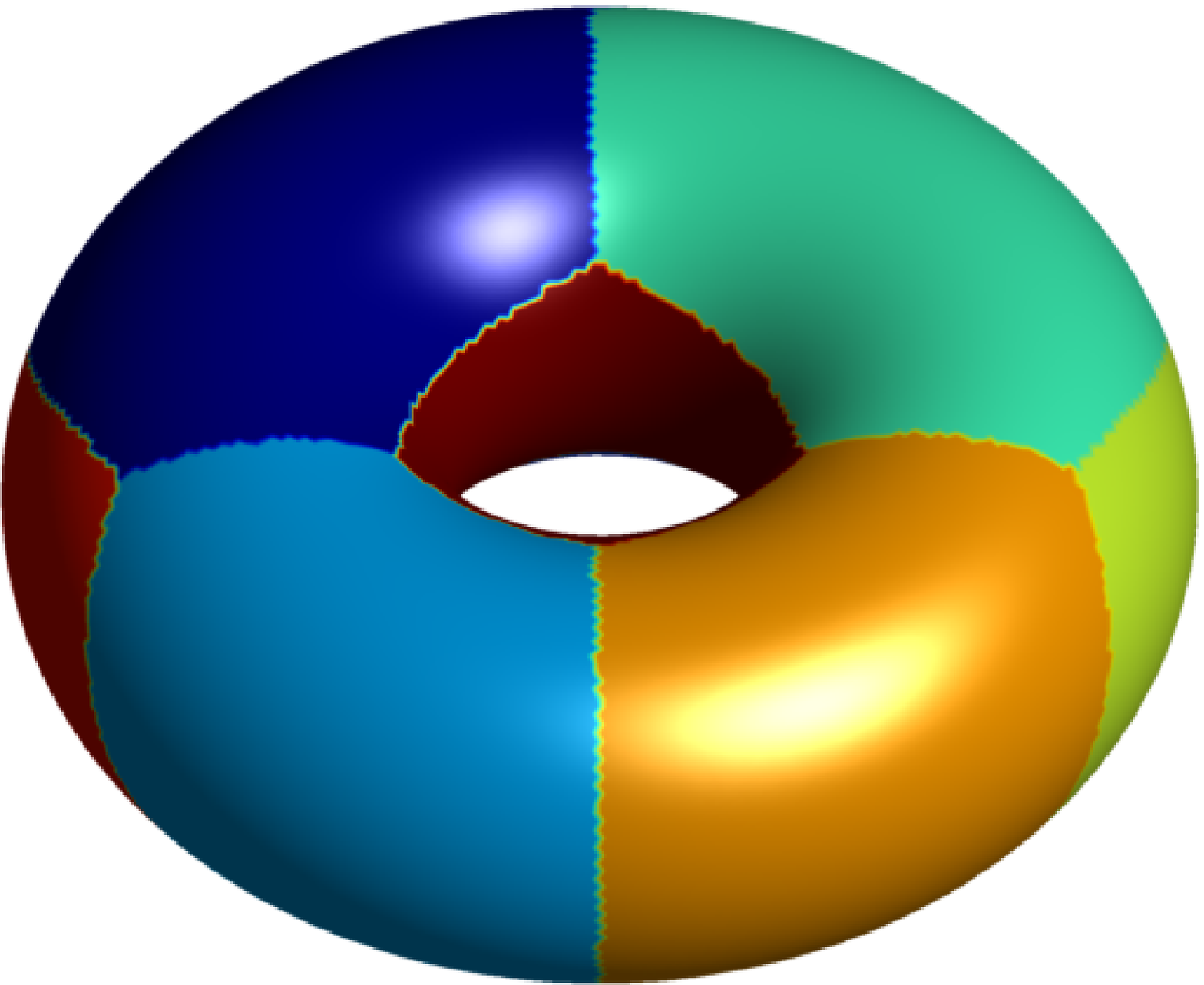}~
\includegraphics[width= 0.19\textwidth]{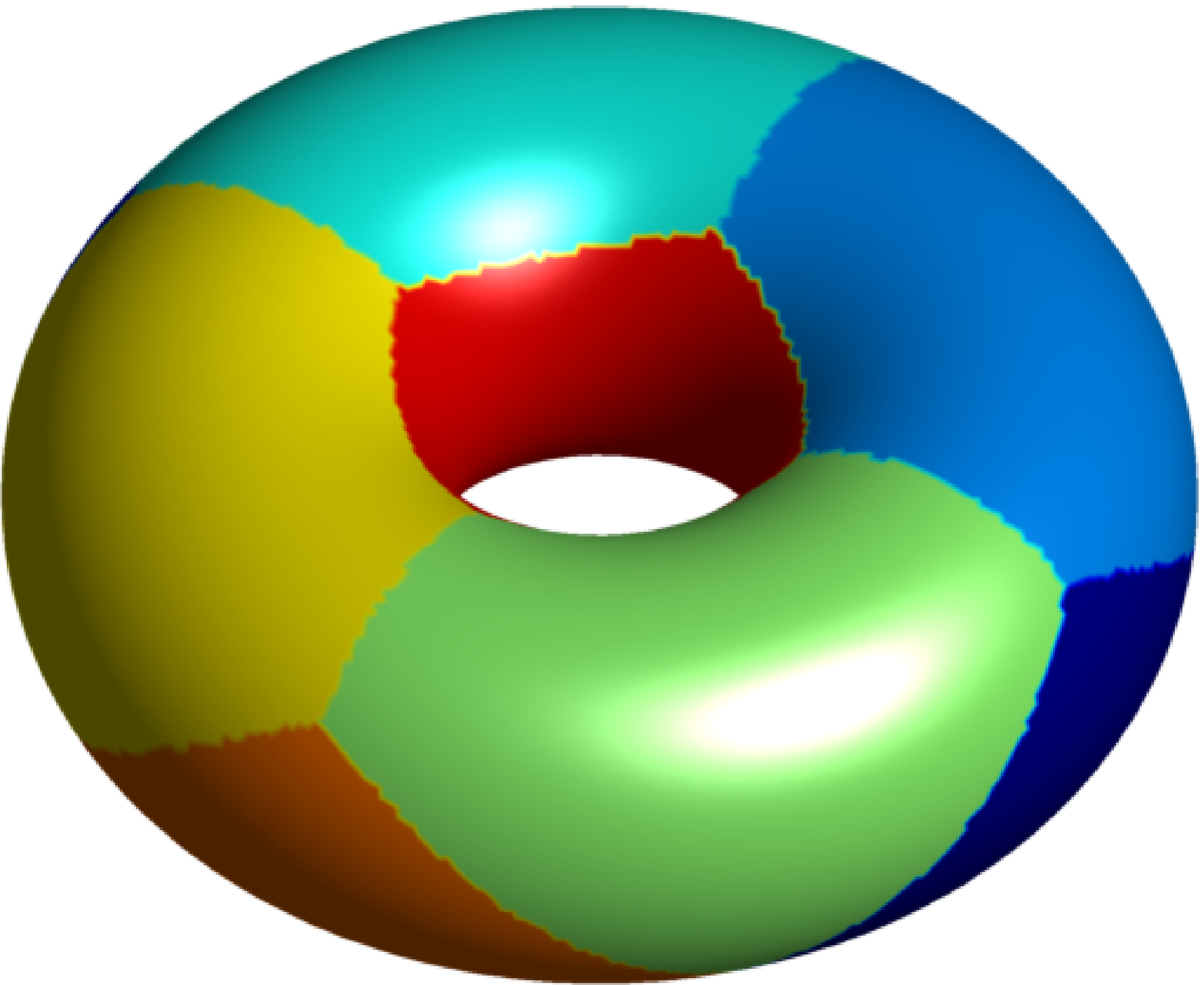}~
\includegraphics[width= 0.19\textwidth]{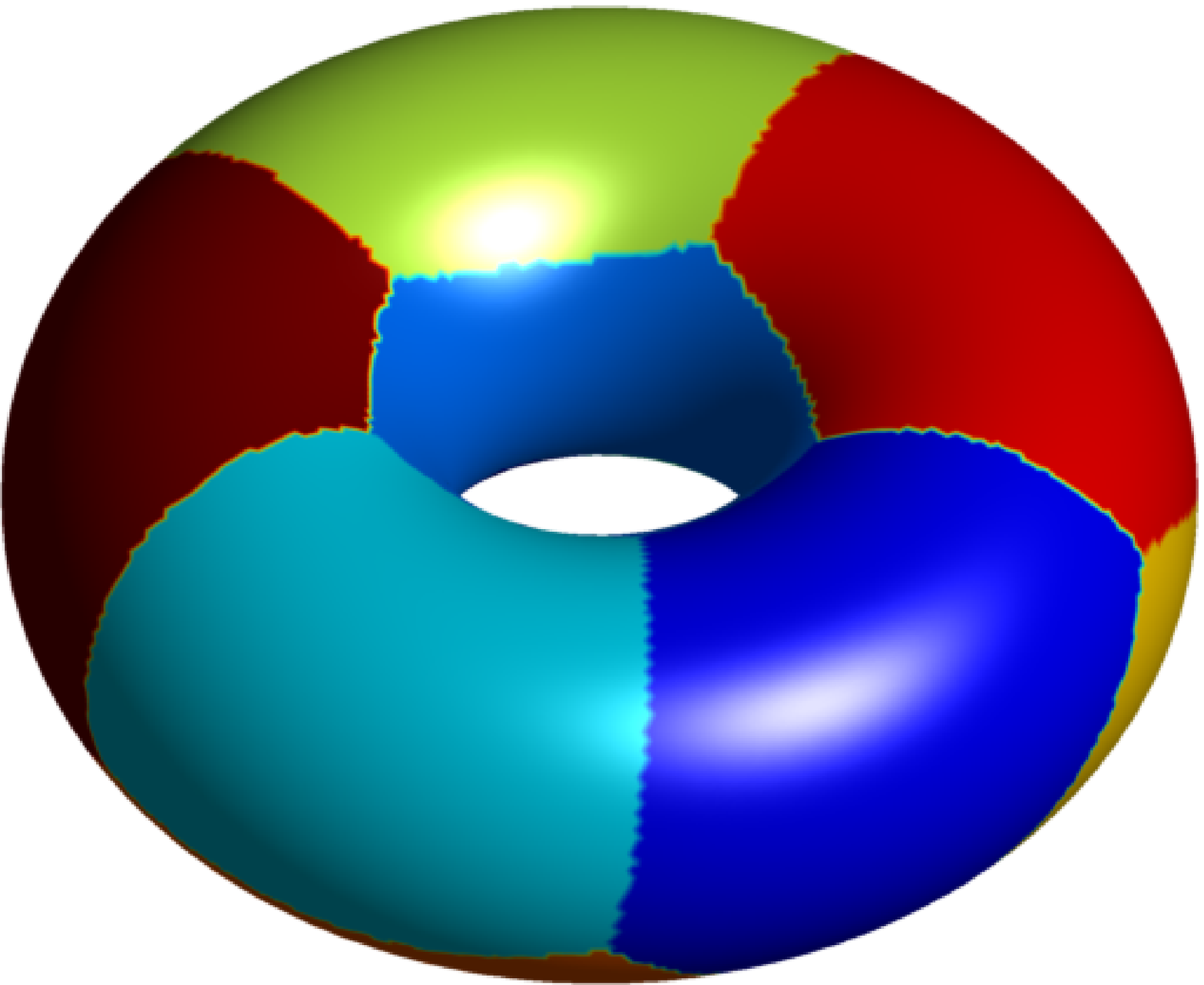}~
\includegraphics[width= 0.19\textwidth]{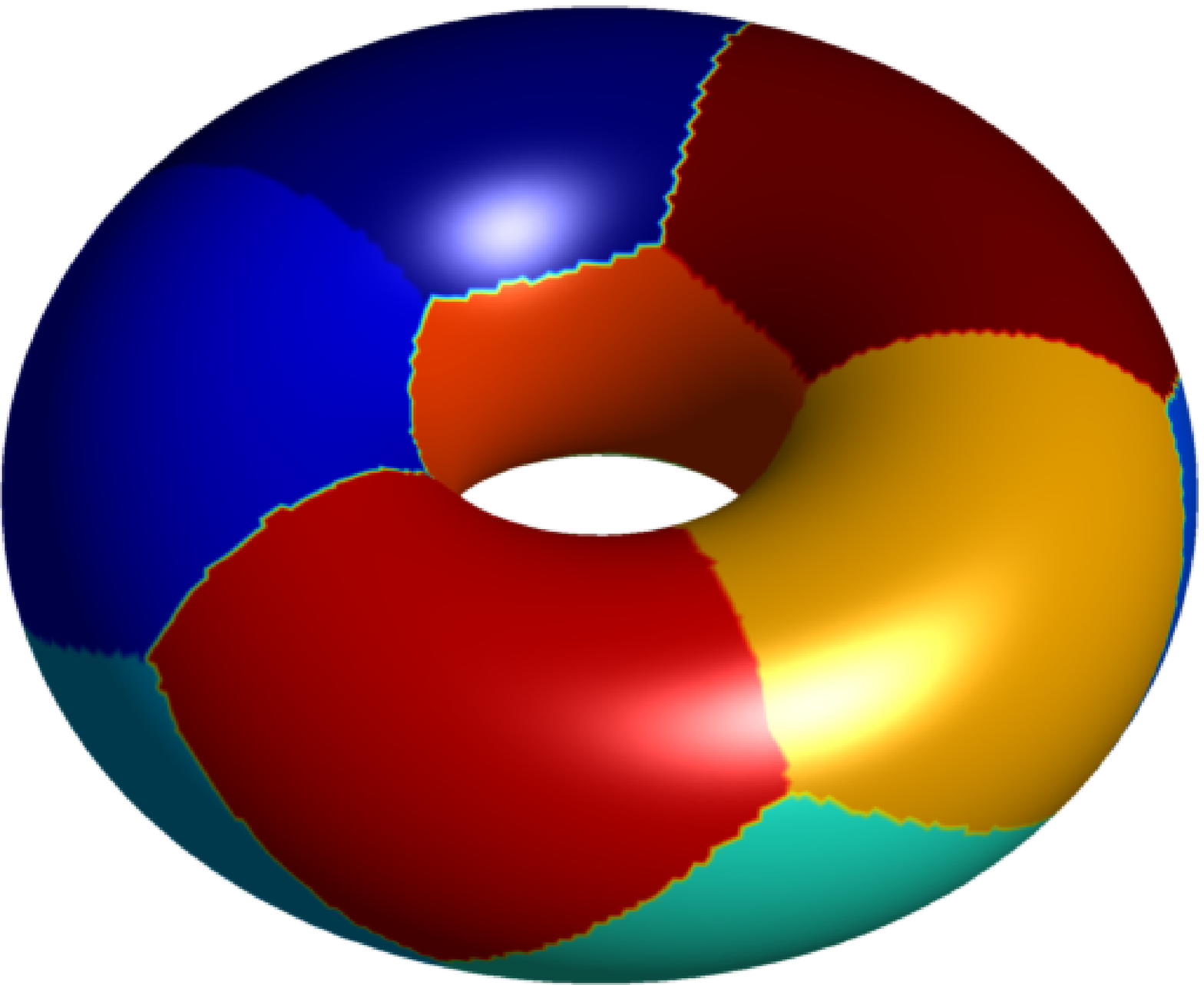}

\includegraphics[width= 0.19\textwidth]{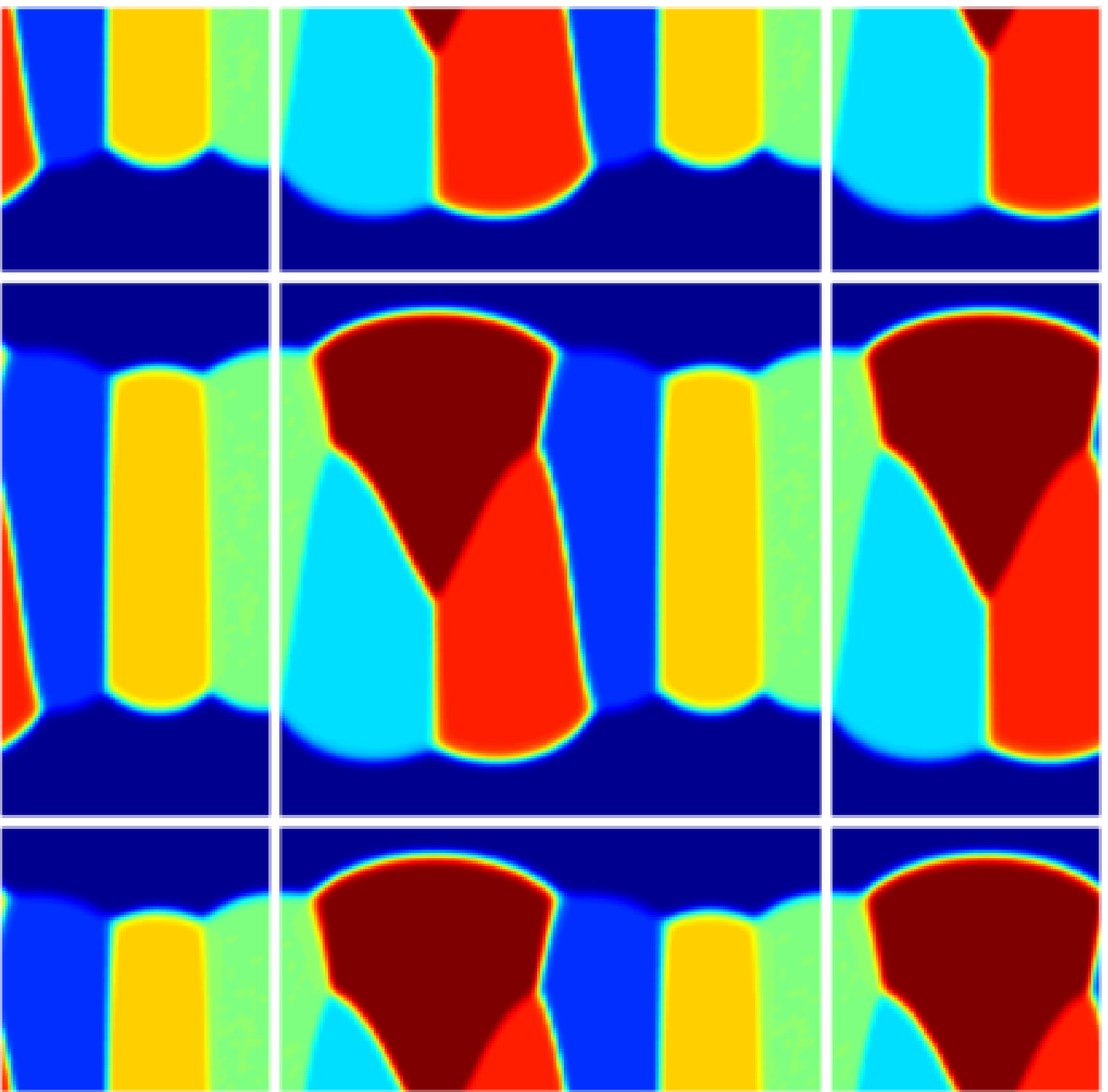}~
\includegraphics[width= 0.19\textwidth]{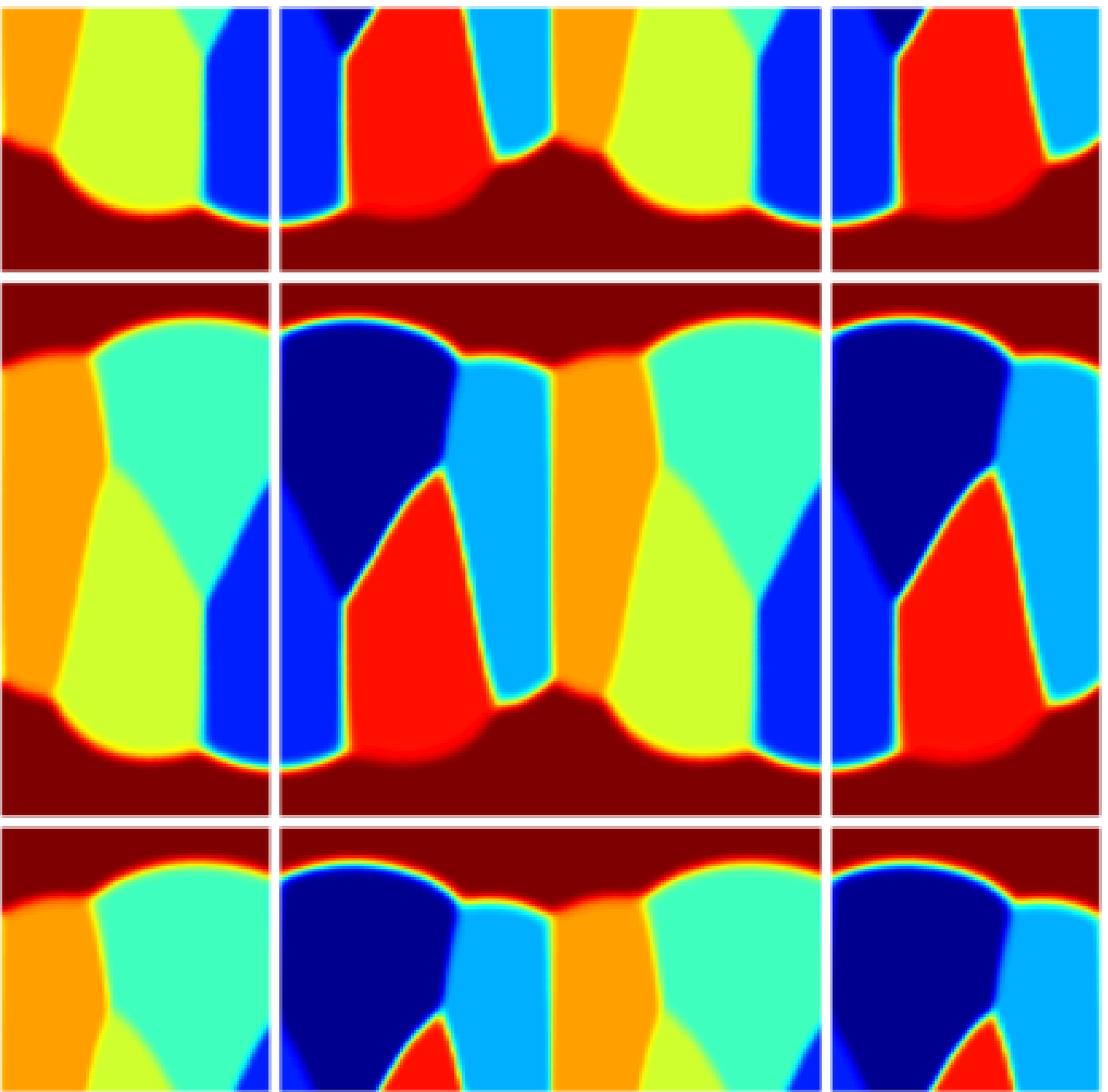}~
\includegraphics[width= 0.19\textwidth]{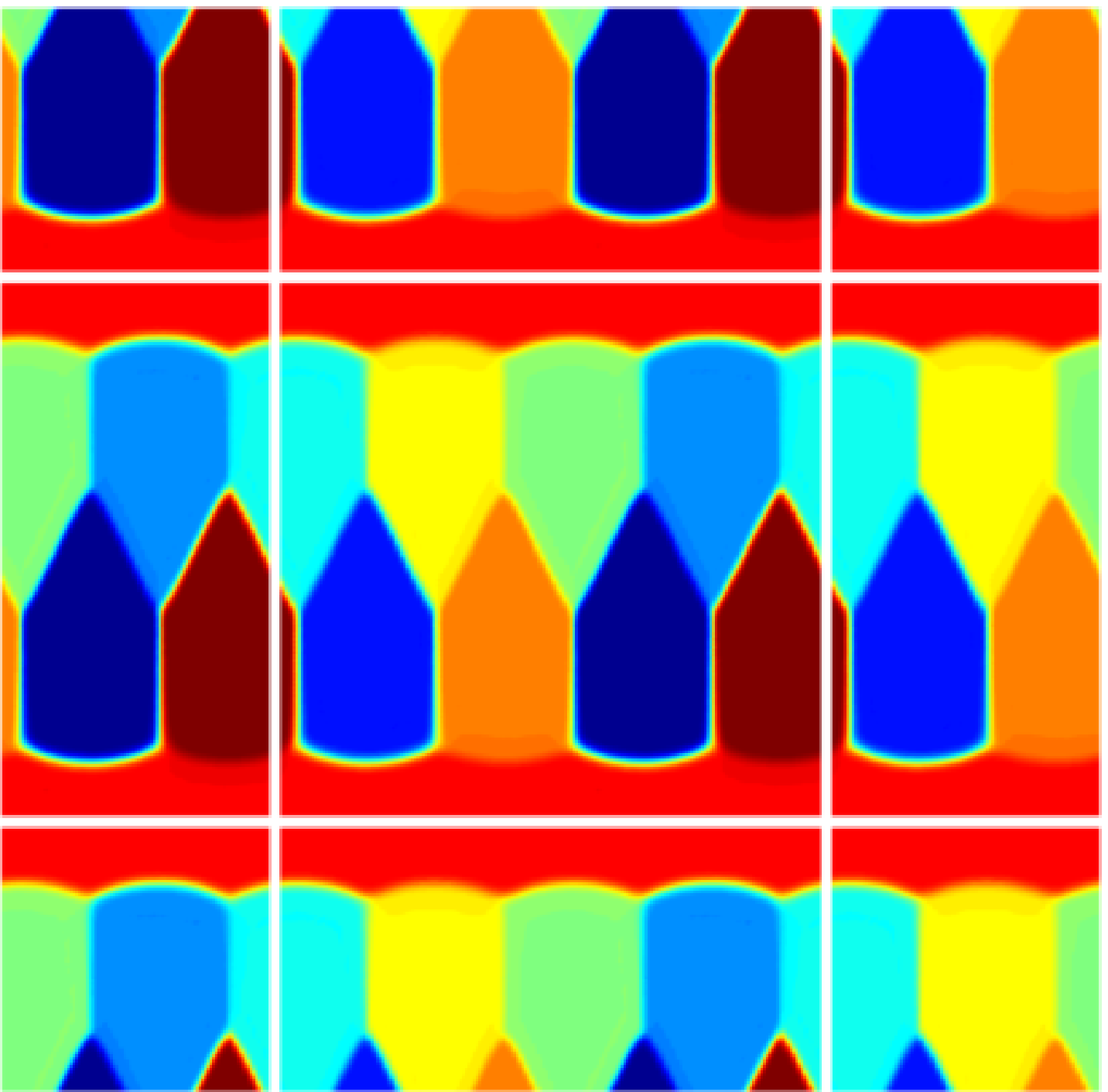}~
\includegraphics[width= 0.19\textwidth]{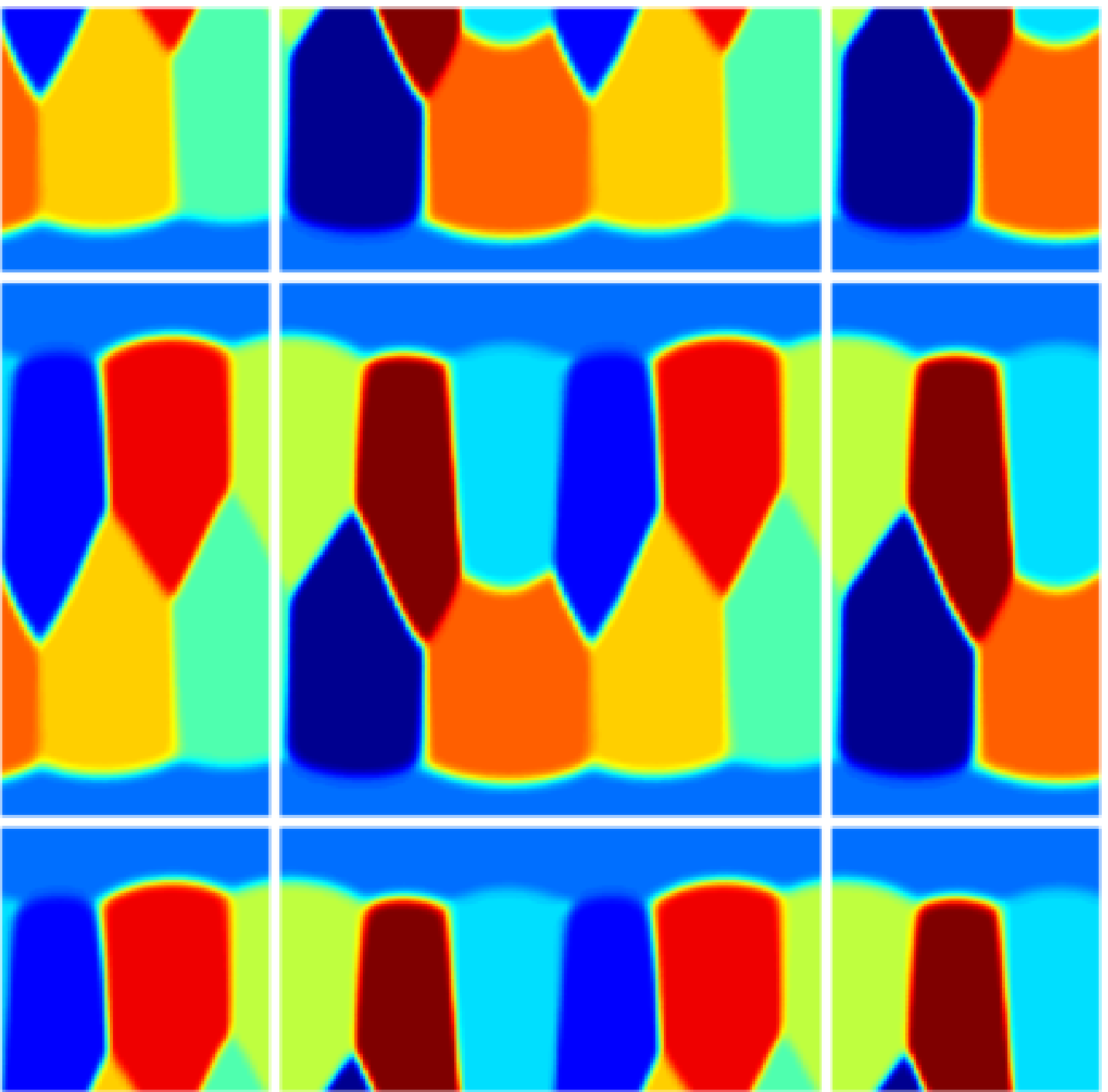}~
\includegraphics[width= 0.19\textwidth]{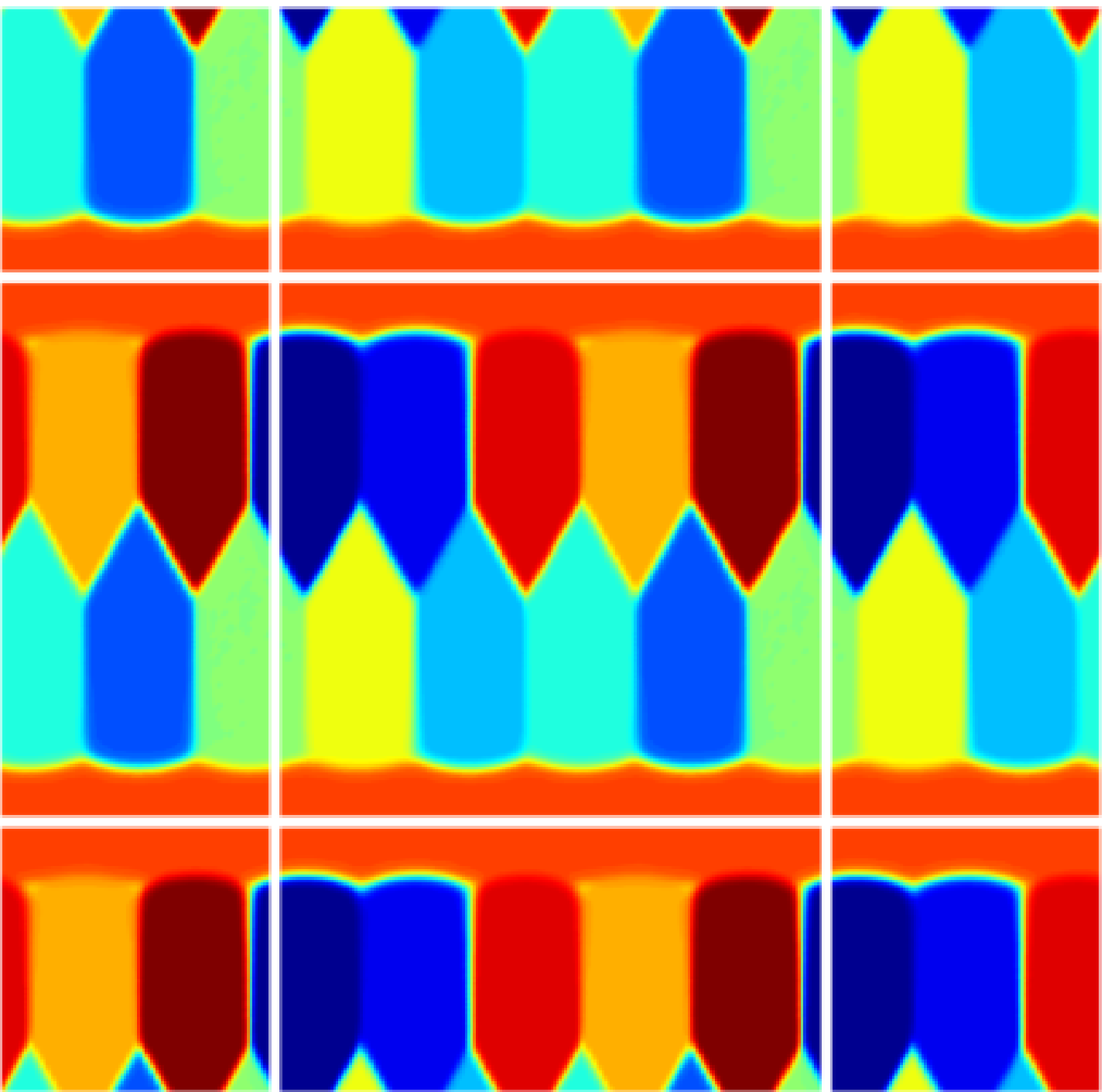}
\caption{Minimal perimeter partitions on the torus with outer radius $R=1$ and inner radius $r = 0.6$ together with their associated flattenings for $n \in [2,11]$. The center rectangle is represents the torus, while periodic continuations are made to easily see the topological structure.}
\label{torus-perim}
\end{figure}

\begin{figure}
\centering
\includegraphics[width= 0.29\textwidth]{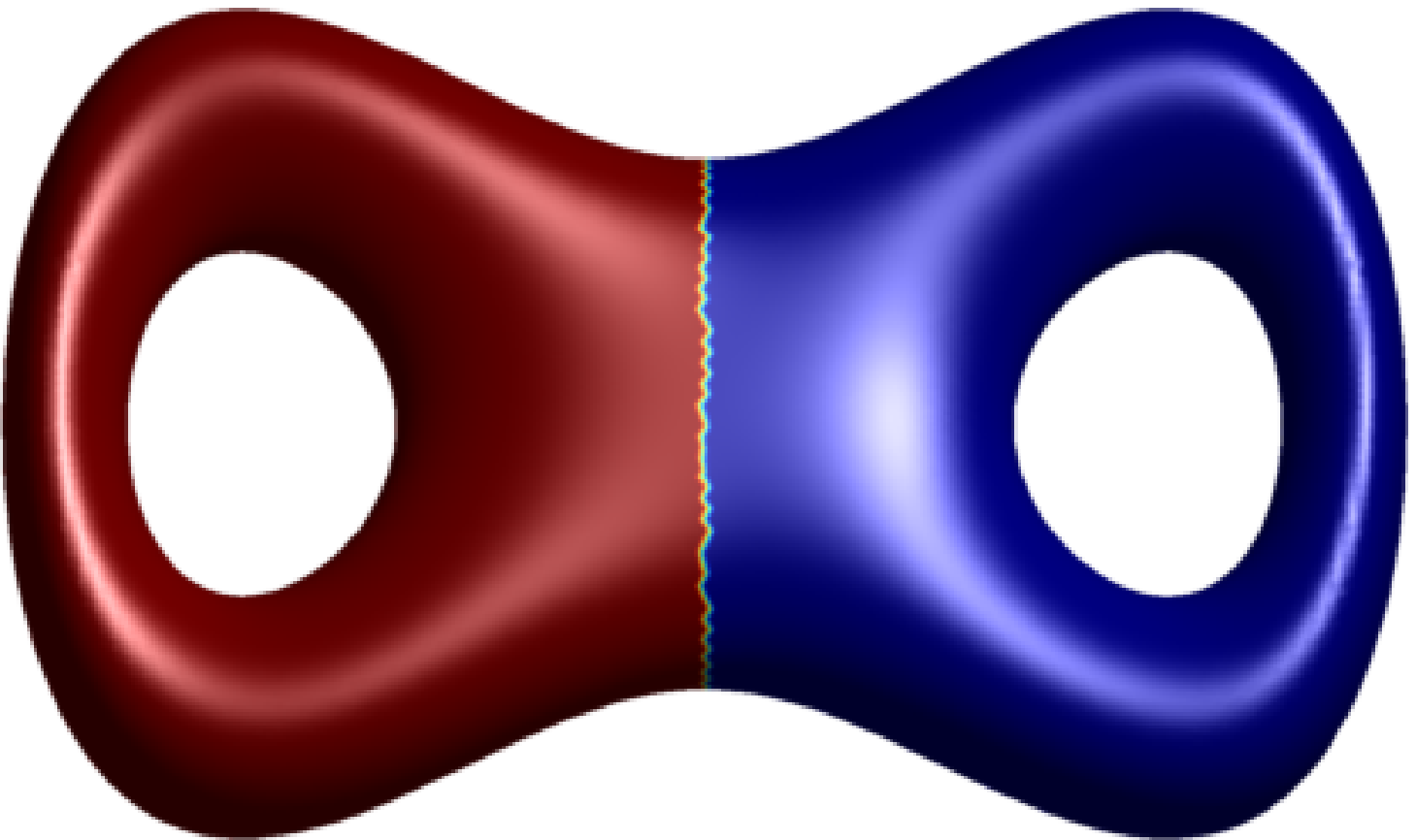}~
\includegraphics[width= 0.29\textwidth]{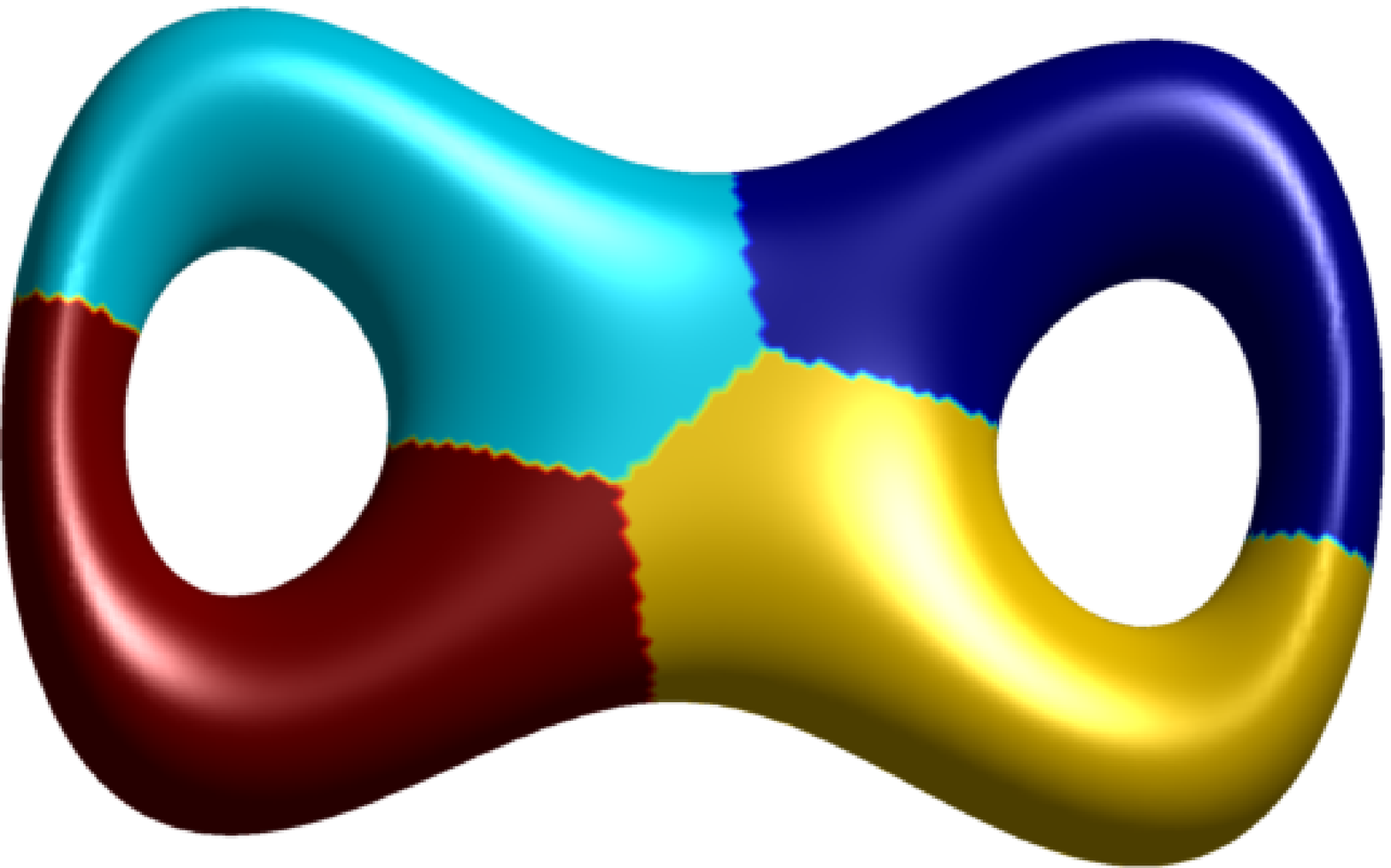}~
\includegraphics[width= 0.29\textwidth]{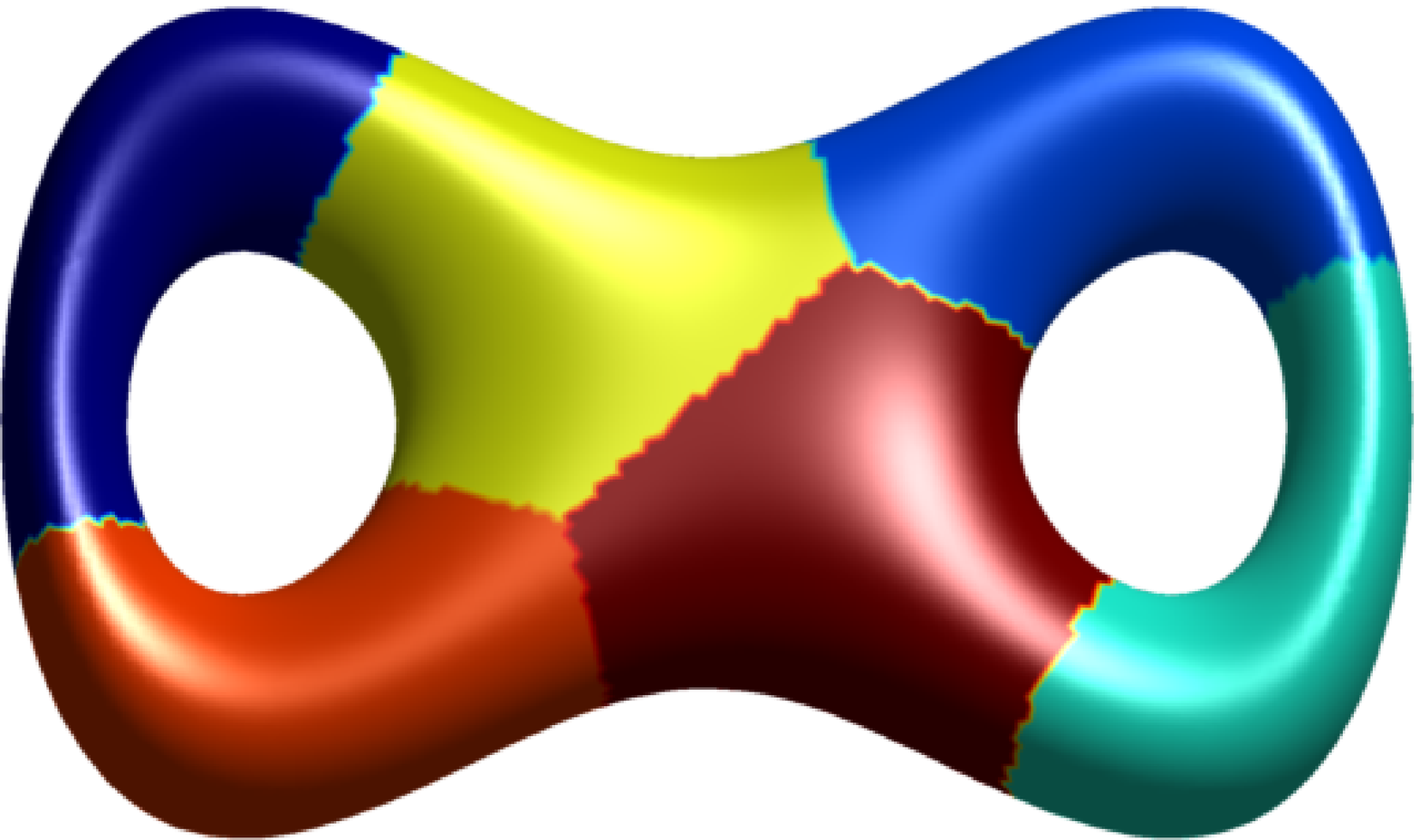}
\caption{Minimal perimeter partitions on a double torus for $n \in \{2,4,6\}$.}
\label{dbtor-perim}
\end{figure}

\begin{figure}
\centering
\includegraphics[width= 0.19\textwidth]{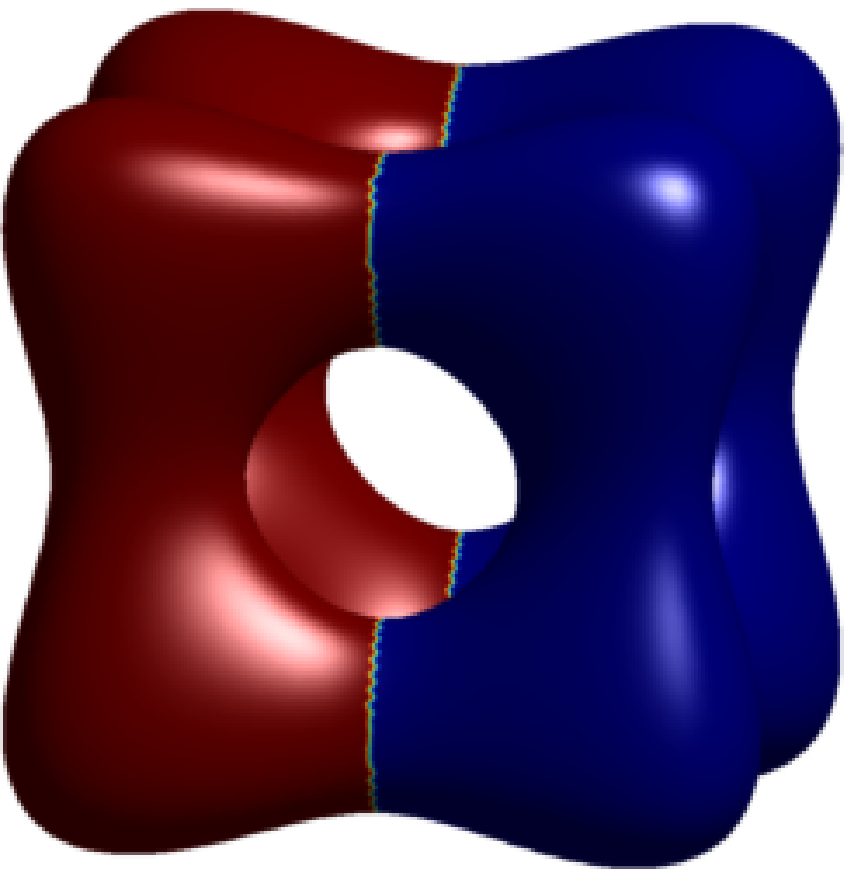}~
\includegraphics[width= 0.19\textwidth]{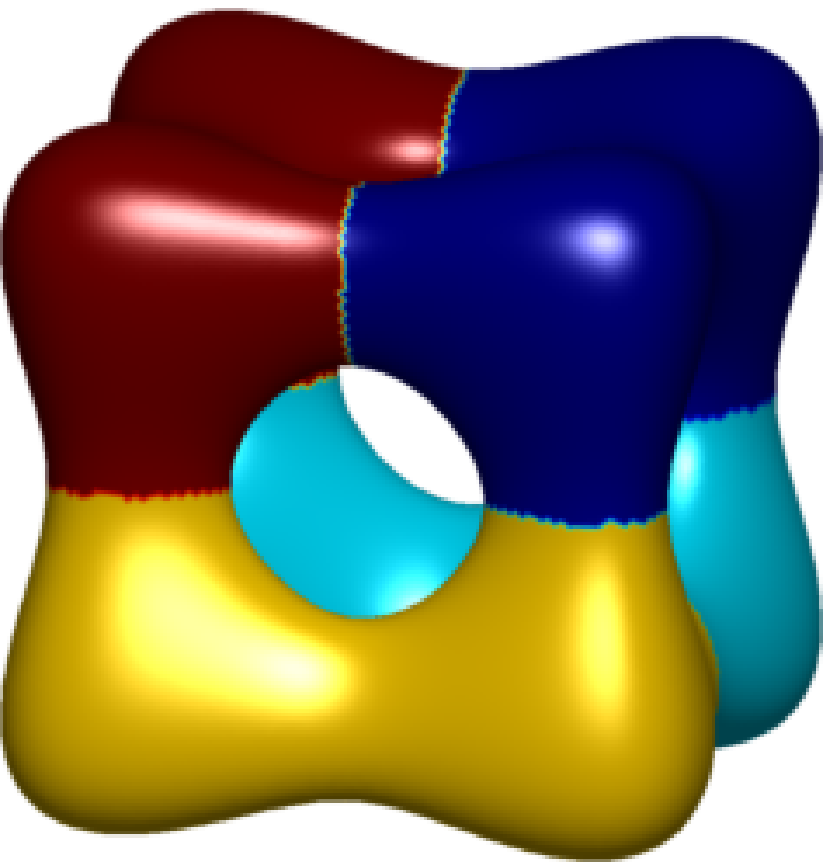}~
\includegraphics[width= 0.19\textwidth]{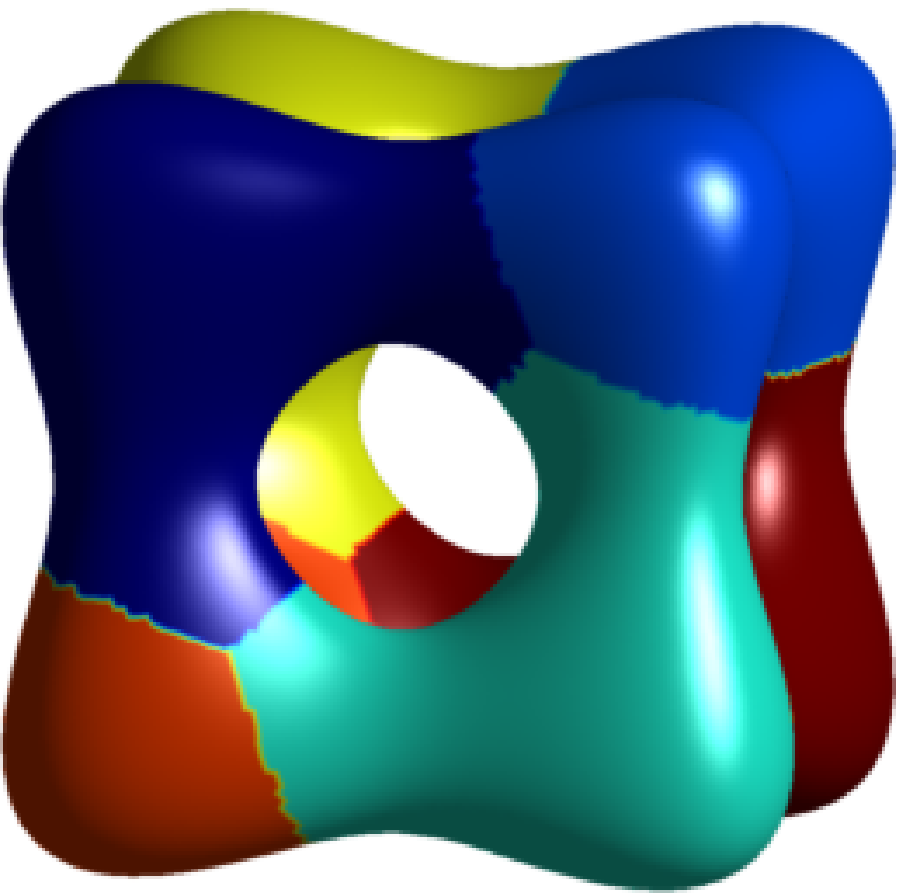}~
\includegraphics[width= 0.19\textwidth]{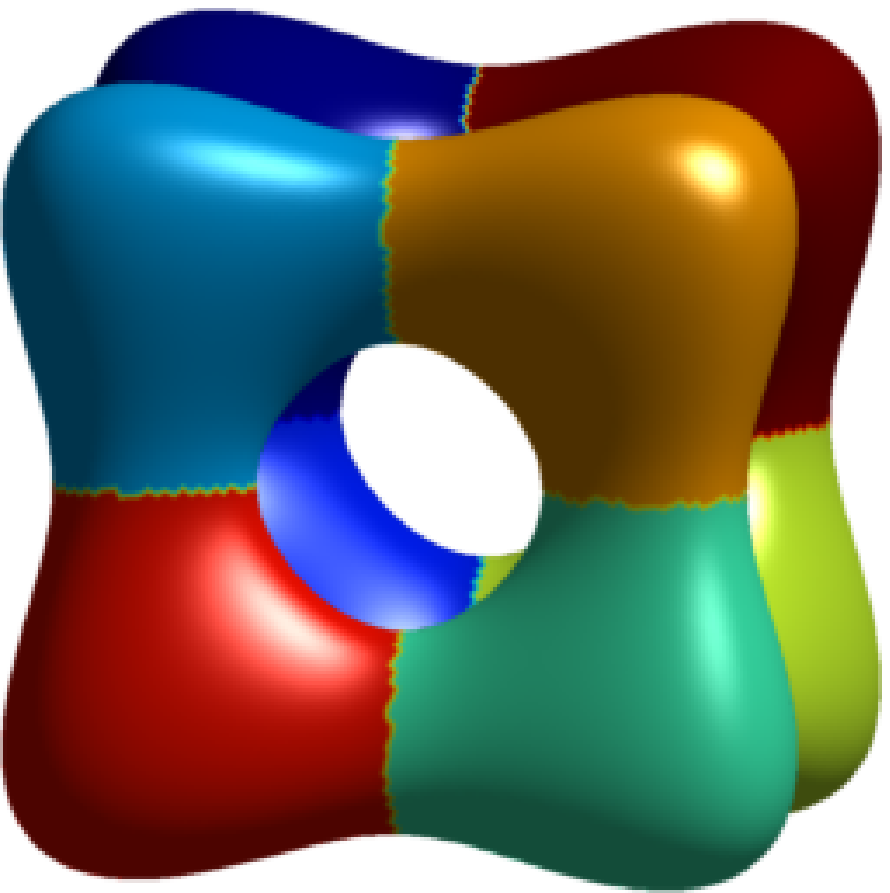}
\caption{Minimal perimeter partitions on a Banchoff-Chmutov surface for $n \in \{2,4,6,8\}$.}
\label{bc-perim}
\end{figure}

As underlined before, our approach allows a direct treatment of any surface, as long as a qualitative triangulation is found. We perform some numerical computations on various shapes like a torus, a double torus, and a more complex surface called Banchoff-Chmutov of order $4$. A few details about the definitions of these surfaces are provided below:
\begin{itemize}
\item We consider a torus of outer radius $R= 1$ and inner radius $0.6$ (see Figure \ref{torus-perim}). This torus is defined as the zero level set of the function
\[ f(x,y,z) = (x^2+y^2+z^2+R^2-r^2)^2 - 4R^2(x^2+y^2).\]
\item The double torus used in the computation (see Figure \ref{dbtor-perim} is given by the zero level set of the function
\[ f(x,y,z) = (x(x-1)^2(x-2)+y^2)^2 +z^2-0.03.\]
\item The complex Banchoff-Chmutov surface (see Figure \ref{bc-perim}) is given by the zero level set of the function
\[ f(x,y,z) = T_4(x)+T_4(y)+T_4(z),\]
where $T_4(X) = 8X^4-8X^2+1$ is the Tchebychev polynomial of order $4$.
\end{itemize}

\section{Refined optimization in the case of the sphere}

The costs associated to the relaxed functional do not provide a good enough approximation of the total length of the boundaries. In this section we propose a method to approximate the optimal cost in the case of the sphere. The results of \cite{morgan-bubbles} state that boundaries of the cells of the optimal partitions have constant geodesic curvature. In the case of the sphere the only such curves are the arcs of circle. See for example \cite[Exercise 2.4.9]{shifrin} for a proof. The results of Cox and Flikkema \cite{cox-partitions} show that optimal configurations are not made of geodesic polygons. In order to perform an optimization procedure which captures this effect they chose to make an initial optimization in the class of geodesic polygons and then divide each geodesic arc into $16$ smaller arcs and restart the procedure with more variable points. They manage to approximate well enough the general optimal structure but they still work in the class of geodesic polygons with additional vertices. Our approach presented below is different in the sense that we consider general circle arcs (not necessarily geodesics) which connect the points. 

The first step is to extract the topology of the partition from the previous density results, i.e. locate the triple points, the edge connections and construct the faces. In order to perform the refined optimization procedure we need to be able to compute the areas of portions of the sphere determined by arcs of circles. This is possible using the Gauss-Bonnet formula. If $M$ is a smooth subset of a surface then
\begin{equation}
\int_M K dA +\int_{\partial M}k_g =2\pi\chi(M),
\label{gauss-bonnet}
\end{equation}
where $K$ is the curvature of the surface, $k_g$ is the geodesic curvature and $\chi(M)$ is the Euler characteristic of $M$. This result extends to piecewise smooth curves and in this case we have
\begin{equation}
\int_M K dA +\int_{\partial M}k_g+\sum \theta_i =2\pi\chi(M),
\label{gauss-bonnet-angles}
\end{equation}
where $\theta_i$ are the \emph{turning angles} between two consecutive smooth parts of the boundary. In the case of a polygon the turning angles are the external angles of the polygon. The formula \eqref{gauss-bonnet-angles} allows the computation of the area of a piece of the sphere bounded by arcs of circle. In this case the Euler characteristic is equal to $1$, the curvature of the unit sphere is $K=1$ and the geodesic curvature is piecewise constant. For more details we refer to \cite[Chapter 4]{docarmo}.

A first consequence of the Gauss-Bonnet theorem in connection to our problem is noting the fact that, apart from cases where we have a certain symmetry like $n \in \{3,4,6,12\}$ the optimal cells are not geodesic polygons. This is made clear in cases where we have a hexagonal cell. If the arcs forming the boundary of such a hexagonal cell would be geodesic polygons then its area would be equal to $6\cdot 2\pi/3-4\pi=0$. Thus a spherical shape bounded by six arcs of circle can never be a geodesic polygon without being degenerate.

In order to perform the optimization we take the vertices as variables and we add one supplementary vertex for each edge. This is enough to contain all the necessary information since an arc of circle is well defined by three distinct points on the sphere. In the sequel we denote $\mathcal P_n$ the set of partitions of the sphere into $n$ cells and with $\mathcal A_n$ the partitions in $\mathcal P_n $ having equal areas. In order to have a simpler numerical treatment of the problem we can incorporate the area constraints in the functional by defining for every partition $(\omega_i) \in \mathcal{P}_n$ the quantity defined for every $\varepsilon>0$ by
\[ G_\varepsilon((\omega_i)) = \sum_{i=1}^n \Per(\omega_i)+\frac{1}{\varepsilon} \sum_{i=1}^{n-1} \sum_{j=i+1}^n (\text{Area}(\omega_i)-\text{Area}(\omega_j))^2.\] 
If we denote 
\[ G((\omega_i)) = \begin{cases} \sum_{i=1}^n \Per(\omega_i) & \text{ if }(\omega_i) \in \mathcal{A}_n \\
 \infty & \text{ if } (\omega_i) \in \mathcal{P}_n \setminus \mathcal{A}_n.
 \end{cases}\]
 then we have the following $\Gamma$-convergence result.
 
 \begin{thm}
 We have that $G_\varepsilon \gconv G$ for the $L^1(\Bbb{S}^2)$ convergence of sets.
 \label{gconv-refined}
 \end{thm}
 
 \emph{Proof:} For the (LI) property consider a sequence $(\omega_i^\varepsilon) \subset \mathcal P_n$ which convergence in $L^1(\Bbb{S}^2)$ to $(\omega_i)$. It is clear that we have $\text{Area}(\omega_i^\varepsilon) \to \text{Area}(\omega_i)$ and the perimeter is lower semicontinuous for the $L^1$ convergence. Thus we have two situations. If $(\omega_i) \in \mathcal P_n \setminus \mathcal A_n$ then $\lim_{\varepsilon \to 0}G_\varepsilon((u_i^\varepsilon)) = \infty$. If $(\omega_i) \in \mathcal{A}_n$ then the lower semicontinuity of the perimeter implies that $\liminf_{\varepsilon \to 0} G_\varepsilon((\omega_i^\varepsilon)) \geq G((\omega_i))$.

 The (LS) property is immediate in this case. Choose $(\omega_i) \in \mathcal A_n$, or else there is nothing to prove. We may choose the recovery sequence equal to $(\omega_i)$ for every $\varepsilon>0$. Thus the property is verified immediately. \hfill $\square$
 
 \begin{rem}
 We note that in the above proof the simplicity of the proof of the (LS) property is due to the fact that the functionals $G_\varepsilon$ are well defined on the space $\{G<\infty\}$, which makes possible the choice of constant recovery sequences. This is not the case in the results proved in Section \ref{theoretical-result}.
 \end{rem}
 
 This $\Gamma$-convergence result proves that minimizers of $G_\varepsilon$ converge to minimizers of $G$. As a consequence, in the numerical computations, we minimize $G_\varepsilon$ for $\varepsilon$ smaller and smaller in order to approach the minimizers of $G$, which are in fact the desired solutions to our problem. 
 
 Since the parameters are of two types: triple points and edge points, we prefer to use an optimization algorithm which is not based on the gradient. The algorithm is described below.
 \begin{itemize}
 \item For each point $P$ consider a family of $m$ tangential directions $(v_i)_{i=1}^m$ chosen as follows: the first direction is chosen randomly and the rest are chosen so that the angles between consecutive directions are $2\pi/m$.
 \item Evaluate the cost function for the new partition obtained by perturbing the point $P$ in each of the directions $v_i$ according to a parameter $\varepsilon$.
 \item Choose the direction which has the largest decrease and update the partition accordingly.
 \item Do the same procedure for each edge point by performing the two possible orthogonal perturbations of the point with respect to the edge.
 \item If there is no decrease for each of the points of the partition, then decrease $\varepsilon$.
 \end{itemize}
 
\begin{figure}
  \centering
  \includegraphics[width = 0.3\textwidth]{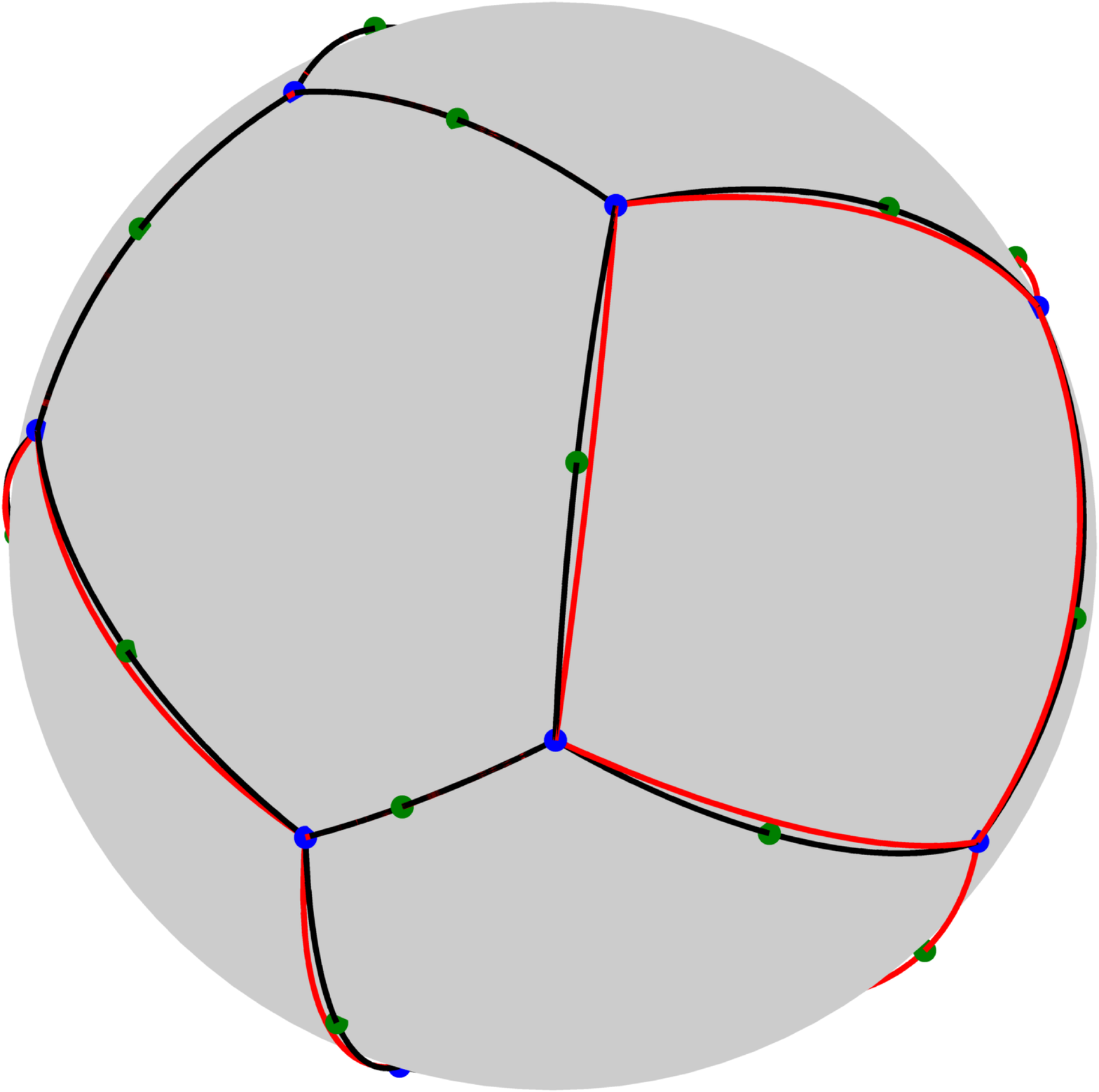}~
  \includegraphics[width = 0.3\textwidth]{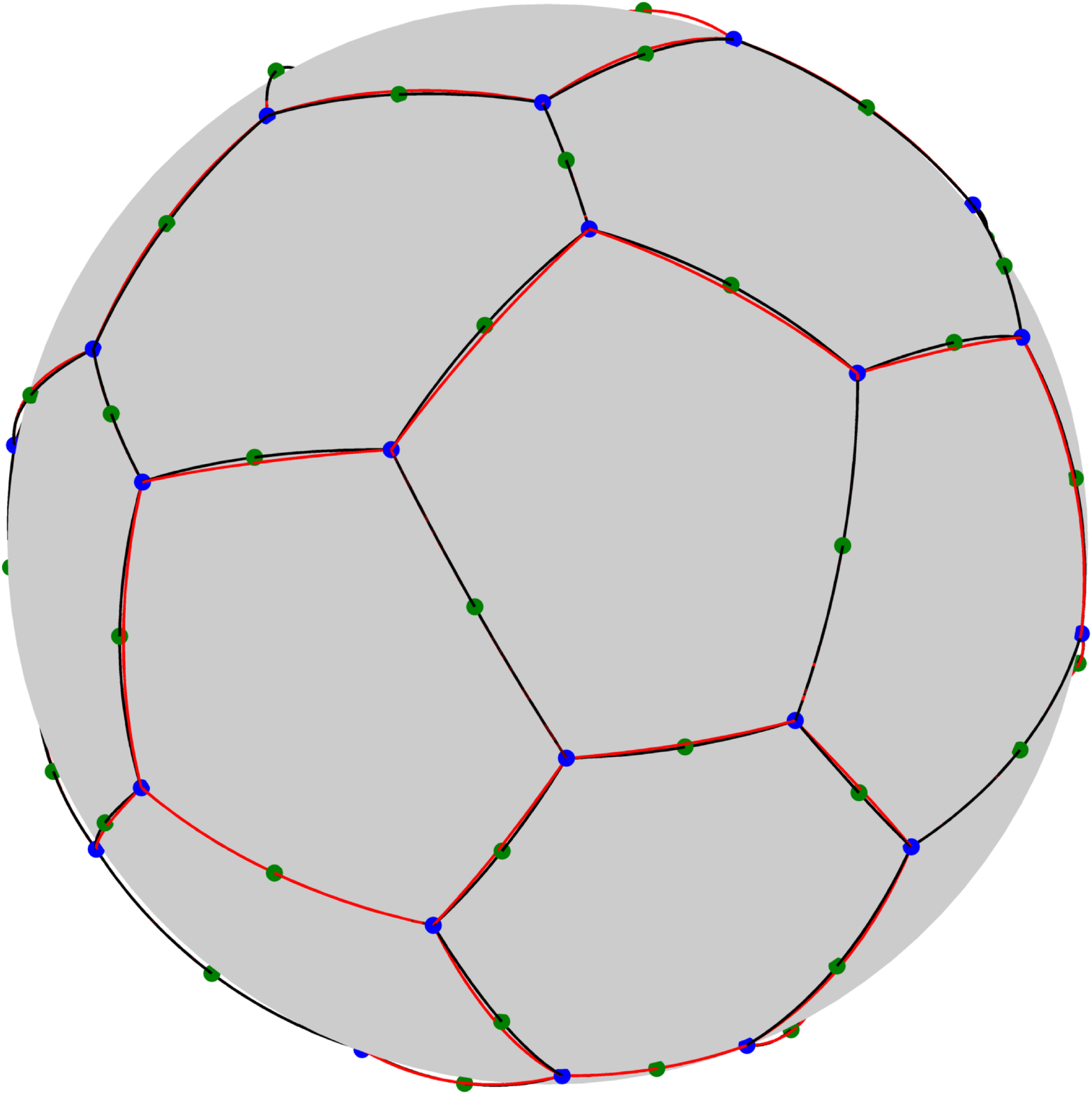}
  \caption{The difference between optimal configuration (black) and the  geodesics connecting the points (red).}
  \label{refined_tests} 
  \end{figure}

 This algorithm converges in each of the test cases and the results are presented in Table \ref{comparison-cox}. In the optimization procedure we start with $\varepsilon=1$ and we reiterate the optimization decreasing $\varepsilon$ by a factor of $10$ at each step until we reach the desired precision on the area constraints. We are able to recover the same results as Cox and Flikkema for $n\in [4,32]$. Furthermore, unlike in the case of geodesic polygons, all triple points consist of boundaries which meet at equal angles of measure $2\pi/3$. In Figure \ref{refined_tests} you can see the results for $n=9$ and $n=20$. The red arcs are geodesic connecting the points and are drawn to visually see that not all the boundaries of the optimal structure are geodesic arcs.
 
 \begin{table}[!htp]
 \centering
 \begin{tabular}[t]{|c||c|c||c|}
 \hline
   & \multicolumn{2}{c||}{our results} & Cox-Flikkema \\
  \hline
 $N$ & non-geo. & area tol. & non-geo. \\
 \hline
 $4$  & $11.4637$ & $5\e{-7}$ & $11.464$ \\
 \hline 
 $5$  & $13.4304$ & $2\e{-7}$ & $13.430$ \\
 \hline
 $6$  & $14.7715$ & $2\e{-7}$ &  $14.772$ \\
 \hline
 $7$  &  $16.3519$ & $3\e{-7}$ & $16.352$ \\
 \hline
 $8$  & $17.6927$  & $3\e{-7}$   &   $17.692$ \\
 \hline
 $9$  & $18.8504$  & $2\e{-7}$  &   $18.850$ \\
 \hline
 $10$ & $19.9997$   & $4\e{-7}$ &   $20.000$ \\
 \hline
 $11$ & $21.1398$  & $4\e{-7}$  &   $21.140$ \\
 \hline
 $12$ & $21.8918$ & $5\e{-7}$  &   $21.892$ \\
 \hline
 $13$ & $23.0953$  & $4\e{-7}$  &   $23.095$ \\
 \hline
 $14$ & $23.9581$  & $3\e{-7}$  &   $23.958$ \\
 \hline
 $15$ & $24.8821$  & $2\e{-7}$  &   $24.882$ \\
 \hline
 
 $16$ & $25.7269$   & $2\e{-7}$  &   $25.727$ \\
 \hline
 
 $17$ & $26.6365$   & $3\e{-7}$  &   $26.637$ \\
 \hline
 
 $18$ & $27.4647$  & $2\e{-7}$  &   $27.465$ \\
 \hline
 \end{tabular}
 \begin{tabular}[t]{|c||c|c||c|}
  \hline
    & \multicolumn{2}{c||}{our results} & Cox-Flikkema \\
   \hline
 $N$ & non-geo. & area tol. & non-geo. \\
 \hline
 $19$ & $28.2735$  & $2\e{-7}$  &   $28.274$ \\
 \hline
 
 $20$ & $28.9992$  & $1\e{-7}$  &   $28.999$ \\
 \hline
 
 $21$ & $29.7748$  & $2\e{-7}$  &   $29.775$ \\
 \hline
 
 $22$ & $30.5094$  & $2\e{-7}$  &   $30.509$ \\
 \hline
 
 $23$ & $31.2260$  & $2\e{-7}$  &   $31.226$ \\
 \hline
 $24$ & $31.9117$  & $3\e{-7}$  &   $31.912$ \\
 \hline
 $25$ & $32.6172$  & $8\e{-8}$  &   $32.617$ \\
 \hline
 $26$ & $33.2675$  & $2\e{-7}$  &   $33.268$ \\
 \hline
 $27$ & $33.8968$  & $9\e{-8}$  &   $33.897$ \\
 \hline
 $28$ & $34.5521$  & $4\e{-7}$  &   $34.552$ \\
 \hline
 $29$ & $35.2065$  & $6\e{-7}$  &   $35.207$ \\
 \hline
 $30$ & $35.8199$  & $5\e{-7}$  &   $35.820$ \\
 \hline
 $31$ & $36.3941$  & $4\e{-6}$  &   $36.394$ \\
 \hline
 $32$ & $36.9310$  & $4\e{-6}$  &   $36.931$ \\
 \hline

 \end{tabular}
 \caption{Comparison between our results and the results of Cox and Flikkema in the case of the sphere.}
 \label{comparison-cox}
 \end{table}
 
 Thus we can conclude that the relaxed formulation presented in the previous section is able to match the best known configurations in the literature. Furthermore for $n \in [5,25]\cup \{32\}$ the algorithm finds the good configuration without much effort, while for $n \in [26,31]$ multiple tries with different initial conditions were needed in order to find the best configuration. The fact that the structure of the partition is not fixed is a great advantage offered by our method.
 
\section{Computing the optimal cost - general surfaces}

The approach used in the previous section cannot be applied to other surfaces than the sphere. Indeed, the general expression of curves of constant curvature is not known explicitly for other types of surfaces. One way to approximate the total perimeter of the partition would be to extract the contours of the optimal densities and evaluate the length of each discrete contour. A natural way to extract a contour corresponding to a density function would be taking a level set, for example the level $0.5$. It is possible to extract such level sets by looking at which triangles contain values which are both above and below the level set. On each triangle which is cut by the contour we make a linear interpolation which determines a segment in the contour of the level set.

\begin{figure}
\centering
v\includegraphics[width = 0.3\textwidth]{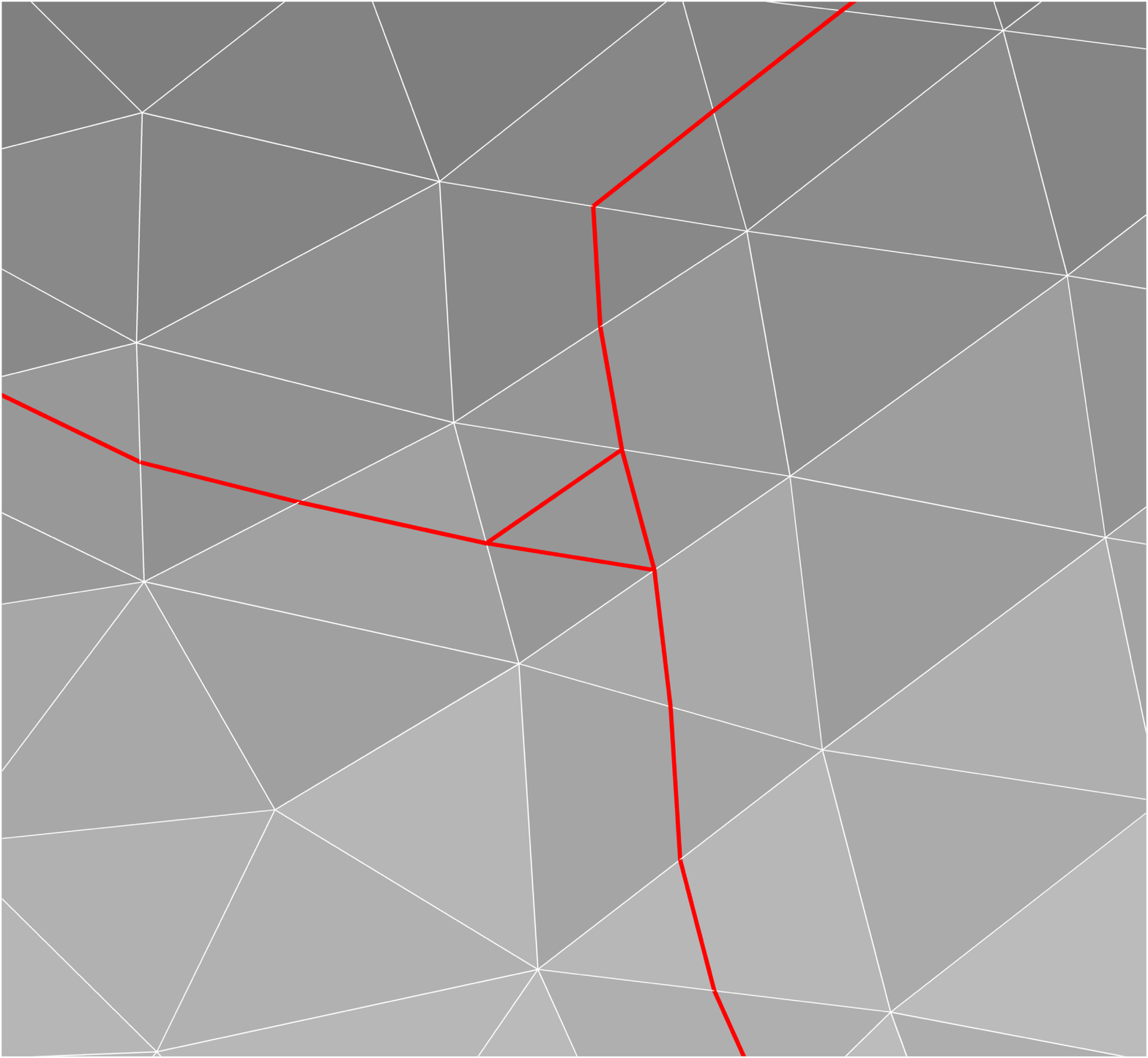}~
\includegraphics[width = 0.3\textwidth]{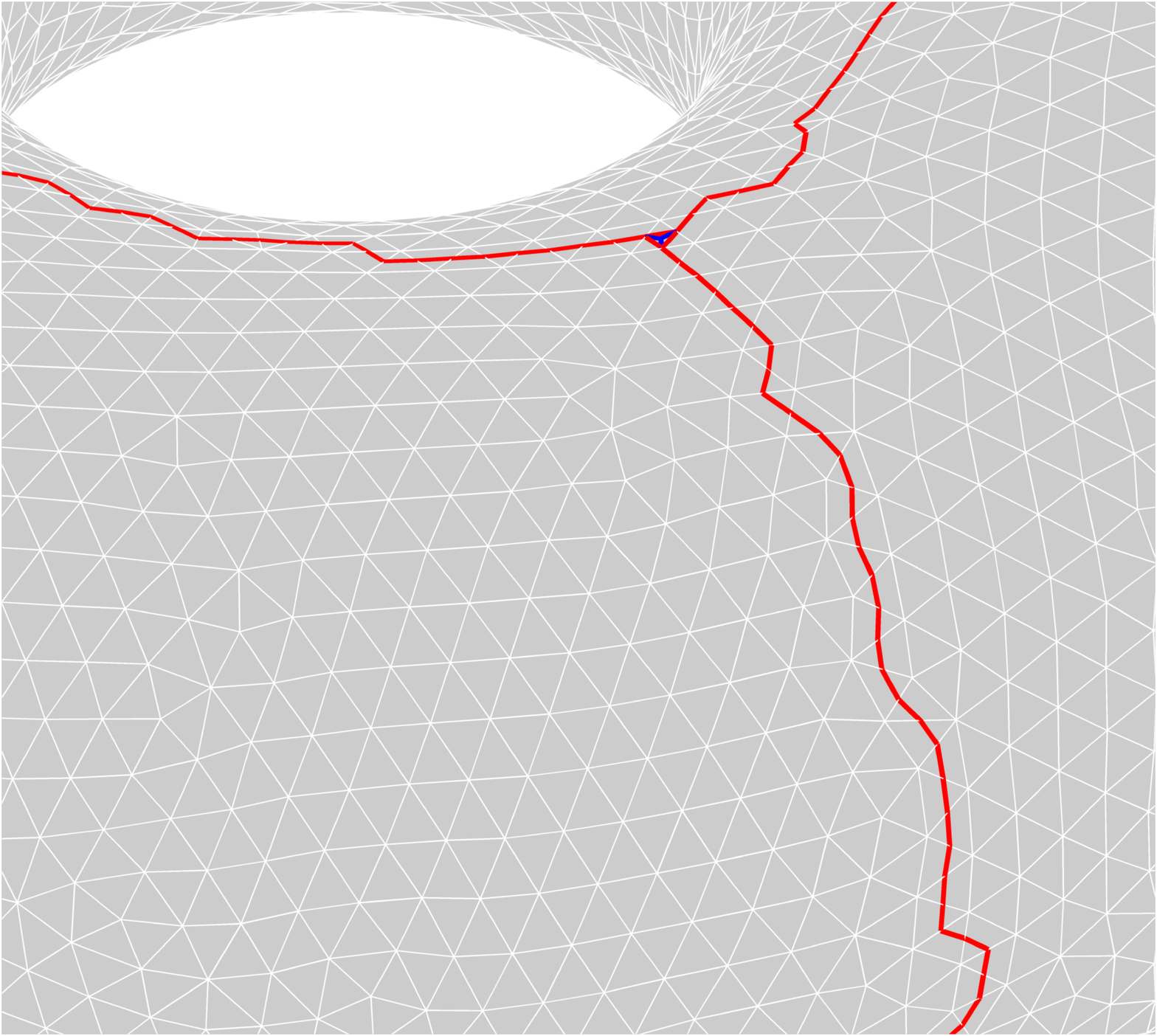}
\caption{A small space left around triple points (left) and the non-regular initial extracted contours (right).}
\label{issues}
\end{figure}

Once we have an idea on how to extract the contours, the first question arises: how to make sure that the level sets extracted form a partition of $S$? We denote by $\mathcal T$ a triangulation of $S$. If we think of extracting the $0.5$ levels of each density, the shapes determined by these contours will not overlap, but around triple points there will be some free space left. One way to make sure that we have extracted a partition is to take the $0.5$ levels of the function defined on the triangulation $\mathcal T$ by
\begin{equation} \phi_i(x) = \begin{cases}
1 & \text{ if } u_i(x)\geq \max_{i\neq j} u_j(x)\\
0 & \text{ otherwise},
\end{cases}
\label{maxphi}
\end{equation}
where $u_i$ are the optimal densities obtained numerically. These contour levels of the functions $\phi_i$ almost realize a partition of $S$ with the following issues:
\begin{enumerate}
\item There is a small void space around each triple point, but this void is included in one of the triangles of the mesh, and can be dealt with.
\item Since we extract the level sets of a function which is either $0$ or $1$ on the vertices of the triangulation, the contour lines will pass through the middle of the edges of the triangles situated at the border between two phases. This creates some contours which are quite zigzagged and whose length is significantly larger than the optimal total perimeter.
\end{enumerate}
We illustrate these two issues in Figure \ref{issues}. 

Nevertheless, once we have extracted these contours it is possible to make a direct optimization of the total length of the boundaries with the constraint of fixed area of the cells. This optimization is made directly on the triangulated surface. We describe the optimization algorithm below.

{\bf Variables and representation of the partitions.} We denote $(x_i)_{i=1}^h$ a generic family of variable points situated each on an edge of the triangulation $\mathcal{T}$ such that each edge contains exactly one variable point. To these points we associate a family of parameters $(\lambda_i)_{i=1}^h$ which gives the position of each point $x_i$ on the corresponding edges. We take this global parametric approach since each of these points belongs to at least two cells and we'll need to evaluate its contribution in the gradient of the area and the for all the cells that contain it. Having a global sets of points avoids having to match points between different contours.

Each cell of the partitions is represented by a structure of pairs of edges of triangles of $\mathcal{T}$ which determine, along with the parameters $(\lambda_i)$, the segments which form the discrete contour of the cell. The pairs of edges is ordered so that the contour is continuous. Contours may have one or more connected components.

{\bf Computation of the perimeters of the cells.} The perimeter of a cell is computed by following the segments forming the contour and incrementally adding their lengths to the total length. If the vertices of the segment are given by $x_i = \lambda_i v_1+(1-\lambda_i)v_2$ and $x_j = \lambda_j v_3+(1-\lambda_j)v_4$ then the length of the segment $[x_i,x_j]$ is
\[ \ell([x_i,x_j]) = \|\lambda_i v_1+(1-\lambda_i)v_2- \lambda_j v_3-(1-\lambda_j)v_4\|,\]
expression which is differentiable if the length is not zero. The derivatives with respect to $\lambda_i$ and $\lambda_j$ are then added to the gradient vector. Note that for the points which are not vertices of some contour the gradient is zero.

{\bf Computation of the areas of the cells.} In order to compute the area of a cell we use the information given by the functions $\phi_i$ defined in \eqref{maxphi}. The function $\phi_i$ shows, among other things, what is the position of each triangle in $\mathcal T$ with respect to the cell $i$. Indeed, denoting by $T$ a triangle in $\mathcal{T}$, we have the following cases:
\begin{enumerate}
\item All the vertices $v$ of the triangle $T$ satisfy $\phi_i(v) = 1$. Then $T$ is completely inside the cell $i$ and we add its area to the total area of the cell.
\item Two vertices $v_1,v_2$ of $T$ satisfy $\phi_i(v_{1,2}) = 1$ and the third satisfies $\phi_i(v_3)=0$. Thus we only add a portion of the area of $T$ to the total area of cell $i$. Note that this value of the area depends linearly of one parameter $\lambda_k$ and of another parameter $\lambda_l$. The derivatives of these contributions are added to the vectors containing the gradient of the area of the cell $i$.
\item Two vertices $v_1,v_2$ of $T$ satisfy $\phi_i(v_{1,2}) = 0$ and the third satisfies $\phi_i(v_3)=1$. Again, we only add a portion of the area of $T$ to the total area of cell $i$ which again depends linearly of one parameter $\lambda_k$ and of another parameter $\lambda_l$. The derivatives of these contributions are added to the vectors containing the gradient of the area of the cell $i$.
\item If all the vertices of $T$ satisfy $\phi_i(v) = 0$ then the triangle is outside the cell and we move on.
\end{enumerate}

\begin{figure}
\includegraphics[width=0.5\textwidth]{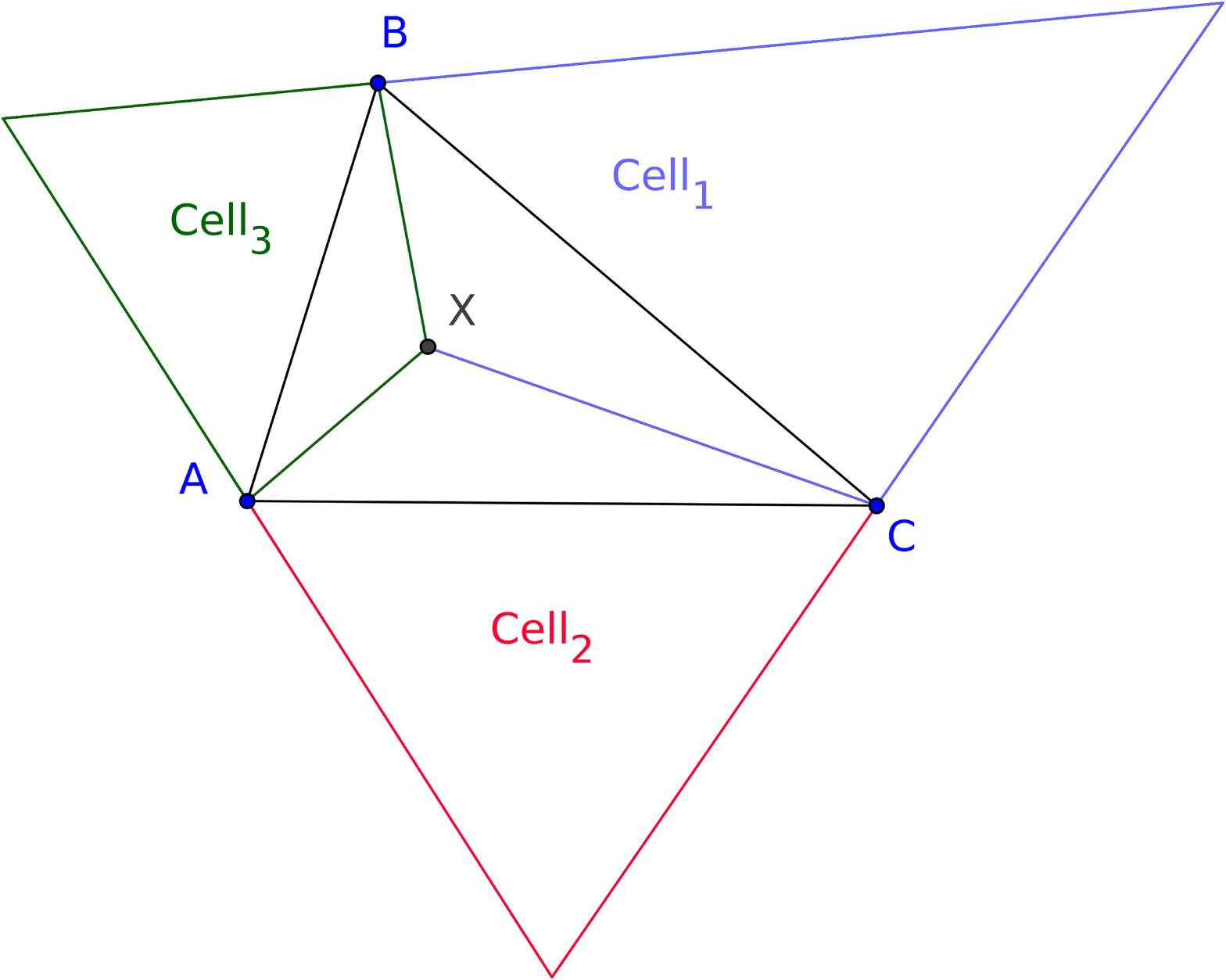}
\caption{Treatment of empty space around triple points. We consider the Fermat point $X$ of the empty triangle $ABC$ and we add corresponding area and perimeters to the corresponding cells. For example the area of $ABX$ is added to Cell 3 and the quantity $AX+BX-AB$ is added to the perimeter of Cell 3.}
\label{steiner}
\end{figure}

{\bf The empty spaces around triple points.} As we have noted above and seen in Figure \ref{issues}, around triple points we have some empty spaces determined by three points which belong to the three sides of some of the triangles in $\mathcal T$. In each configuration of this type we add a Steiner tree corresponding to the three variable points. Each of the three area regions which are formed are added to the corresponding cell while the perimeter is modified with the length of two adjacent segments in the Steiner tree. See Figure \ref{steiner} for further details. In order to find the gradient corresponding to the lengths and area changes due to the addition of these Steiner points we use a finite differences approximation. 

\begin{figure}
\includegraphics[width = 0.3\textwidth]{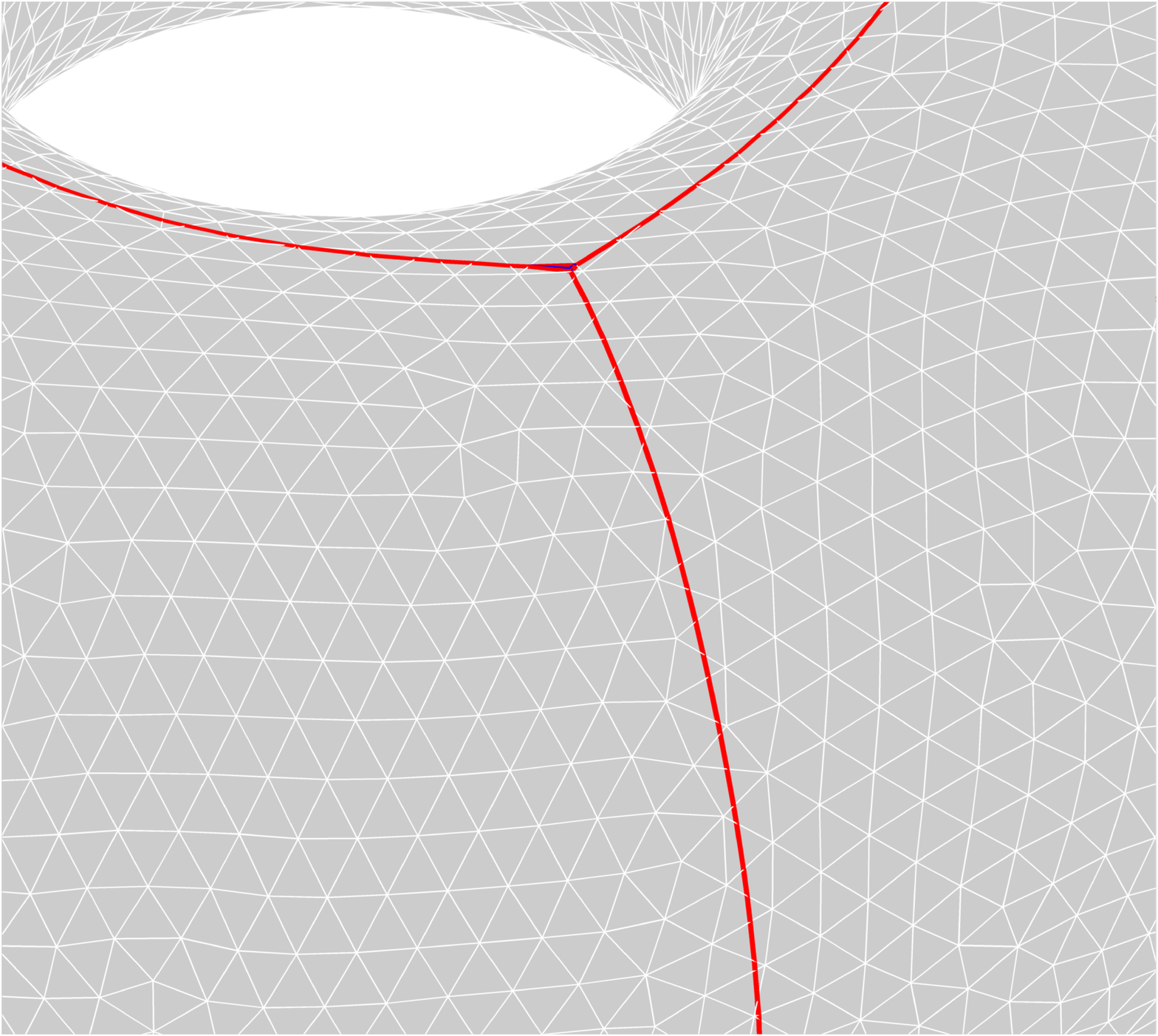}~
\includegraphics[width = 0.3\textwidth]{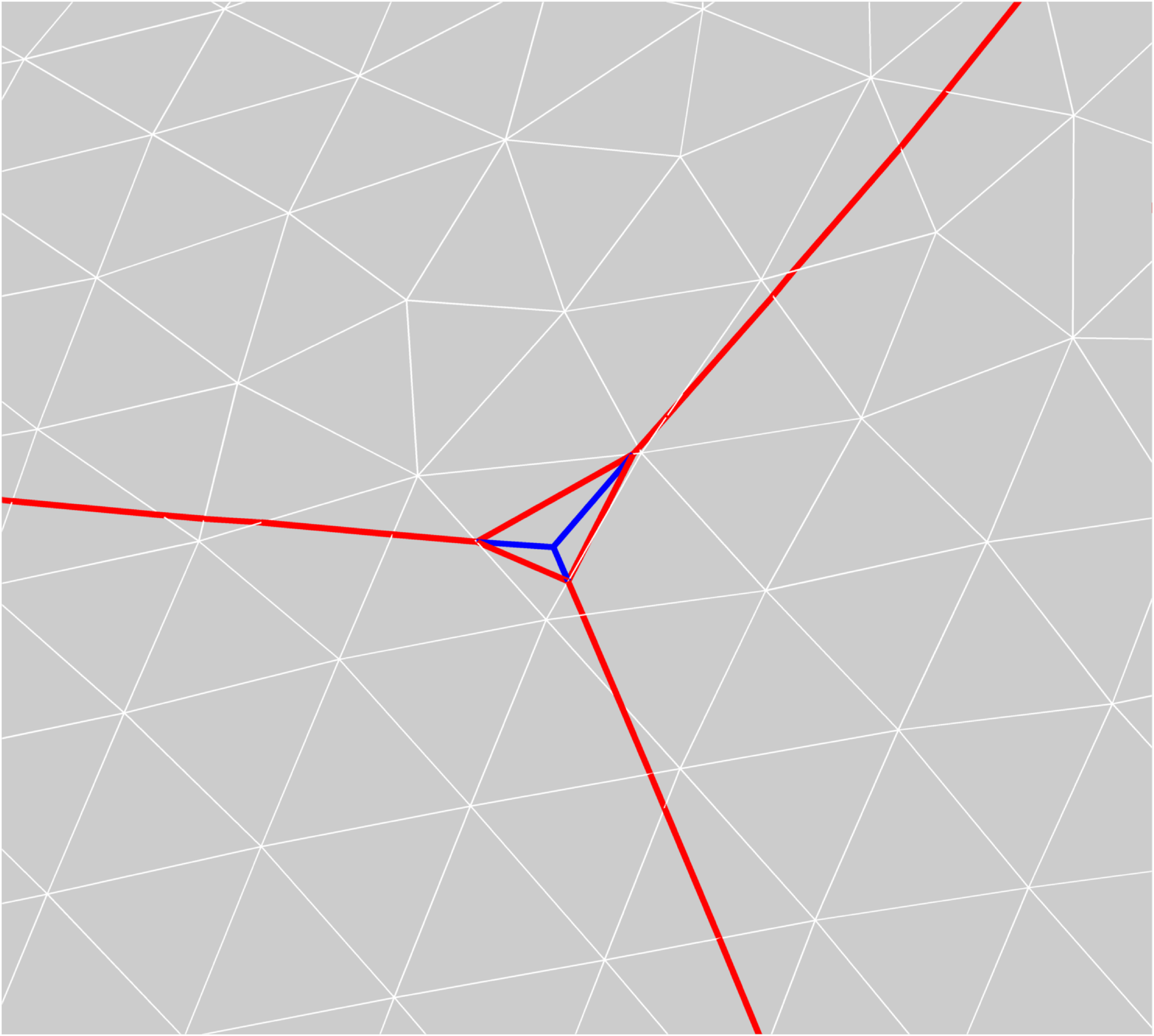}
\caption{Contours after the constrained optimization algorithm. You can also see a zoom around the triple points: the segments which join the Fermat points align themselves with the rest of the contour.}
\label{smooth}
\end{figure}

{\bf Constrained optimization algorithm.} We have the expressions and the gradients of the perimeters and areas of the cells as functions of the parameters $(\lambda_i)_{i=1}^h$. This allows us to use the algorithm \texttt{fmincon} from the Matlab Optimization Toolbox in order to implement the constrained optimization algorithm. We use the \emph{interior-point} algorithm with a low-memory hessian approximation given by an \emph{LBFGS} algorithm. The initial values of the parameters $(\lambda_i)_{i=1}^h$ are all set to $0.5$. The algorithm manages to satisfy the constraints at machine precision while minimizing the perimeter and thus smoothing the zigzagged initial contours (like the ones in Figure \ref{issues}). An example of result may be seen in Figure \ref{smooth}.

It may be the case that some vertices of the contour would "like" to switch to another side. This can be the case if at the end of the optimization one of the parameters $\lambda_i$ is close to $0$ or $1$ or a triple point in one of the constructed Steiner trees is on the boundary of the corresponding mesh triangle. In this cases we modify the initial contours taking into the account these results and we restart the optimization procedure. The modification is done in the following way.
\begin{enumerate}
\item If one of the $\lambda_i$ is equal to $0$ or $1$ then we add the corresponding point to the adjacent cell and restart the algorithm.
\item If one of the triple points arrives on the edge of its corresponding mesh triangle then we allow it to move to the adjacent triangle.
\end{enumerate}
After a finite number of switches the configuration stabilizes and a local minimum is found. 

We test the presented algorithm on the results obtained in previous sections. In the case of the sphere we obtain the same values found in Table \ref{comparison-cox}. The approximations of the optimal costs for partitions presented in Figure \ref{torus-perim} for a torus of radii $R=1, r=0.6$ in Table \ref{tor-values}.

\begin{center}
\begin{table}
\begin{tabular}[t]{|c|c|}
\hline
$n$ & Minimal length  \\
\hline\hline
$2$ &  $15.07$     \\
\hline
$3$ &  $22.61$       \\
\hline
$4$ &  $30.15$        \\
\hline
$5$ & $37.25$         \\
\hline
$6$ &  $41.93$     \\
\hline
\end{tabular}
\begin{tabular}[t]{|c|c|}
\hline
$n$ & Minimal length  \\
\hline\hline
 $7$    & $47.12$      \\
\hline
 $8$    & $50.77$      \\
\hline
 $9$    & $53.37$      \\
\hline
 $10$    & $56.80$      \\
\hline
\end{tabular}
\caption{Approximation of the optimal costs for minimal partitions of a torus into equal area cells. These partitions are represented in Figure \ref{torus-perim}}
\label{tor-values}
\end{table}
\end{center}

\section{Conclusions}
We propose an algorithm for finding numerically the partitions which divide a surface into cells of prescribed areas and minimize the sum of the corresponding perimeters. This algorithm is rigorously justified by a $\Gamma$-convergence result which is a generalization of the Modica-Mortola theorem in the case of smooth $(d-1)$-dimensional manifolds.

In the case of the sphere we are able to recover all the results presented in the article of Cox and Flikkema \cite{cox-partitions}. The optimal costs of the spherical partitions are precisely evaluated by using the qualitative results in \cite{morgan-bubbles}, which imply that the boundaries of the cells are arcs of circles. We recover the same optimal costs as the ones presented in \cite{cox-partitions}. We underline that one of the advantages of this relaxed method is the fact that we do not need to set the polyhedral configuration of the partition \emph{a priori}. The cells emerge from random density configurations and place themselves in the best positions.

The $\Gamma$-convergence method is not limited to the case of the sphere. Once we have triangulated a surface the same algorithm applies. We present a few test cases of more complex surfaces. While the relaxed optimal partitions can easily be obtained, computing the optimal costs is not straightforward since the relaxed costs are not precise enough. In order to be able to compute an approximation of these optimal costs we extract the contours of the optimal densities and we perform a constrained optimization on the triangulated surface.

\bibliography{./master}

\begin{thebibliography}{10}

\bibitem{gammaconvalberti}
Giovanni Alberti.
\newblock Variational models for phase transitions, an approach via
  gamma-convergence.
\newblock 1998.

\bibitem{ambrosio}
Luigi Ambrosio and Andrea Braides.
\newblock Functionals defined on partitions in sets of finite perimeter. {II}.\
  {S}emicontinuity, relaxation and homogenization.
\newblock {\em J. Math. Pures Appl. (9)}, 69(3):307--333, 1990.

\bibitem{ambrosiofuscopallara}
Luigi Ambrosio, Nicola Fusco, and Diego Pallara.
\newblock {\em Functions of bounded variation and free discontinuity problems}.
\newblock Oxford Mathematical Monographs. The Clarendon Press, Oxford
  University Press, New York, 2000.

\bibitem{baldo-manifolds}
Sisto Baldo and Giandomenico Orlandi.
\newblock Cycles of least mass in a {R}iemannian manifold, described through
  the ``phase transition'' energy of the sections of a line bundle.
\newblock {\em Math. Z.}, 225(4):639--655, 1997.

\bibitem{bernstein-sphere}
Felix Bernstein.
\newblock \"{U}ber die isoperimetrische {E}igenschaft des {K}reises auf der
  {K}ugeloberfl\"ache und in der {E}bene.
\newblock {\em Math. Ann.}, 60(1):117--136, 1905.

\bibitem{braides2}
Andrea Braides.
\newblock {\em Approximation of {F}ree-{D}iscontinuity Problems}.
\newblock Springer, 1998.

\bibitem{evolver}
Kenneth~A. Brakke.
\newblock The surface evolver.
\newblock {\em Experiment. Math.}, 1(2):141--165, 1992.

\bibitem{buttazzogconv}
Giuseppe Buttazzo.
\newblock Gamma-convergence and its {A}pplications to {S}ome {P}roblems in the
  {C}alculus of {V}ariations.
\newblock {\em School on {H}omogenization {ICTP}, {T}rieste, {S}eptember 6-17},
  1993.

\bibitem{cox-partitions}
S.~J. Cox and E.~Flikkema.
\newblock The minimal perimeter for {$N$} confined deformable bubbles of equal
  area.
\newblock {\em Electron. J. Combin.}, 17(1):Research Paper 45, 23, 2010.

\bibitem{docarmo}
Manfredo~P. do~Carmo.
\newblock {\em Differential geometry of curves and surfaces}.
\newblock Prentice-Hall, Inc., Englewood Cliffs, N.J., 1976.
\newblock Translated from the Portuguese.

\bibitem{engelstein-four}
Max Engelstein.
\newblock The least-perimeter partition of a sphere into four equal areas.
\newblock {\em Discrete Comput. Geom.}, 44(3):645--653, 2010.

\bibitem{hales}
Thomas~C. Hales.
\newblock The honeycomb conjecture.
\newblock {\em Discrete \& Computational Geometry}, 25(1):1--22, 2001.

\bibitem{hales-sphere}
Thomas~C. Hales.
\newblock The honeycomb problem on the sphere, 2002.

\bibitem{henrot-pierre}
Antoine Henrot and Michel Pierre.
\newblock {\em Variation et optimisation de formes}, volume~48 of {\em
  Math\'ematiques \& Applications (Berlin) [Mathematics \& Applications]}.
\newblock Springer, Berlin, 2005.
\newblock Une analyse g{\'e}om{\'e}trique. [A geometric analysis].

\bibitem{masters-sphere}
Joseph~D. Masters.
\newblock The perimeter-minimizing enclosure of two areas in {$S^2$}.
\newblock {\em Real Anal. Exchange}, 22(2):645--654, 1996/97.

\bibitem{morgan-bubbles}
Frank Morgan.
\newblock Soap bubbles in {${\bf R}^2$} and in surfaces.
\newblock {\em Pacific J. Math.}, 165(2):347--361, 1994.

\bibitem{oudet}
{\'E}douard Oudet.
\newblock Approximation of partitions of least perimeter by
  {$\Gamma$}-convergence: around {K}elvin's conjecture.
\newblock {\em Exp. Math.}, 20(3):260--270, 2011.

\bibitem{shifrin}
Theodore Shifrin.
\newblock Differential geometry - a first course in curves and surfaces.

\bibitem{lbfgs}
Liam Stewart.
\newblock Matlab lbfgs wrapper.
\newblock http://www.cs.toronto.edu/~liam/software.shtml.

\end{thebibliography}
\bibliographystyle{plain}

\end{document}